\documentclass[twoside,11pt]{article}

\usepackage{blindtext}

\usepackage[abbrvbib, preprint]{jmlr2e}

\usepackage{amsmath}
\usepackage{mathrsfs} 
\usepackage[cal=boondox]{mathalfa}
\usepackage{tablefootnote}
\usepackage[applemac]{inputenc} 		
\usepackage{cleveref}
\usepackage{enumitem}
\usepackage{mathtools}
\usepackage{comment}
\usepackage{tcolorbox}
\usepackage{xcolor,colortbl}
\usepackage{float}
\usepackage{indentfirst,latexsym,bm}
\usepackage{algorithm}
\usepackage{algorithmicx}
\usepackage{algpseudocode}
\usepackage{nicefrac} 
\usepackage{mdframed}

\usepackage{tabularx}
\usepackage{longtable}
\usepackage{multirow}
\usepackage{dsfont}
\usepackage{nicematrix}
\usepackage{pifont}	
\usepackage{ragged2e}
\usepackage{booktabs}

\usepackage{tikz,pgfplots}
\usetikzlibrary{arrows.meta,positioning}
\pgfplotsset{compat=newest}
\usepgfplotslibrary{fillbetween}
\usetikzlibrary{shapes,decorations}
\usetikzlibrary{fit}
\usetikzlibrary{calc}

\definecolor{ivory}{RGB}{218,215,203}

\definecolor{MyBlue}{RGB}{123,144,210} 	
\definecolor{MyRed}{RGB}{224,60,138} 	
\definecolor{DeepGreen}{RGB}{92,172,129} 	
\definecolor{MyOrange}{RGB}{226,148,59} 	
\definecolor{DeepBlue}{RGB}{0,92,175}  
\definecolor{DeepRed}{RGB}{167,67,67}  
\definecolor{LightBlue}{RGB}{229,232,247}   
\definecolor{MyGray}{RGB}{100,100,100} 
\definecolor{lavender}{rgb}{0.9, 0.9, 0.98}

\definecolor{blackp}{RGB}{0,0,0} 
\definecolor{redp}{RGB}{255,0,0}
\definecolor{orangep}{RGB}{255,128,0}
\definecolor{brownp}{RGB}{128,77,0}
\definecolor{yellowp}{RGB}{255,230,0}
\definecolor{greenp}{RGB}{128,230,0}
\definecolor{bluep}{RGB}{0,128,255}
\definecolor{purplep}{RGB}{152,24,147}
\definecolor{pinkp}{RGB}{230,0,128}     
\definecolor{turqp}{RGB}{45, 84, 174} 

\hypersetup{
	colorlinks=true,
	frenchlinks=false,
	pdfborder={0 0 0},
	naturalnames=false,
	hypertexnames=false,
	breaklinks,
	allcolors = blue,
	urlcolor = red!80,
}

\newcommand{\red}[1]{\textcolor{red}{#1}}

\newcommand{\coltab}[1]{\textcolor{turqp}{#1}}

\newcommand{\rmn}[1]{\textup{\textrm{#1}}}

\newcommand{\prt}[1]{\Big(#1 \Big)}
\newcommand{\bra}[1]{\Big[#1 \Big]}
\newcommand{\cbra}[1]{\Big\{#1 \Big\}}

\newcommand{\R}{\mathbb{R}}
\newcommand{\N}{\mathbb{N}}
\newcommand{\Rn}{\mathbb{R}^n}

\newcommand{\Prob}{\mathbb{P}}
\newcommand{\Exp}{\mathbb{E}}

\newcommand{\sA}{{\sf A}}

\newcommand{\sL}{{\sf L}}

\newcommand{\cO}{\mathcal O}
\newcommand{\zetarc}{\bar{\zeta}_r}
\newcommand{\xirc}{\bar{\xi}_r}

\newcommand{\SGD}{{\sf SGD }}
\newcommand{\RR}{{\sf RR }}

\newcommand{\xmark}{\color{DeepRed}\ding{55}}

\newcommand{\ukc}{\zeta}

\newcommand{\sct}[1]{\boldsymbol{#1}}

\newcommand{\vg}{\sct{g}}

\newcommand{\vx}{\sct{x}}

\mdfdefinestyle{questionbox}{
backgroundcolor=LightBlue!10,
    linecolor=blue!60!black, 
    linewidth=1pt,
    topline=true,
    bottomline=true,
    rightline=true,
    leftmargin=0pt,
    innertopmargin=8pt,
    innerbottommargin=8pt,
    innerrightmargin=24pt,
    innerleftmargin=24pt,
    shadow=false,
    shadowcolor=gray!10,
    shadowsize=2pt,
    roundcorner=0pt,
}

\newtheorem{thm}{Theorem}[section]
\newtheorem{lem}[thm]{Lemma}
\newtheorem{cor}[thm]{Corollary}
\newtheorem{prop}[thm]{Proposition}
\newtheorem{rem}[thm]{Remark}
\newtheorem{assmp}[thm]{Assumption}
\newtheorem{defn}[thm]{Definition}

\newtheorem{exmp}[thm]{Example}
\crefalias{prop}{Proposition}

\usepackage{lastpage}
\jmlrheading{23}{2026}{1-\pageref{LastPage}}{1/21; Revised 5/22}{9/22}{21-0000}{}

\ShortHeadings{A Generalized Version of Chung's Lemma}{Jiang, Li, Milzarek and Qiu}
\firstpageno{1}

\begin{document}

\title{A Generalized Version of Chung's Lemma and its Applications}

\author{\name Li Jiang \email lijiang@link.cuhk.edu.cn \\
       \AND
       \name Xiao Li \email lixiao@cuhk.edu.cn \\
       \AND
       \name Andre Milzarek \email andremilzarek@cuhk.edu.cn \\
       \addr 
       School of Data Science \\
       The Chinese University of Hong Kong, Shenzhen \\
       Guangdong 518172, P.R. China  \\
       \AND
       \name Junwen Qiu \email jwqiu@nus.edu.sg \\
       \addr 
       Industrial Systems Engineering \& Management \\
       National University of Singapore \\
       Singapore, 119077, Singapore
       }

\editor{My editor}

\maketitle

\begin{abstract}
Chung's Lemma is a classical tool for establishing asymptotic convergence rates of (stochastic) optimization methods under strong convexity-type assumptions and appropriate polynomial diminishing step sizes. In this work, we develop a generalized version of Chung's Lemma, which provides a simple non-asymptotic convergence framework for a more general family of step size rules. We demonstrate broad applicability of the proposed generalized lemma by deriving tight non-asymptotic convergence rates for a large variety of stochastic methods. In particular, we obtain partially new non-asymptotic complexity results for stochastic optimization methods, such as Stochastic Gradient Descent ($\SGD$) and Random Reshuffling ($\RR$), under a general $(\theta,\mu)$-Polyak-{\L}ojasiewicz (PL) condition and for various step sizes strategies, including polynomial, constant, exponential, and cosine step sizes rules. Notably, as a by-product of our analysis, we observe that exponential step sizes exhibit superior adaptivity to both landscape geometry and gradient noise; specifically, they achieve optimal convergence rates without requiring exact knowledge of the underlying landscape or separate parameter selection strategies for noisy and noise-free regimes. Our results demonstrate that the developed variant of Chung's Lemma offers a versatile, systematic, and streamlined approach to establish non-asymptotic convergence rates under general step size rules.
\end{abstract}


\section{Introduction}

Chung's Lemma \citep{Chu54} is widely used to derive convergence rates of stochastic optimization methods; see, e.g., \citep{bottou2018optimization,gurbu2015,nguyen2020unified,gurbu2019,li2021convergence}. In his seminal paper, \cite{Chu54} analyzed the stochastic approximation method introduced in \citep{RobMon51} for solving a univariate equation, which is considered the precursor to the modern stochastic gradient method ($\SGD$), under the polynomial step size rule $\alpha_k \sim 1/k^p$. Let $\{x^k\}_k$ represent the sequence of iterates generated by the method in \citep{RobMon51}, and let $x^*$ denote the corresponding solution. The recursion for $a_k := \Exp[|x^k - x^*|^2]$ analyzed in \citep{Chu54} can then be expressed as
\begin{equation}
\label{eq:classical-chung-recursion}
a_{k+1} \leqslant \Big({ 1 - \frac{c}{k^p} }\Big) a_k + \frac{d}{k^{p+q}}, \quad \text{where} \quad p\in(0,1], \quad c,d,q>0.
\end{equation}
A crucial task in the convergence analysis is to study how fast $\{a_k\}_k$ converges to 0. Chung provided the following asymptotic estimate.
\begin{tcolorbox}[arc=0pt,outer arc=0pt,colbacktitle=blue!1!white, title filled=true,titlerule=0.2mm,toptitle=0mm,bottomtitle=0mm,top=0mm,
	boxrule=.2mm,
	colback=white,coltitle=black,left=1mm,right=1mm,
	title= Chung's Lemma (Informal),fonttitle=\bfseries]
    \vspace{1mm}
	Let $\{a_k\}_k\subseteq \R_{++}$ be given as in \eqref{eq:classical-chung-recursion}. If $0<p<1$ or $p =  1$, $c>q$, then $a_k = \cO(k^{-q})$.
	\vspace{-1mm}
\end{tcolorbox}
The classical version of Chung's Lemma only applies to polynomial step sizes. In order to go beyond such a specific step size rule and to encompass a broader class of strategies (e.g., exponential and cosine step size rules), we consider the following general recursion:
\begin{equation}
	\label{eq:basic-recursion}
	a_{k+1} \leqslant \Big({1 - \frac{1}{s(b_k)}}\Big) a_k + \frac{1}{t(b_k)},
\end{equation}
where $\{b_k\}_k\subseteq \R_{++}$ is given and $s,t:\R_{++}\to\R_{++}$ are continuously differentiable functions. Clearly, $s(b_k)$ and $t(b_k)$ generalize the terms $k^p/c$ and $k^{p+q}/d$ in \eqref{eq:classical-chung-recursion}. Existing strategies for solving recursions of the form \eqref{eq:basic-recursion} mainly rely on explicit expansions:
	\begin{equation}
		\label{eq:expansion}
		a_{K+1} \leqslant a_0 \prod_{k = 0}^K \Big(1-\frac{1}{s(b_k)} \Big) + \sum_{k=0}^K \Big[{\prod}_{i = k+1}^K \Big(1-\frac{1}{s(b_i)}\Big)\Big]\frac{1}{t(b_k)},
	\end{equation}
	where $\prod_{i = K+1}^K (1-\frac{1}{s(b_i)}):= 1$. While this methodology can seemingly deal with any choice of $s$ and $t$, it requires to derive tight bounds of the different appearing product and sum terms, which can be highly nontrivial and case-dependent; see, e.g., \citep{li2021second}. Based on these observations, our work is motivated by the following core question:


\begin{mdframed}[style=questionbox]
    Can we provide a \emph{simple, non-asymptotic} upper bound for the sequence $\{a_k\}_k$ defined in \eqref{eq:basic-recursion} (as in the classic Chung's Lemma) for more \emph{general step size rules}, i.e., for general $s(b_k),t(b_k)$?
\end{mdframed}

\subsection{Contributions}
In this work, we address the latter core question affirmatively by imposing an additional (mild) structure on $s$ and $t$. 
We now list our main contributions: \\[1mm]
\noindent$\blacktriangleright$ \textbf{Generalized Chung's Lemma.}  We derive a novel general variant of Chung's Lemma to characterize the behavior of $\{a_k\}_k$. An informal version is shown below:
\begin{tcolorbox}[arc=0pt,outer arc=0pt, colbacktitle=blue!1!white, title filled=true,titlerule=0.2mm,toptitle=0mm,bottomtitle=0mm,top=0mm,
	boxrule=.2mm,
	colback=white,coltitle=black,left=1mm,right=1mm,
	title= Generalized Chung's Lemma (Informal \Cref{main-thm}),fonttitle=\bfseries]
	Suppose that $\{a_k\}_k\subseteq \R_{++}$ follows the recursion \eqref{eq:basic-recursion} with $r:=s/t$ being convex on some interval $I$ and $\{b_k\}_k\subseteq I$. Then, under certain conditions, we have for all $K \geqslant 0$
		\begin{equation}\label{eq:informal}
		a_{k+1} ~\leqslant ~\underbracket{C_1 \cdot r(b_{k+1})\vphantom{ \Big(1-\frac{1}{s(b_i)} \Big)}}_{\text{$\mathcal T$-induced rate}}~ + ~\underbracket{C_2 a_0 \cdot {\prod}_{i=0}^{k}\Big(1-\frac{1}{s(b_i)} \Big)}_{\text{$\mathcal S$-induced rate}},\quad \text{for some }\; C_1, C_2>0.
		\end{equation}
\end{tcolorbox} 

As indicated in the informal theorem, our generalized variant of Chung's Lemma naturally decomposes the rate bound into two components, a $\mathcal T$- and an $\mathcal S$-induced rate. The $\mathcal T$-induced (often ``slower'') rate is linked with the ``error term'' ${1}/{t(b_k)}$ which, in turn, is typically associated with {stochastic gradient errors or other types of errors}. The $\mathcal S$-induced (often ``faster'') part represents the convergence rate derived from recursion \eqref{eq:basic-recursion} without the error term ${1}/{t(b_k)}$ and is solely connected to the mapping $s$. 
In Corollary \ref{coro:chung-1}, we present evidence demonstrating that the $\mathcal S$-induced term usually converges at a faster rate compared to the $\mathcal T$-induced part. Thus, the convergence rate of $\{a_k\}_k$ in \eqref{eq:basic-recursion} is $\mathcal{O}(s(b_k)/t(b_k)) = \mathcal{O}(r(b_k))$, which generalizes Chung's classical result $\mathcal{O}(k^{p}/k^{p+q}) = \mathcal{O}(k^{-q})$ for recursion \eqref{eq:classical-chung-recursion}. Overall, our generalized Chung's Lemma provides a simple tool to obtain \emph{non-asymptotic convergence rates} based on the recursion \eqref{eq:basic-recursion}. As a direct application, we utilize \Cref{main-thm} to derive a non-asymptotic extension of the classical Chung's Lemma (cf. Lemma~\ref{lem:non-asymptotic-Chung-main}). \\[1mm]
\noindent$\blacktriangleright$ \textbf{$(\theta,\mu)$-PL Condition, Exponential and Cosine Step Sizes.} Building on the generalized Chung's Lemma, we provide non-asymptotic analyses for $\SGD$ and the random reshuffling method ($\RR$, \Cref{algo:rr}) with exponential and cosine step size rules under a general $(\theta,\mu)$-PL condition; see Theorems~\ref{thm:sgd_exp}, \ref{thm:rr_exp}, \ref{thm:sgd_cos}, and \ref{thm:rr_cos}. Note that the existing analyses for these step size rules primarily focus on $\SGD$ with $\theta = \frac12$. We recover existing complexity bounds and show that, under the standard $(\frac{1}{2},\mu)$-Polyak {\L}ojasiewicz (PL) condition, the convergence rate for cosine step sizes derived in \cite{li2021second} can be improved from $\cO(K^{-2/3})$ to $\cO(K^{-1})$, up to a logarithmic factor, by employing a specialized parameter selection strategy. Moreover, we establish new results for general PL-exponents $\theta \in (\frac{1}{2},1]$. A detailed comparison of our findings to the existing ones is displayed in \Cref{tab:results}. From a technical standpoint, there are two main challenges. The first challenge is that under the general $(\theta,\mu)$-PL condition, we need to analyze a more intricate recursion,
    \begin{equation}\label{eq:new recur}
    y_{k+1} \leqslant (1+\ell_1 \alpha_k^\tau) y_k - \ell_2 \alpha_k y_k^{2\theta} + \ell_3 \alpha_k^\tau,
    \end{equation}
which renders the generalized Chung's Lemma inapplicable when $\theta\in (\frac{1}{2},1]$. We tackle this issue by constructing auxiliary sequences, for which the generalized Chung's Lemma offers important instructional guidance, to simplify \eqref{eq:new recur} to a recursion of the form \eqref{eq:basic-recursion}. We refer to \Cref{sec:recursion&roadmap} for more details. The second challenge arises from the fact that the assumptions required to apply the generalized Chung's Lemma hold only for a subset of the iterates when using exponential and cosine step sizes. We propose a novel splitting technique by combining the generalized Chung's Lemma with an extension lemma (see \Cref{lem:extension}) to solve this problem. See \Cref{sec:exp&cos} for more details. \\[1mm]
%
%
%
\noindent$\blacktriangleright$ \textbf{Polynomial and Constant Step Size Rules.} For completeness, we also provide non-asymptotic analyses for $\SGD$ and $\RR$ with polynomial and constant step size rules under the $(\theta,\mu)$-PL condition based on our generalized Chung's Lemma; see Theorems~\ref{thm:sgd_poly}, \ref{thm:rr_poly}, \ref{thm:sgd_const}, and \ref{thm:rr_const}. Our findings not only corroborate the established results in the literature, but also introduce several new ones, particularly for scenarios where $\theta \in (\frac{1}{2},1]$. A comparison between our results and the existing body of work is provided in \Cref{tab:results}. 

\noindent$\blacktriangleright$ \textbf{Landscape Adaptivity and Noise Adaptivity of Exponential Step Sizes.} Our convergence analysis reveals that exponential step sizes behave distinctly from the other three variants, particularly regarding noise and landscape adaptivity. Most notably, they are the only schedule to exhibit landscape adaptivity---automatically attaining optimal convergence rates without prior knowledge of the PL exponent $\theta$. We validate this finding through numerical simulations. Furthermore, unlike the other three step sizes, which require separate parameter strategies to achieve optimality in noisy versus noiseless settings, exponential step sizes require no such adjustments.

\begin{table}[t]
\centering
\setlength{\tabcolsep}{5pt}
\NiceMatrixOptions{cell-space-limits=1pt}
\begin{NiceTabular}{|p{1.2cm}p{1.5cm}|p{2.3cm}p{2.7cm}|p{2.45cm}p{2.85cm}|}%
 [ 
   code-before = 
    \rectanglecolor{lavender!40}{3-1}{10-1}
    \rectanglecolor{lavender!40}{3-3}{10-3}
    \rectanglecolor{lavender!40}{3-5}{10-5}
 ]
\toprule
\Block[c]{1-2}{\textbf{Algorithms}} & & \Block[c]{1-2}{$\SGD$} & & \Block[c]{1-2}{$\RR$}  \\ \Hline
\Block{1-1}{$\theta$} & \Block{1-1}{step size} & \Block{1-1}{existing} & \Block{1-1}{ours} & \Block{1-1}{existing} & \Block{1-1}{ours} \\ \Hline 
\Block{4-1}{$\frac12$} & \Block{1-1}{Exp} & \Block{1-1}{$\cO(\frac{\log(K)}{K})$\,${}^{\textcolor{blue}{\text{(a)}}}$} & \Block{1-1}{$\cO(\frac{\log(K)}{K})$} & \Block{1-1}{\xmark} & \Block{1-1}{\coltab{$\cO(\frac{\log(K)}{NK^2})$}} \\[2mm]
& \Block{1-1}{Cos\,${}^{\textcolor{blue}{\text{(b)}}}$} & \Block{1-1}{$\cO(\frac{1}{K^{2/3}})$\,${}^{\textcolor{blue}{\text{(c)}}}$} & \Block{1-1}{\coltab{$\cO(\frac{(\log(K))^{\frac{1}{2p+1}}}{K})$}} & \Block{1-1}{\xmark} & \Block{1-1}{\coltab{$\cO(\frac{(\log(NK))^{\frac{2}{2p+1}}}{NK^2})$}}\\[2mm]
& \Block{1-1}{Const} & \Block{1-1}{$\cO(\frac{\log(K)}{K})$\,${}^{\textcolor{blue}{\text{(d)}}}$} & \Block{1-1}{$\cO(\frac{\log(K)}{K})$} & \Block{1-1}{$\cO(\frac{\log^2(NK)}{NK^2})$\,${}^{\textcolor{blue}{\text{(e)}}}$} & \Block{1-1}{$\cO(\frac{\log^{2}(NK)}{NK^2})$} \\[2mm]
& \Block{1-1}{Poly} & \Block{1-1}{$\cO({K^{-1}})$\,${}^{\textcolor{blue}{\text{(f)}}}$} & \Block{1-1}{$\cO({K^{-1}})$} & \Block{1-1}{\hspace{4mm}$\star$\,${}^{\textcolor{blue}{\text{(g)}}}$} & \Block{1-1}{\coltab{$\cO(\frac{\log^{2}(N)}{NK^2})$}} \\[2mm] \Hline
\Block{4-1}{$(\frac12,1]$\,${}^{\textcolor{blue}{\text{(h)}}}$} & \Block{1-1}{Exp} & \Block{1-1}{\xmark} & \Block{1-1}{\coltab{$\cO([{\frac{\log(K/\beta)}{K}}]^{\omega_s})$}} & \Block{4-1}{\xmark} & \Block{1-1}{\coltab{$\cO([{\frac{\log(K/\beta)}{\sqrt{N}K}}]^{\omega_r})$}} \\[2mm]
& \Block{1-1}{Cos} & \Block{1-1}{\xmark} & \Block{1-1}{\coltab{$\cO([{\frac{(\log(K))^{v_s}}{K}}]^{\omega_s})$}} & & \Block{1-1}{\coltab{$\cO([{\frac{(\log(NK))^{v_r}}{\sqrt{N}K}}]^{\omega_r})$}} \\[2mm]
& \Block{1-1}{Const} & \Block{1-1}{\xmark} & \Block{1-1}{\coltab{$\cO([{\frac{\log(K)}{K}}]^{\omega_s})$}} & & \Block{1-1}{\coltab{$\cO([{\frac{\log(NK)}{\sqrt{N}K}}]^{\omega_r})$}} \\[2mm]
& \Block{1-1}{Poly} & \Block{1-1}{$\cO({K}^{-\omega_s})$\,${}^{\textcolor{blue}{\text{(i)}}}$} & \Block{1-1}{{$\cO({K}^{-\omega_s})$}} & & \Block{1-1}{\coltab{$\cO([{\frac{\log(N)}{\sqrt{N}K}}]^{\omega_r})$}} \\[1mm] 
\bottomrule
\end{NiceTabular}
\vspace{1ex}
\caption{Comparison of our results with existing literature. \endgraf
\vspace{1ex} 
\setlength{\parindent}{1ex} \setlength{\baselineskip}{10.5pt}
\noindent \hspace{1pt} ${\footnotesize\textcolor{blue}{\text{(a)}}}$ {\footnotesize See \cite[Theorem 4.1]{wang2021convergence}. The result in \citep{wang2021convergence} is established under the assumptions that the function $ f $ is strongly convex and that the stochastic gradient is bounded. We include it here because their analysis ultimately reduces to a recursion similar to ours. Comparable results are also reported in \citep{li2021second, vaswani2022towards}, albeit with an additional logarithmic factor.}\\
\hspace*{1pt} ${\footnotesize\textcolor{blue}{\text{(b)}}}$ {\footnotesize We use cosine step sizes of the form $\alpha_k = \alpha [\frac{1+\cos(\frac{k\pi}{K})}{2}]^p$, $p>0$; Existing literature only covers $p = 1$.} \\[.5mm]
\hspace*{1pt} ${\footnotesize\textcolor{blue}{\text{(c)}}}$ {\footnotesize See \cite[Theorem 2]{li2021second}. We improve this rate by suggesting a particular strategy for selecting $\alpha$ in the cosine step sizes.} \\[.5mm]
\hspace*{1pt} ${\footnotesize\textcolor{blue}{\text{(d)}}}$ {\footnotesize See \cite[Theorem 4.6]{gower2021sgd}, \cite[Theorem 4]{KarNutSch16}. Several studies have obtained similar results under strong convexity or its variants; see \citep{ShaZha2013,needell2014stochastic,bottou2018optimization,gower2019sgd}.} \\[.5mm]
\hspace*{1pt} ${\footnotesize\textcolor{blue}{\text{(e)}}}$ {\footnotesize See \cite[Theorem 1]{nguyen2020unified}.} \\[.5mm]
\hspace*{1pt} ${\footnotesize\textcolor{blue}{\text{(f)}}}$ {\footnotesize See \cite[Theorem 4]{KarNutSch16}, \cite[Theorem 3]{khaled2020better}, \cite[Theorem 4.8]{gower2021sgd}. It is worth noting that both \cite{khaled2020better,gower2021sgd} did not analyze polynomial step sizes directly, but instead studied a hybrid step size rule that combines constant and polynomial step sizes. Similar results are obtained under strong convexity or its variants; see \citep{Moulines2011,rakhlin2012making,bottou2018optimization,gower2019sgd}.} \\[.5mm]
\hspace*{1pt} ${\footnotesize\textcolor{blue}{\text{(g)}}}$ {\footnotesize 
For polynomial step sizes, $\alpha_k = 2(k+1+\frac1N)^{-1}$, and $\theta = \frac{1}{2}$, \cite[Theorem 7]{nguyen2020unified} establish the complexity bound $\mathcal{O}(\frac{\log(K)}{NK^2})$. This result requires an explicit assumption on the Lipschitz constant $\sL$ of $\nabla f$, ${\sL\sqrt{2(\Theta/N+1)}} \leqslant1 $, where $\Theta$ is a parameter in the variance bound used in \citep{nguyen2020unified}.} \\[.5mm]
\hspace*{1pt} ${\footnotesize \textcolor{blue}{\text{(h)}}}$ {\footnotesize Here, we use the notation $\omega_s = \frac{1}{4\theta-1}$, $\omega_r = \frac{1}{3\theta-1}$, $\varrho_s = \frac{2\theta}{4\theta-1}$, $\varrho_r = \frac{\theta}{3\theta-1}$, $v_s = \frac{\varrho_s}{2p+\varrho_s}$ and $v_r = \frac{\varrho_r}{2p+\varrho_r}$.} \\[.5mm]
\hspace*{1pt} ${\footnotesize \textcolor{blue}{\text{(i)}}}$ {\footnotesize See \cite[Theorem 12]{fontaine2021convergence}; \cite[Corollary 1]{Fatkhullin2022}.} 
\endgraf}
\label{tab:results}
\end{table}

\subsection{Related Works}\label{sec:related-work}
Recursions of the form \eqref{eq:basic-recursion} are abundant in stochastic optimization; see, e.g., \citep{Mis2020,nguyen2020unified,bottou2018optimization,shamir2016,gurbu2015,gurbu2019,ahn2020sgd,huang2021improved}. As mentioned, the functions  $s,t$ in \eqref{eq:basic-recursion} are often related to the step size dynamics. Our overview focuses on Chung's Lemma and its extensions, different step size rules, convergence analyses based on PL-type conditions, and the adaptivity of step sizes. \\[1mm]
\noindent$\blacktriangleright$ \textbf{Chung's Lemma and Extensions.} Chung's Lemma \citep{Chu54} is a fundamental tool extensively used in the literature to analyze the convergence rate of stochastic approximation and optimization algorithms; see \citep{Chu54, derman1956application,Sac58,fabian1968asymptotic,NemJudLanSha09,Moulines2011,gurbu2015,KarNutSch16,bottou2018optimization,gurbu2019,ahn2020sgd,nguyen2020unified,Fatkhullin2022}.  Several extensions of Chung's Lemma are derived in \citep{Moulines2011,ahn2020sgd,nguyen2020unified}, that are tailored to polynomial step size rules. In addition, some studies have considered the recursion \eqref{eq:basic-recursion} with specific non-polynomial forms of $s(b_k)$ and $t(b_k)$; see \citep{li2021second,wang2021convergence}. However, these analyses are conducted only on a case-by-case basis. To the best of our knowledge, there is no generalization of Chung's Lemma that accommodates general terms $s(b_k)$ and $t(b_k)$ in a unified manner. \\[1mm] 
%
%
\noindent$\blacktriangleright$ \textbf{Exponential and Cosine Step Sizes.}
Exponential step sizes \citep{li2021second}, also known as geometrically decaying step sizes \citep{goffin1977convergence,ge2019step}, were used by \cite{goffin1977convergence} for subgradient methods. This step size strategy is also widely used in machine learning, and is currently one of the learning rate decay schedules in TensorFlow \citep{tensorflow2015} and PyTorch \citep{paszke2019pytorch}. 
The theoretical analysis of exponential step size is limited. \cite{ge2019step} provide derivations for least-squares regression problems. Theoretical and broader guarantees for a more general class of functions are shown in \citep{li2021second, wang2021convergence,vaswani2022towards}. {Cosine step sizes}, initially proposed by \cite{loshchilov2016sgdr}, are used as a popular warm restart scheme in performing stochastic algorithms. Recently, cosine step sizes have been widely employed in the optimization and training of large language models. For instance, the Llama-2 model family \citep{touvron2023llama} is pre-trained using an AdamW optimizer with a cosine step size rule. Its fine-tuning procedure uses the same step size scheme. Therefore, understanding the convergence properties of representative stochastic methods with this step size rule is crucial.
Yet, the convergence behavior of stochastic methods under cosine step sizes has not been well understood. The first complexity bounds for $\SGD$ with cosine step sizes are presented in \cite[Theorem 2]{li2021second}. In the general nonconvex setting, extensions to an $\RR$-type method (with momentum) are considered in \cite[Corollary 4]{tran2021smg}. \\[1mm]
\noindent$\blacktriangleright$ \textbf{Stochastic Optimization under PL-type Conditions.} The Polyak-{\L}ojasiewicz (PL) condition \citep{polyak1963,lojasiewicz1959,lojasiewicz1963} provides a measure of the ``strong convexity-like'' behavior of a function without requiring it to be convex. The recent frameworks discussed in \citep{fontaine2021convergence,Fatkhullin2022} rely on certain generalizations of the PL condition. In this paper, we consider the $(\theta,\mu)$-PL condition (\Cref{def:PL}), which is identical to Assumption 3 in \citep{Fatkhullin2022}, and encompasses the PL condition when $\theta = \frac{1}{2}$. The additional PL exponent, $\theta \in [\frac{1}{2},1]$, enables this condition to hold for a wider range of functions, particularly those relevant to machine learning applications (see \cite[Appendix A]{Fatkhullin2022} for pertinent examples). Specifically, in the field of reinforcement learning, the objective function can satisfy the $(1,\mu)$-PL condition for some $\mu > 0$; see \citep{wang2022policy,agarwal2021theory,mei2020global}. These evidences indicate that analyzing the cases beyond the standard setting $\theta = \frac{1}{2}$ holds significance in real applications.\\[1mm]
\noindent$\blacktriangleright$ \textbf{Adaptivity of Step Sizes in Stochastic Optimization Methods.} Adaptivity of algorithms is one of the central concerns in stochastic optimization, which seeks to achieve strong performance without relying on additional, often impractical, information. A substantial body of literature has investigated algorithmic adaptivity in this context, including but not limited to \citep{duchi2011adaptive, li2021second, kingma2014adam, orabona2021parameter, loizou2021stochastic, vaswani2022towards}; for a comprehensive overview, see \citep{attia2024free}. In this paper, we focus primarily on the noise adaptivity and landscape adaptivity of various offline step size schedules for stochastic gradient descent ($\SGD$) and random reshuffling ($\RR$). The concept of noise adaptivity, first introduced by \cite{li2021second}, refers to a step size that achieves optimal convergence rates in both noisy and noiseless settings without requiring prior knowledge of the noise level. Under the $(\frac{1}{2},\mu)$-PL condition, \cite{li2021second} demonstrates that the exponential step size exhibits noise adaptivity, whereas the cosine step size achieves only a suboptimal convergence rate in the noisy setting. A similar discussion on exponential step sizes also appears in \citep{vaswani2022towards}. A comprehensive understanding of noise adaptivity under the more general $(\theta,\mu)$-PL condition remains unexplored. In addition, we introduce the concept of landscape adaptivity, which describes the ability of a step size schedule to achieve the optimal convergence rate under the $(\theta,\mu)$-PL condition without requiring explicit knowledge of $\theta$. We demonstrate that, among the four step size schedules considered, only the exponential step size exhibits this form of adaptivity.

\section{Generalized Chung's Lemma} 
The following generalized version of Chung's Lemma provides a new tool for analyzing the non-asymptotic upper bound for sequences $\{a_k\}_k$ satisfying
\begin{equation}
	\label{eq:general-recursion}
	a_{k+1} \leqslant \Big({1 - \frac{1}{s(b_k)}}\Big)a_k + \frac{1}{t(b_k)},\quad \forall~ k \in [K]:=\{0,1,\dots,K\}.
\end{equation}

\begin{thm}[Generalized Chung's Lemma]\label{main-thm}
	Let $\{a_k\}_k \subseteq \R,\, \{b_k\}_k \subseteq \R $ and functions $s,t: \R \to \R$ be given. Assume that $\{a_k\}_k$ follows \eqref{eq:general-recursion} and there is an interval $I$ such that $b_k \in I$ for all $k\in[K]$, $s(x)\geqslant 1$ and $t(x)>0$ for all $x\in I$, and $r: = s/t$ is convex on $I$. If either of the following statements holds:
    \vspace{1mm}
	\begin{enumerate}[label=\textup{\textrm{(\alph*)}},leftmargin=8ex,topsep=2pt,itemsep=0ex,partopsep=0ex,nosep]
        \item There exists a subgradient $g_k \in \partial r(b_k)$, and a non-decreasing sequence $\{\lambda_k\}_k \subseteq \R_{++}$ such that for all $k\in [K-1]$, 
        \begin{equation}\label{condition1}
			(b_{k+1}-b_k) \cdot g_k \, t(b_k) \geq -1 + \frac{1}{\lambda_{k+1}}.
		\end{equation}
        \vspace{1mm}
		\item The function $r$ is differentiable on $I$ and defining $u:= (\log r)'s = r^\prime t$, there exists a non-decreasing sequence $\{\lambda_k\}_k \subseteq \R_{++}$ such that for all $k\in [K-1]$,
            \vspace{-2mm}
		\begin{equation}\label{condition2}
			(b_{k+1}-b_k) \cdot u(b_k) \geq -1 + \frac{1}{\lambda_{k+1}},
		\end{equation}
            \vspace{-5mm}
	\end{enumerate}
        then we have
        \vspace{-3mm}
		\begin{equation} \label{eq:gen-chung}
		    a_{k+1} \leqslant \lambda_{k+1} r(b_{k+1}) + ({a_0 - \lambda_0 r(b_0)}) \prod_{i = 0}^k \Big({1 - \frac{1}{s(b_i)}}\Big),\quad \forall~ k\in [K-1].
		\end{equation}
		\vspace{-5mm}
\end{thm}
\begin{proof}
    Let us set $\beta_k := b_{k+1}-b_k$. By the convexity of the mapping $r$, we have
	\begin{equation}\label{eq:proof-main-thm1}
	    r(b_{k+1}) \geqslant r(b_k) + g_k \beta_k = r(b_k)\Big(1+\frac{g_k}{r(b_k)}\beta_k\Big) = r(b_k)\Big(1+\frac{g_k t(b_k)}{s(b_k)}\beta_k\Big).
	\end{equation}
If statement (a) holds, substituting the latter estimate into \eqref{condition1}, this yields
\begin{align*}
    \lambda_{k+1} r(b_{k+1}) &\geqslant \lambda_{k+1}r(b_k)\Big(1+\frac{-1 + \frac{1}{\lambda_{k+1}}}{s(b_k)}\Big) = r(b_k)\frac{\lambda_{k+1}(s(b_k)-1)+1}{s(b_k)}\\[1mm]& \geqslant r(b_k)\frac{\lambda_{k}(s(b_k) -1) + 1}{s(b_k)} = \lambda_{k} r(b_k) 
	\Big({1 - \frac{1}{s(b_k)}}\Big) + \frac{1}{t(b_k)}.
\end{align*}
	Hence, invoking the recursion \eqref{eq:general-recursion}, we obtain $a_{k+1} - \lambda_{k+1} r(b_{k+1}) \leqslant (1 - \frac{1}{s(b_k)})(a_{k} - \lambda_{k} r(b_{k}))$. This verifies \eqref{eq:gen-chung}. When statement (b) is satisfied, we have $g_k = r'(b_k)$. Noting $u = \frac{r'}{r} s = r' t$, it holds that $u(b_k) = g_k t(b_k)$. Therefore, condition \eqref{condition2} immediately implies condition \eqref{condition1}, from which \eqref{eq:gen-chung} follows.
\end{proof}

The most significant setting in \Cref{main-thm} arises when the sequence $\{\lambda_k\}_k$ can be chosen to be constant. This will form the basis for the majority of our subsequent analysis. We now provide a simple example to briefly illustrate the application of \Cref{main-thm}.

\begin{exmp}\label{ex:chung1}
    Let us consider the recursion $a_{k+1}\leqslant (1-\alpha)a_k +  \beta$, where $\alpha\in (0,1)$, $\beta>0$. Setting $I = \R$, $s\equiv \frac{1}{\alpha}$, $t\equiv \frac{1}{\beta}$, $b_k = k$, $\lambda = 1$, the rate mapping $r\equiv \frac{\beta}{\alpha}$ is convex and we have $(b_{k+1}-b_k)u(b_k) = 0 = -1+ \frac{1}{\lambda}$ for all $k$. Consequently, \Cref{main-thm} is applicable and we have $a_{k}\leqslant \frac{\beta}{\alpha} + (1-\alpha)^k(a_0-\frac{\beta}{\alpha})$ for all $k$. 
\end{exmp}


\begin{rem}\label{rem:convexity-r}
	The convexity requirement on $r:=s/t$ in \Cref{main-thm} is not strict since there is flexibility in choosing $\{b_k\}_k$. Consider the example $r(b_k) = 1+\cos(\frac{k\pi}{K})$, $k \in [K]$. The direct choice $b_k := k$, $r(x):= 1+\cos(\frac{x\pi}{K})$, and $I=[0,K]$ does not meet the convexity condition. However, we can ensure the convexity of $r$ by setting $b_k :=1+\cos(\frac{k\pi}{K})$, $r(x):= x$, and $I = \R$. This strategy proves valuable in the analysis of cosine step sizes, see \Cref{sec:exp&cos}.
\end{rem}

\noindent\textbf{Non-asymptotic Extension of Chung's Lemma.} As a quick application, we use \Cref{main-thm} to obtain a non-asymptotic extension of the classical Chung's Lemma. 

\begin{lem}[Non-asymptotic Chung's Lemma] \label{lem:non-asymptotic-Chung-main}
	Let $\{a_k\}_k\subseteq \R_{++}$ be given satisfying \begin{equation*}
	a_{k+1} \leqslant \Big(1-\frac{c}{(k+\gamma)^\nu}\Big) a_k + \frac{d}{(k+\gamma)^{\nu+q}},\quad \nu\in(0,1],~ q>0, ~c,d > 0,~\gamma \geqslant c^{\frac{1}{\nu}}.
\end{equation*}
        Then, the following statements hold:
	\begin{enumerate}[label=\textup{\textrm{(\alph*)}},leftmargin=8ex,topsep=2pt,itemsep=0ex,partopsep=0ex]		
		\item If $\nu\in (0,1)$, $\gamma \geqslant c^{\frac{1}{\nu}}$ and $\gamma > (\frac{q}{c})^{\frac{1}{1-\nu}}$, then 
		\begin{align*}
			a_{k+1} \leqslant \frac{\lambda d}{c} (k+1+\gamma)^{-q} +  \Big(a_0 - \frac{\lambda d}{c\gamma^{q}}\Big)^+	\exp\Big(\frac{c\gamma^{1-\nu}}{1-\nu}\Big) \exp\Big(-\frac{c(k+1+\gamma)^{1-\nu}}{1-\nu}\Big),
		\end{align*}
		holds for all $k \in [K]$, where $\lambda = \frac{c\gamma^{1-\nu}}{c\gamma^{1-\nu}-q}$ and $(x)^+:=\max\{x,0\}$. 
		\item If $\nu = 1$, $c>q$, and $\gamma \geqslant c$, then
		\[
              a_{k+1} \leqslant \frac{d}{c-q} (k+1+\gamma)^{-q} + \gamma^c\Big(a_0 - \frac{ d}{(c-q)\gamma^{q}}\Big)^+ (k+1+\gamma)^{-c}.
		\]
	\end{enumerate}
\end{lem}

\begin{proof}
We apply \Cref{main-thm} (b) to prove the stated non-asymptotic version of the classical Chung's Lemma. To make use of \Cref{main-thm} (b), we first specify the following sequences and functions: 
		\[
		b_k := k+\gamma,\quad  s(x) := {x^\nu}/{c} , \quad \text{and}\quad t(x) := {x^{\nu+q}}/{d}.
		\] 
		We only need to check if the conditions in \Cref{main-thm}\,(b) are fulfilled. Let us note that the mapping $r(x) := s(x)/t(x) = \frac{d}{c} x^{-q}$ is convex on $\mathbb{R}_{++}$. Then $(\log(r))'(x) = -\frac{q}{x}$, and it holds that
		\[
		u(x) : = (\log(r))'(x) s(x) = \Big(-\frac{q}{x}\Big)\cdot \Big(\frac{x^\nu}{c}\Big) = -\frac{q}{c} x^{\nu-1}
		\]
        and $s(x)\geqslant 1$ for all $x\geqslant c^{{1}/{\nu}}$. By assumption, we have $\gamma \geqslant c^{{1}/{\nu}}$ implying $b_k \in I := [c^{{1}/{\nu}},\infty) \subseteq \R_{++}$ for all $k \in [K]$. Now, it only remains to verify the condition \eqref{condition2}. Since $b_{k+1}-b_k = 1$, we need to find a non-decreasing sequence $\{\lambda_k\}_k \subseteq \R_{+}$ such that 
            \begin{equation}
                \label{condition2:poly-Chungs-main}
                -\frac{q}{c}(k+\gamma)^{\nu-1} \geqslant -1 + \frac{1}{\lambda_{k+1}},\quad \forall\;k\in[K].
            \end{equation}
		Based on the parameter $p \in (0,1]$, we now discuss different cases. \\[1mm] 
        \noindent\textbf{Case I:} $\nu\in(0,1)$. Let us set $\lambda_k \equiv \lambda = \frac{c\gamma^{1-\nu}}{c\gamma^{1-\nu}-q}$. Per assumption, we further have $\gamma > (\frac{q}{c})^{\frac{1}{1-\nu}}$ and thus, we can infer $\lambda > 0$ and 
        \begin{equation} \label{eq:nice-am-01} -\frac{q}{c}(k+\gamma)^{\nu-1} \geqslant -\frac{q}{c}\gamma^{\nu-1} = -1 + \frac{c\gamma^{1-\nu}-q}{c\gamma^{1-\nu}} = -1 + \frac{1}{\lambda} \quad \forall~k\in[K]. \end{equation}
        Hence, the condition \eqref{condition2:poly-Chungs-main} is satisfied for this choice of $\lambda_k$ and $\lambda$. Using \Cref{main-thm}\,(b), the inequality $1-x \leqslant \exp(-x)$, and the integral test (\Cref{prop:int}), we obtain
        \begingroup
        \allowdisplaybreaks
		\begin{align*}
			a_{k+1} & \leqslant \lambda r(b_{k+1}) + (a_0 - \lambda r(\gamma)) \cdot  {\prod}_{i=0}^k \Big(1-\frac{1}{s(b_i)}\Big) \\ & = \frac{\lambda d}{c} (k+1+\gamma)^{-q} + \Big(a_0 - \frac{\lambda d}{c\gamma^{q}}\Big) {\prod}_{i = 0}^k \Big(1-\frac{c}{(i+\gamma)^\nu}\Big)\\ 
			&\leqslant \frac{\lambda d}{c} (k+1+\gamma)^{-q} + \Big(a_0 - \frac{\lambda d}{c\gamma^{q}}\Big)^+ \exp\Big(-{\sum}_{i=0}^{k} \frac{c}{(i+\gamma)^\nu}\Big)\\[2mm]
			&\leqslant \frac{\lambda d}{c} (k+1+\gamma)^{-q} + \Big(a_0 - \frac{\lambda d}{c\gamma^{q}} \Big)^+	\exp\Big(\frac{c\gamma^{1-\nu}}{1-\nu}\Big) \exp\Big(-\frac{c(k+1+\gamma)^{1-\nu}}{1-\nu}\Big),
		\end{align*}
        \endgroup
		where the last line is due to $\sum_{i=0}^{k}\, (i+\gamma)^{-\nu} \geq \int_0^{k+1} (x+\gamma)^{-\nu} \rmn{d}x = \frac{1}{1-\nu}[(k+1+\gamma)^{1-\nu}-\gamma^{1-\nu}]$. \\[1mm]
		\noindent\textbf{Case II:} $\nu=1$, $c>q$, and $\gamma \geqslant c$. Setting $\lambda_k \equiv \lambda = \frac{c}{c-q}$, the condition \eqref{condition2:poly-Chungs-main} is obviously satisfied (cf. \eqref{eq:nice-am-01}). Repeating the core steps in \textbf{Case I}, we obtain
		\[
		\begin{aligned}
			a_{k+1} \leqslant ~&\frac{\lambda d}{c} (k+1+\gamma)^{-q} + \Big(a_0 - \frac{\lambda d}{c\gamma^{q}}\Big) {\prod}_{i = 0}^k \Big(1-\frac{c}{i+\gamma}\Big)\\
			\leqslant ~& \frac{d}{c-q} (k+1+\gamma)^{-q} + \gamma^c\Big(a_0 - \frac{d}{(c-q)\gamma^q}\Big)^+ (k+1+\gamma)^{-c}.
		\end{aligned}
		\]
		Notice that the difference in the analysis only stems from the integral test (Proposition~\ref{prop:int}): ${\sum}_{i=0}^{k}\, \frac{1}{i+\gamma} \geqslant {\int}_{\hspace{-1mm}0}^{k+1} \frac{1}{x+\gamma} \;\rmn{d}x = \log(k+1+\gamma) - \log(\gamma)$.
		\end{proof}


 

\noindent\textbf{Extension Lemma.} We now introduce an additional lemma to enhance \Cref{main-thm}(b) for cases where constructing $\{\lambda_k\}_k$ is only possible for a subset of $[K]$, such as with exponential and cosine step sizes (\Cref{sec:exp&cos}). Our strategy involves partitioning $[K]$: we first find the maximal scalar $\overline{K}$ that satisfies \eqref{condition2} to bound $a_{\overline{K}+1}$ by \Cref{main-thm}(b); in a second step, \Cref{lem:extension} can be used to extend this estimate to $a_K$. 

\begin{lem}[Extension Lemma] \label[lemma]{lem:extension}
			Let $\{a_k\}_k \subseteq \R_{++}$, $\{b_k\}_k \subseteq \R_+ $ and the functions $s,t: \R \to \R_{++}$ be given. Suppose that $\{a_k\}_k$ follows the recursion \eqref{eq:general-recursion}. Assume that $s(b_k) \geqslant 1$ for all $k\in[K]$. If there exist $k_0,\overline{K}\in [K-1]$, $k_0\leqslant \overline{K}$, $B, C \in \R$ such that $r(b_k)\leqslant B$ for all $k \geqslant \overline{K}+1$ and
			\begin{equation}\label{eq:extesion1}
			    a_{\overline{K}+1} \leqslant B + C \prod_{i = k_0}^{\overline{K}} \Big({1 - \frac{1}{s(b_i)}}\Big),
			\end{equation}
			then, we have $a_{K} \leqslant B + C \prod_{i = k_0}^{K-1} ({1 - \frac{1}{s(b_i)}})$.
		\end{lem}
        
 \begin{proof}
There is nothing to show in the case $\overline{K} = K-1$. We can verify the result for $\overline{K} < K-1$ via induction. For any $k = \overline{K}+1,\dots,K-1$, if 
			$a_{k} \leqslant B + C \prod_{i = k_0}^{k-1} ({1 - \frac{1}{s(b_i)}})$,
and by applying \eqref{eq:general-recursion}, we obtain
            \begingroup
            \allowdisplaybreaks
			\begin{align*}
				a_{k+1} - C {\prod}_{i = k_0}^{k} \Big({1 - \frac{1}{s(b_i)}}\Big) \leqslant ~& \Big({1 - \frac{1}{s(b_k)}}\Big)a_k + \frac{1}{t(b_k)} 
				- C {\prod}_{i = k_0}^{k} \Big({1 - \frac{1}{s(b_i)}}\Big) \\
				\leqslant~& \Big({1 - \frac{1}{s(b_k)}}\Big)B +\frac{1}{s(b_k)} r(b_k) 
				\leqslant B,
			\end{align*}
            \endgroup
			where we used $r(b_k) \leqslant B$ for any $k\geqslant \overline{K}+1$ in the last step. Hence, by induction, it follows $a_{K} \leqslant B + C \prod_{i = k_0}^{K-1} ({1 - \frac{1}{s(b_i)}})$.
		\end{proof}

%
\noindent\textbf{Forgetting Initial Conditions.} The following corollary  offers insight into the relationship between the $\mathcal S$- and $\mathcal T$-induced components in \eqref{eq:informal}.
\begin{cor} \label[corollary]{coro:chung-1}
	Let the conditions in \Cref{main-thm}\,\textup{(b)} hold with $\lambda_k \equiv \lambda >0$. Then, we have
	\[
	a_{k+1} \leqslant \lambda r(b_{k+1}) +  \Big(\frac{a_0}{r(b_0)} - \lambda\Big)^+ \prod_{i=0}^{k} \Big(1-\frac{1/\lambda}{s(b_i)-1+1/\lambda}\Big) r(b_{k+1}),\quad \forall\; k\in [K].
	\]
\end{cor}

\begin{proof}
		From \Cref{main-thm}\,(b), we have for all $k \in  [K]$,
        \begingroup
        \allowdisplaybreaks
		\begin{align*}
			a_{k+1} \leqslant ~&\lambda r(b_{k+1}) + \Big(\frac{a_0}{r(b_0)} - \lambda\Big)^+{\prod}_{i = 0}^k \Big(1 - \frac{1}{s(b_i)}\Big) r(b_0)\\
			= ~& \lambda r(b_{k+1}) + \Big(\frac{a_0}{r(b_0)} - \lambda\Big)^+{\prod}_{i = 0}^k \Big[\Big(1 - \frac{1}{s(b_i)}\Big)\frac{r(b_i)}{r(b_{i+1})}\Big] r(b_{k+1}).
		\end{align*}
  \endgroup
		Combining equation \eqref{eq:proof-main-thm1} (where $u(b_k) = g_kt(b_k)$) and \eqref{condition2}, we obtain
		\[
		r(b_{i+1}) \geqslant \frac{(b_{i+1}-b_i) u(b_i)+s(b_i)}{t(b_i)} \geqslant \frac{s(b_i)-1+{1}/{\lambda}}{t(b_i)} > 0.
		\]
		Thus, we can establish the upper bound
		\[
		\Big(1 - \frac{1}{s(b_i)}\Big)\frac{r(b_i)}{r(b_{i+1})} \leqslant  \Big(1 - \frac{1}{s(b_i)}\Big) \frac{s(b_i)}{s(b_i)-1+{1}/{\lambda}}
		=  1 - \frac{1/\lambda}{s(b_i)-1+{1}/{\lambda}}.
		\]
		Substituting this bound into the first estimate finishes the proof.
	\end{proof}
 
\Cref{coro:chung-1} illustrates that under \Cref{main-thm}(b), the $\mathcal S$-induced rate outpaces the $\mathcal T$-induced rate by a factor of $\prod_{i=0}^{k} (1-\frac{1/\lambda}{s(b_i)-1+1/\lambda})$. Since the initial point affects only the $\mathcal{S}$-induced rate, this phenomenon is often termed ``forgetting initial conditions''. While prior literature addresses this property only in specific settings---such as for polynomial step sizes in \citep{Moulines2011}---\Cref{coro:chung-1} holds in a more general context. 

\section{\texorpdfstring{Convergence Rates of Stochastic Methods under the ($\theta$,$\mu$)-Polyak {\L}ojasiewicz (PL) Condition}{Convergence Rates of Stochastic Methods under the PL Condition}} \label{sec:app-2}

In this section, we establish non-asymptotic rates for stochastic methods under various step size strategies by applying the generalized Chung's Lemma---\Cref{main-thm}. We consider the problem
\begin{equation} \label{eq:sec03:problem}
    \min_{x\in\Rn} \; f(x),
\end{equation}
where $f:\Rn\to\R$ is a continuously differentiable function. 

\subsection{Basic Assumptions and Descent Properties}\label{sec:assum&descent}

In \Cref{algo:sgd,algo:rr}, we summarize the main algorithmic steps of the stochastic gradient descent ($\SGD$) method and the random reshuffling ($\RR$) method to solve the optimization problem \eqref{eq:sec03:problem}.  We note that, in $\RR$, the objective function $f$ is assumed to have the finite-sum structure $f(x):=\frac{1}{N}\sum_{i=1}^N f_i(x)$. 

\begin{algorithm}[t]
			\caption{Stochastic Gradient Descent ($\SGD$) }
			\label{algo:sgd}
			\begin{algorithmic}[1] 
				\Require Choose an initial point $x^0\in\Rn$ and step sizes $\{\alpha_k\}_{k} \subset \R_{++}$;
				\While{stopping criterion is not met}
				\State Generate a new stochastic gradient $g^k\approx\nabla f(x^k)$;           
				\State Perform the update $x^{k+1}=x^{k}-\alpha_k g^k$;             
				\EndWhile
			\end{algorithmic}
		\end{algorithm}
		
		\begin{algorithm}[t]
			\caption{Random Reshuffling ($\RR$) }
			\label{algo:rr}
			\begin{algorithmic}[1] 
				\Require Choose an initial point $x^0\in\Rn$ and step sizes $\{\alpha_k\}_{k} \subset \R_{++}$;
				\While{stopping criterion is not met}
				\State Generate a uniformly random permutation $\pi^{k}$ of $\{1,\ldots, N\}$ and set $\tilde x_0^{k}=x^{k}$; 
				\For{$i = 1,2,\ldots, N$}
				\State Update the inner iterate via 
				$
				\tilde x_{i}^{k}=\tilde x_{i-1}^{k}-\frac{\alpha_{k}}{N}\, \nabla {f}_{\pi^{k}_{i}}(\tilde x_{i-1}^{k});
				$        
				\EndFor
				\State Perform the update $x^{k+1} = \tilde x_{N}^{k}$;
				\EndWhile
			\end{algorithmic}
		\end{algorithm}
        
Our analysis will be based on a filtered probability space $(\Omega,\mathcal F,\{\mathcal F_k\}_k,\Prob)$ to which the stochastic processes $\{\vx^k\}_k$ generated by these two algorithms are adapted. We make the following standard assumptions for $\SGD$ and $\RR$.

\begin{assmp}[Conditions for $\SGD$]\label[assumption]{assum:sgd} We consider the assumptions:
    \begin{enumerate}[label=\textup{\textrm{(A.\arabic*)}},topsep=2pt,itemsep=0ex,partopsep=0ex,leftmargin=8ex]
	\item \label{A1} The objective $f$ is $\sL$-smooth, i.e., $\nabla f:\Rn\to\Rn$ is $\sL$-Lipschitz continuous.
    \item \label{A2} The stochastic gradient $\vg^k:\Omega \to \Rn$ is an unbiased estimator of $\nabla f(\vx^k)$, and there exist $\sA,\sigma\geqslant 0$, $\bar f \in \R$ such that $\Exp[\|\vg^k - \nabla f(\vx^k)\|^2 \mid \mathcal F_{k-1}] \leqslant \sA(f(\vx^k)-\bar f) +  \sigma^2$ a.s.\,.
    \end{enumerate}
\end{assmp} 

\begin{assmp}[Conditions for $\RR$] \label[assumption]{assum:rr1} We work with the conditions:
    \begin{enumerate}[label=\textup{\textrm{(B.\arabic*)}},topsep=2pt,itemsep=0ex,partopsep=0ex,leftmargin=8ex]
    \item \label{B1} The individual component functions $f_1,f_2,\dots,f_N$ are all $\sL$-smooth.
    \item \label{B2} There exist constants $\sA,\sigma>0$, $\bar f \in \R$ such that for all $x\in\Rn$, it holds that
	\[ \frac{1}{N} \, {\sum}_{i=1}^{N}\|\nabla f_i(x) - \nabla f(x)\|^2  \leqslant \sA(f(x)-\bar f) +  \sigma^2. \]
    \item \label{B3} The permutations $\{\pi^k\}_k$ are sampled independently (for each $k$) and uniformly without replacement from $\{1,\dots,N\}$.
    \end{enumerate}
\end{assmp}

The Lipschitz smoothness in \ref{A1} and \ref{B1} is ubiquitous in the analysis of optimization methods, see, e.g., \citep{BolSabTeb14,gurbu2015,haochen2019}. The assumptions \ref{A2} and \ref{B2} have been extensively employed in the analysis of stochastic gradient-based methods \citep{nguyen2020unified, Mis2020,shamir2016,sohl2014fast}. In the case $\sA=0$, \ref{A2} and \ref{B2} reduce to bounded variance conditions. Furthermore, \ref{B2} is applicable in a wide range of scenarios. For instance, \ref{B2} is satisfied when each individual component function $f_i$ is Lipschitz smooth and bounded from below, see \cite[Proposition 2]{Mis2020}. 


For both $\SGD$ and $\RR$, we employ a fundamental property that describes the geometric structure of $f$, largely following \cite[Assumption 3]{Fatkhullin2022}.

\begin{defn}[($\theta$,$\mu$)-PL Condition] \label[definition]{def:PL}
	Given $\theta\in[\frac12,1]$, $\mu>0$, the function $f$ is said to satisfy the ($\theta$,$\mu$)-PL condition if $f$ is bounded from below by $f^*:= \inf_{x\in\Rn} f(x)>-\infty$ and  
	\[
	\|\nabla f(x)\| \geqslant \sqrt{2\mu} (f(x)-f^*)^{\theta}, \quad \forall\, x\in \R^n.
	\]
\end{defn}

\Cref{def:PL} covers the well-known, classical Polyak-{\L}ojasiewicz condition \citep{polyak1963} as a special case when $\theta = \frac{1}{2}$. 

\begin{assmp}[PL Condition for $f$] We consider the following assumption:
\begin{enumerate}[label=\textup{\textrm{(P)}},topsep=2pt,itemsep=0ex,partopsep=0ex,leftmargin=8ex]
    \item \label{P} For some $\theta\in [\frac{1}{2},1]$, $\mu>0$, the objective function $f$ satisfies the ($\theta,\mu$)-PL condition.
\end{enumerate}    
\end{assmp}

Under assumption \ref{P} and without loss generality, we set $\bar f = f^*$ in \ref{A2} and \ref{B2}. We start with (approximate) descent-type properties for $\SGD$ and $\RR$.

\begin{lem}[Descent-type Properties]\label[lemma]{lem:descent-property}
	Suppose assumption \ref{P} holds.
	\begin{enumerate}[label=\textup{\textrm{(\alph*)}},leftmargin=8ex,topsep=2pt,itemsep=0ex,partopsep=0ex]		
		\item Let $\{\vx^k\}_k$ be generated by \Cref{algo:sgd} with step sizes $\alpha_k\in(0,\frac{1}{\sL}]$ and let the assumptions \ref{A1}--\ref{A2} be satisfied. Then, it holds that
		\begin{equation}\label{eq:sgd-recursion}
		\mathbb{E}[f(\vx^{k+1})-f^*] \leqslant \Big(1+\frac{\sA \sL \alpha_k^2}{2}\Big)\mathbb{E}[f(\vx^{k})-f^*]  - \mu \alpha_k (\mathbb{E}[f(\vx^{k})-f^*] )^{2\theta} + \frac{\sL \sigma^2\alpha_k^2}{2}.
		\end{equation}
		\item Let the sequence $\{\vx^k\}_k$ be generated by \Cref{algo:rr} with step sizes $\alpha_k\in(0,\frac{1}{2\sL}]$ and let \ref{B1}--\ref{B3} hold. Then, we have
		\begin{equation}\label{eq:rr-recursion}
		\Exp[f(\vx^{k+1})-f^*] \leqslant \Big(1+ \frac{\sA\sL^2\alpha_k^3}{2N}\Big)\mathbb{E}[f(\vx^{k})-f^*]  -  \frac{\mu\alpha_k}{2} (\mathbb{E}[f(\vx^{k})-f^*])^{2\theta} + \frac{\sL^2 \sigma^2\alpha_k^3}{2N}.
		\end{equation}
	\end{enumerate}
\end{lem}

\begin{proof}
                For \Cref{algo:sgd}, with assumptions \ref{A1}--\ref{A2}, and $\alpha_k\in (0,\frac{1}{\sL}]$, we have
                \[
            	\Exp[f(\vx^{k+1})-f^*] \leqslant \Big(1+ \frac{\sA\sL\alpha_k^2}{2}\Big)\mathbb{E}[f(\vx^{k})-f^*]  -  \frac{\alpha_k}{2} \mathbb{E}[\|\nabla f(\vx^k)\|^2] + \frac{\sL \sigma^2\alpha_k^2}{2}.
                \]
                For \Cref{algo:rr}, we refer to the result in \cite[Appendix 9.5, Proof of Theorem 4, eq. (47)]{Mis2020}. With \ref{B1}--\ref{B3}, for $\alpha_k\in (0,\frac{1}{2\sL}]$, we have
            	\[
            	\Exp[f(\vx^{k+1})-f^*] \leqslant \Big(1+ \frac{\sA\sL^2\alpha_k^3}{2N}\Big)\mathbb{E}[f(\vx^{k})-f^*]  -  \frac{\alpha_k}{4} \mathbb{E}[\|\nabla f(\vx^k)\|^2] + \frac{\sL^2 \sigma^2\alpha_k^3}{2N}.
            	\]
                Combining $\|\nabla f(x)\|^2 \geqslant 2\mu (f(x)-f^*)^{2\theta}$ in \ref{P} and H$\text{\"o}$lder's inequality, we obtain
                \[
                \Exp [\|\nabla f(x)\|^2] \geqslant 2\mu \Exp[(f(x)-f^*)^{2\theta}] \geqslant 2\mu \Exp[(f(x)-f^*)]^{2\theta},
                \]
                which finalizes the proof.
            \end{proof}

\subsection{Main Recursion and Analysis Roadmap}\label{sec:recursion&roadmap}

From Lemma \ref{lem:descent-property}, it is evident that under the $(\theta,\mu)$-PL framework, the problem reduces to analyzing the recursion with the given form:
\begin{equation}\label{PL-recursion}
	y_{k+1} \leqslant (1+\ell_1 \alpha_k^\tau) y_k - \ell_2 \alpha_k y_k^{2\theta} + \ell_3 \alpha_k^\tau,\quad k\in [K-1],
\end{equation}
where $y_k = \Exp[f(\vx^{k})-f^*]$. The constants for each method are given by:
\[
\SGD:~\ell_1 = \frac{\sA \sL}{2},~\ell_2 = \mu,~\ell_3 = \frac{\sL \sigma^2}{2},~\tau = 2;\quad\RR:~\ell_1 = \frac{\sA\sL^2}{2N},~\ell_2 = \frac{\mu}{2},~\ell_3 = \frac{\sL^2 \sigma^2}{2N},~\tau = 3.
\]
From now on, we assume $\ell_3>0$. Recursion \eqref{PL-recursion} generally falls outside the scope of the generalized Chung's Lemma when $\theta \neq \frac{1}{2}$. To address this issue, we construct an auxiliary sequence $\{u_k\}_k$ that transforms \eqref{PL-recursion} into a recursion of the form \eqref{eq:general-recursion}. The generalized Chung's Lemma is employed not only to analyze the resulting simplified recursion but also to guide the design of $\{u_k\}_k$. Here, we provide a general roadmap for the analysis.

\vspace{2mm}
\noindent\textbf{Step 1: Linearization and Relaxation}. 
            When $\theta = \frac{1}{2}$, \eqref{PL-recursion} can be reformulated as
            \begin{equation*}
            	y_{k+1} \leqslant (1+\ell_1 \alpha_k^\tau - \ell_2 \alpha_k) y_k  + \ell_3 \alpha_k^\tau,\quad k\in [K-1].
            \end{equation*}
            Define
            \begin{equation}\label{eq:def-maxstep}
                \bar{\alpha}=\max \{\alpha_k:0\leqslant k\leqslant K-1\}.
            \end{equation}
            If $\ell_1 \alpha_k^\tau \leqslant \frac{\ell_2}{2} \alpha_k$, which is equivalent to $\bar{\alpha}\leqslant (\ell_2/2\ell_1)^{\frac{1}{\tau-1}}$, then
            \begin{equation}\label{PL-recursion:1/2}
                y_{k+1} \leqslant \Big(1- \frac{\ell_2}{2} \alpha_k\Big) y_k  + \ell_3 \alpha_k^\tau,\quad k\in [K-1].
            \end{equation}
            Furthermore, if $\alpha_k \leqslant \frac{2}{\ell_2}$, recursion \eqref{PL-recursion:1/2} can be analyzed using Theorem \ref{main-thm}. The requirements on $\bar{\alpha}$ can therefore be summarized as
            \begin{equation}\label{PL1/2-condition}
                \bar{\alpha} \leqslant \min\Big\{\Big(\frac{\ell_2}{2\ell_1}\Big)^{\frac{1}{\tau-1}},\frac{2}{\ell_2}\Big\}.
            \end{equation}
            However, when $\theta \neq \frac{1}{2}$, Theorem \ref{main-thm} is no longer directly applicable. This issue can be addressed by replacing the nonlinear term ``$y_k^{2\theta}$'' with a suitable linear approximation. Since the function $h_\theta(y) := y^{2\theta}$ is convex on $\mathbb{R}_+$, we can apply the inequality
            \[ y^{2\theta} \geqslant x^{2\theta} + 2\theta x^{2\theta-1}(y-x) \quad \forall~ x,y\geqslant0. \]
            to construct a linear lower bound.
            By choosing $x = u_k$ and $y = y_k$, inequality \eqref{PL-recursion} can be relaxed to
            \begin{equation}\label{eq:roadmap-step2-1}
            \begin{aligned}
            	y_{k+1} & \leqslant (1 + \ell_1 \alpha_k^\tau)  y_k  - \ell_2 \alpha_k [u_k^{2\theta}+2\theta u_k^{2\theta-1}(y_k-u_k)]+ \ell_3 \alpha_k^\tau \\
             &= [1+\ell_1\alpha_k^\tau - 2\theta \ell_2\alpha_k u_k^{2\theta-1}] y_k + [(2\theta-1)\ell_2 \alpha_k u_k^{2\theta} +\ell_3 \alpha_k^\tau],
            \end{aligned}
            \end{equation}
            where $\{u_k\}_{k}$ is a (free) non-negative auxiliary sequence. 

            \vspace{2mm}
            \noindent\textbf{Step 2: Optimal Choice of $\{u_k\}_k$}. The sequence $\{u_k\}_k$ should be selected such that the following conditions are satisfied:
            \begin{itemize}
            \item[1.] The term ``$\theta \ell_2\alpha_ku_k^{2\theta-1}$'' is not smaller than ``$\ell_1\alpha_k^\tau$'' (for all $k$). Then, \eqref{eq:roadmap-step2-1} can be further relaxed to
            \begin{equation}\label{eq:roadmap-step2-2}
                y_{k+1} \leqslant [1-\theta \ell_2\alpha_k u_k^{2\theta-1}] y_k + [(2\theta-1)\ell_2\alpha_k u_k^{2\theta}+\ell_3 \alpha_k^\tau].
            \end{equation}
            \item[2.] For \eqref{eq:roadmap-step2-2}, the expected rate $r_k \sim \frac{(2\theta-1)\ell_2 \alpha_k u_k^{2\theta}+\ell_3 \alpha_k^\tau}{\alpha_k u_k^{2\theta-1}} = (2\theta-1)\ell_2 u_k + \ell_3 \frac{\alpha_k^{\tau-1}}{u_k^{2\theta-1}}$ is optimal (as fast as possible), cf. \Cref{main-thm}. 
            \end{itemize}
            To satisfy the first requirement, it is sufficient to ensure that  
            \begin{equation}\label{roadmap:uk-lb}
            u_k\geqslant ({\ell_1}/{\theta \ell_2})^{\frac{1}{2\theta-1}} \alpha_k^{\frac{\tau-1}{2\theta-1}}.
            \end{equation}
            Meanwhile, to achieve tight estimates, $r_k$ should be minimized. Optimizing $ r_k $ as a function of $ u_k $ over $ \mathbb{R}_{++} $, we obtain the candidate $({\ell_3}/{\ell_2})^{\frac{1}{2\theta}}\alpha_k^{\frac{\tau-1}{2\theta}}$ for $u_k$. However, this choice does not necessarily satisfy the condition \eqref{roadmap:uk-lb}. To address this issue, we choose $ u_k $ to be proportional to $ \alpha_k^{\frac{\tau - 1}{2\theta}} $, i.e., $u_k := \ukc \alpha_k^{\frac{\tau - 1}{2\theta}}$ for some $\ukc > 0$. Observe that for all $\theta\in (\frac{1}{2},1]$, we have
            \[
            \Big(\frac{\ell_1}{\theta \ell_2}\Big)^{\frac{1}{2\theta-1}} \alpha_k^{\frac{\tau-1}{2\theta-1}} = \Big(\frac{\ell_1}{\theta \ell_2}\Big)^{\frac{1}{2\theta-1}} \alpha_k^{\frac{\tau-1}{2\theta}} \bar{\alpha}^{\frac{\tau-1}{2\theta(2\theta-1)}} \Big(\frac{\alpha_k}{\bar{\alpha}}\Big)^{\frac{\tau-1}{2\theta(2\theta-1)}} \leqslant \Big(\frac{\ell_1 \bar{\alpha}^{\frac{\tau-1}{2\theta}}}{\theta \ell_2}\Big)^{\frac{1}{2\theta-1}} \alpha_k^{\frac{\tau-1}{2\theta}}
            \]
            To guarantee \eqref{roadmap:uk-lb}, it suffices to ensure that $(\frac{\ell_1 \bar{\alpha}^{\frac{\tau-1}{2\theta}}}{\theta \ell_2})^{\frac{1}{2\theta-1}} \leqslant (2\theta-1)\delta$ and $\zeta\geqslant (2\theta-1)\delta$ for some $\delta>0$. Define 
            $$\mathfrak{p}(\theta) = \begin{cases}
                (2\theta-1)^{2\theta-1},\quad &\text{if $\theta\in (\frac{1}{2},1]$},\\
                1,\quad &\text{if $\theta=\frac{1}{2}$.}
            \end{cases}
            $$ 
            It is straightforward to verify that $\mathfrak{p}(\theta)\in [e^{-\frac{1}{e}},1]$ for all $\theta\in [\frac{1}{2},1]$, and that $\mathfrak{p}$ is continuous at $\frac{1}{2}$. Assume that $\bar{\alpha} \leqslant (\frac{\theta \mathfrak{p}(\theta) \ell_2 \delta^{2\theta-1}}{\ell_1})^{\frac{2\theta}{\tau - 1}}$ (recall that $ \bar{\alpha} $ is defined in equation \eqref{eq:def-maxstep}), which yields $(\frac{\ell_1 \bar{\alpha}^{\frac{\tau - 1}{2\theta}}}{\theta \ell_2})^{\frac{1}{2\theta - 1}} \leqslant (2\theta - 1)\delta$ for all $\theta\in (\frac{1}{2},1]$. Then, we can choose $ u_k $ as follows:
            \begin{equation}\label{eq:uk-def}
            u_k := \ukc \alpha_k^{\frac{\tau - 1}{2\theta}},\quad\text{where $\ukc := \max\{(2\theta - 1)\delta, ({\ell_3}/{\ell_2})^{\frac{1}{2\theta}}\}$.} 
            \end{equation}
            Here, $\delta > 0$ is a free (algorithmic) parameter that will be specified later. By construction, this choice of $ u_k $ satisfies \eqref{roadmap:uk-lb}. Moreover, the condition $\ukc \geqslant ({\ell_3}/{\ell_2})^{\frac{1}{2\theta}}$ ensures that $\ell_3 \leqslant \ell_2 \zeta^{2\theta}$. Substituting $u_k$ into inequality \eqref{eq:roadmap-step2-2} and using the bound $\ell_3 \leqslant \ell_2 \zeta^{2\theta}$, we obtain: 
                \begin{equation}\label{eq:relaxed-PL-final0}
                    y_{k+1} \leqslant (1 - \theta \ell_2 \ukc^{2\theta - 1} \alpha_k^{1 / \varrho}) y_k + 2\theta \ell_2 \ukc^{2\theta} \alpha_k^\tau,
                \end{equation}
                where 
                \[ \varrho := \Big(1+(2\theta-1)\cdot \frac{\tau-1}{2\theta}\Big)^{-1} = \frac{2\theta}{(2\theta - 1)\tau + 1}.
                \]
                Denoting 
                \[
                \xi := \theta \ell_2 \ukc^{2\theta-1} = \max\{\theta\mathfrak{p}(\theta) \ell_2 \delta^{2\theta-1},  \theta \ell_2^{\frac{1}{2\theta}}\ell_3^{\frac{2\theta-1}{2\theta}}\},
                \]
                inequality \eqref{eq:relaxed-PL-final0} can be rewritten as
                \begin{equation}\label{eq:relaxed-PL-final}
                    y_{k+1} \leqslant (1 - \xi \alpha_k^{1 / \varrho}) y_k + 2\zeta\xi \alpha_k^\tau,
                \end{equation}
                To apply Theorem~\ref{main-thm} to inequality~\eqref{eq:relaxed-PL-final}, it is additionally required that $\xi\bar{\alpha}^{1/\varrho} \leqslant 1$, i.e., $\bar{\alpha}\leqslant \xi^{-\varrho}$. 
                Therefore, the overall conditions on $ \bar{\alpha} $ are given by
                \begin{equation}\label{eq:relaxed-PL-condition}
                    \bar{\alpha} \leqslant \min\Big\{\Big(\frac{\theta \mathfrak{p}(\theta)\ell_2\delta^{2\theta-1}}{\ell_1}\Big)^{\frac{2\theta}{\tau - 1}},\xi^{-\varrho}\Big\}.
                \end{equation}
                Notably, the inequalities~\eqref{eq:relaxed-PL-final} and \eqref{eq:relaxed-PL-condition} reduce to~\eqref{PL-recursion:1/2} and~\eqref{PL1/2-condition}, respectively, when $\theta=\frac{1}{2}$. This unifies the analysis across the entire range $\theta\in [\frac{1}{2},1]$, which is not a coincidence, but rather a result of the deliberate design of $\{u_k\}_k$. The simplified and relaxed recursion \eqref{eq:relaxed-PL-final} is the basis of our final step and of the application of the generalized Chung's Lemma. 

            \vspace{2mm}
            \noindent\textbf{Step 3: Application of Chung's Lemma to \eqref{eq:relaxed-PL-final}}. Suppose that the step sizes $\alpha_k$ can be represented as a function $\eta : \R \to \R_+$ evaluated at $b_k$, i.e., $\alpha_k = \eta(b_k)$. For simplicity, we assume that $\eta$ is differentiable. Setting $s(x) = \frac{1}{\xi}  \eta(x)^{-1/\varrho}$, $t(x) = \frac{1}{2\ukc\xi} \eta(x)^{-\tau}$, we then have
            \begin{equation}\label{roadmap_step3_notation-choice-2}
            r(x) = 2\ukc \eta(x)^{\frac{\tau-1}{2\theta}},~~u(x) = \frac{\tau-1}{2\theta\xi} \eta(x)^{-(1+\frac{1}{\varrho})} \eta'(x).
            \end{equation}
            The rest of the derivation involves determining the appropriate function $\eta$ and sequence $\{b_k\}_k$ to use for different types of step sizes and verifying the conditions stated in \Cref{main-thm} for such $\eta$ and $\{b_k\}_k$.

    \begin{rem}
        A straightforward approach to ensure condition \eqref{roadmap:uk-lb} is to directly assume that $({\ell_1 \bar{\alpha}^{\frac{\tau - 1}{2\theta}}}/{\theta \ell_2})^{\frac{1}{2\theta - 1}} \leqslant ({\ell_3}/{\ell_2})^{\frac{1}{2\theta}}$ and to choose $u_k = ({\ell_3}/{\ell_2})^{\frac{1}{2\theta}} \alpha_k^{\frac{\tau-1}{2\theta}}$. However, this leads to an unnecessarily stringent constraint on $\bar{\alpha}$. In particular, under this strategy, the upper bound on $\bar{\alpha}$ tends to zero as the noise term $\sigma^2$ approaches zero (as $\ell_3$ depends on $\sigma^2$). This is undesirable, as smaller noise should allow the use of a larger range of step sizes.
    \end{rem}

\subsection{Notations and Important Constants}
This subsection introduces the notations and key constants that will be used consistently throughout the remainder of the paper. We use $K$ to represent the total iterations in $\SGD$ or total epochs in $\RR$. For all $k\in [K]$, we denote
\[
  y_k := \Exp[f(\vx^k) - f^*].
\]
We introduce two sets of notations for analyzing the performance of $\SGD$ (see \Cref{algo:sgd}) and $\RR$ (see \Cref{algo:rr}). For $\SGD$, we have $\tau = 2$ and choosing $\delta = 1$ in \eqref{eq:uk-def}, we define: 
\begin{equation}\label{eq:notation-sgd}
\varrho_s = \frac{2\theta}{4\theta-1},\quad\omega_s = \frac{1}{4\theta-1},\quad \zeta_s = \max\Big\{2\theta - 1, \Big(\frac{\sL \sigma^2}{2\mu}\Big)^{\frac{1}{2\theta}}\Big\},\quad \xi_s = \theta \mu \zeta_s^{2\theta-1}.
\end{equation}
For $\RR$, recalling $\tau = 3$ and setting $\delta = N^{-\frac{1}{2\theta}}$ in \eqref{eq:uk-def}, we define:
\begin{equation}\label{eq:notation-rr}
\begin{gathered}
    \varrho_r = \frac{\theta}{3\theta-1},\quad \omega_r = \frac{1}{3\theta-1},\\
    \zetarc = \max\Big\{2\theta - 1, \Big(\frac{\sL^2 \sigma^2}{\mu}\Big)^{\frac{1}{2\theta}}\Big\},\quad \zeta_r = N^{-\frac{1}{2\theta}}\zetarc,\quad\xirc = \frac{\theta \mu \zetarc^{2\theta-1}}{2},\quad \xi_r = \frac{\xirc}{N^{1-\frac{1}{2\theta}}}.
\end{gathered}
\end{equation}
The explicit expressions for $\xi_s$ and $\bar{\xi}_r$ are given by:
\[
\xi_s = \max\Big\{\theta\mathfrak{p}(\theta) \mu,  \theta \mu^{\frac{1}{2\theta}}\Big(\frac{\sL\sigma^2}{2}\Big)^{\frac{2\theta-1}{2\theta}}\Big\},\quad \bar{\xi}_r = \max\Big\{\frac{\theta\mathfrak{p}(\theta)\mu}{2},  \theta \Big(\frac{\mu}{2}\Big)^{\frac{1}{2\theta}}\Big(\frac{\sL^2\sigma^2}{2}\Big)^{\frac{2\theta-1}{2\theta}}\Big\},
\]
where we adopt the convention $0^0 = 1$ when $\theta = \frac{1}{2}$. To characterize the specific range of step sizes typically required by non-asymptotic convergence results, we define the following constants for $\SGD$ and $\RR$, respectively:
\begin{equation}\label{eq:stepsize-ub}
   \bar{\alpha}_s = \min\Big\{\Big(\frac{2\theta \mathfrak{p}(\theta)\mu}{\sA\sL}\Big)^{\frac{2\theta}{\tau-1}}, \xi_s^{-\varrho_s},\frac{1}{\sL}\Big\},~~  \bar{\alpha}_r = \min\Big\{\Big(\frac{\theta \mathfrak{p}(\theta)\mu N^{\frac{1}{2\theta}}}{\sA\sL^2}\Big)^{\frac{2\theta}{\tau-1}}, \Big(\frac{N^{1-\frac{1}{2\theta}}}{\xirc}\Big)^{\varrho_r},\frac{1}{2\sL}\Big\}.
\end{equation}
Our main results require $\alpha_k \leqslant \bar{\alpha}_s$ for $\SGD$ and $\alpha_k \leqslant \bar{\alpha}_r$ for $\RR$ (cf.\ \eqref{eq:relaxed-PL-condition} and \Cref{lem:descent-property}). 

\subsection{The Analysis of Exponential and  Cosine Step Sizes}\label{sec:exp&cos}
\noindent\textbf{Exponential and Cosine Step Size.}  We will use the following formats for exponential step sizes:
\begin{equation*}\label{eq:exp_step}
    \alpha_k = \alpha\gamma^{k},\quad\text{where} \quad \gamma = \Big(\frac{\beta}{K}\Big)^{\frac{p}{K}}, \quad \alpha, p>0, \quad K>\beta>0,
    \end{equation*}
and for cosine step sizes:
\begin{equation*}\label{eq:cos_step}
    \alpha_k = \alpha\Big[{\frac{1+\cos(\frac{k\pi}{K})}{2}}\Big]^p,\quad\text{where} \quad \alpha, p>0.
\end{equation*}

\vspace{1mm}
\noindent\textbf{The Splitting Technique.} A direct application of \Cref{main-thm} fails for exponential and cosine step sizes. To illustrate this, consider $\SGD$ with exponential step sizes, setting $\theta = \frac{1}{2}$, $\beta = 1$, and $p=1$ for simplicity. Following the previous roadmap, we only need to analyze recursion \eqref{eq:relaxed-PL-final}. Focusing on \eqref{eq:relaxed-PL-final} and applying \eqref{roadmap_step3_notation-choice-2} with $b_k = k$ and $\eta(x) = \alpha K^{-\frac{x}{K}}$, we obtain $\varrho_s = 1$ and 
\[
(b_{k+1}-b_k) u(b_k) = \frac{1}{\xi_s} \eta(k)^{-2} \eta^\prime(k)= -\frac{1}{\xi_s \alpha} K^{-(1-\frac{k}{K})} {\log(K)}.
\]
Recall that \Cref{main-thm} (b) requires $(b_{k+1}-b_k) u(b_k) \geqslant -1 + 1/\lambda_{k+1}$ for some $\lambda_{k+1}>0$. As $k$ approaches $K$ (e.g., $k = K-1$), the left-hand side scales as $-\Theta(\log(K))$, which violates such lower bound for large $K$. Therefore, the condition cannot be satisfied for the entire sequence. However, it holds with $\lambda_k \equiv {2}$ for the partial index set $k \leqslant \overline{K}$, where $\overline{K}:= \lfloor \log(\frac{\xi_s \alpha K}{2\log(K)})\frac{K}{\log(K)} \rfloor$. This allows us to apply \Cref{main-thm} to estimate the intermediate iterate $y_{\overline{K}+1}$. Crucially, using the extension lemma (\Cref{lem:extension}), we can propagate this estimate to the final iterate $y_K$, establishing the desired convergence result without requiring the condition to hold globally.

\vspace{1mm}

\noindent\textbf{Results on Exponential and Cosine Step Sizes.} Following the general roadmap and using the above splitting technique, we are now in the position to show convergence rates for both $\SGD$ and $\RR$ under the general $(\theta,\mu)$-PL condition when using exponential and cosine step sizes. We first provide an important technical lemma.

\begin{lem}\label[lemma]{lem:exp&cos0}
                Consider a non-negative sequence $\{y_k\}_{0\leqslant k\leqslant K}$ satisfying  
                \begin{equation}\label{recursion-exp&cos0} 
            		y_{k+1} \leqslant (1 -c \alpha_k^\nu)  y_k  + d \alpha_k^{\nu+q}
                \end{equation}
            	for all $k\in [K-1]$, where $c,d,\nu,q>0$, $K\in \mathbb{N}_{+}$.
\begin{enumerate}[label=\textup{\textrm{(\alph*)}},leftmargin=8ex,topsep=0pt,itemsep=0ex,partopsep=0ex]
\item Let $\alpha_k = \alpha\gamma^{k}$ be given with $\gamma = (\frac{\beta}{K})^{\frac{p}{K}}$, $K>\beta>0$, $p>0$, and $0 < \alpha \leqslant c^{-\frac{1}{\nu}}$; then
              \[
                y_K \leqslant \frac{2d}{c}\max\cbra{ \bra{\frac{2pq\log(K/\beta)}{cK}}^{\frac{q}{\nu}},{\alpha^q}\prt{\frac{\beta}{K}}^{pq}}+ y_0 \exp\prt{-\frac{c\alpha^\nu}{p\nu}\frac{[1-(\beta/K)^{p\nu}] K}{\log(K/\beta)}}.
                \]
\item Let $\alpha_k = \alpha [\frac{1+\cos(\frac{k\pi}{K})}{2}]^p$ be given with $\alpha,p>0$ and $\alpha \leqslant c^{-\frac{1}{\nu}}$; then
                \[
                y_{K}\leqslant \frac{d}{c} \max\cbra{\Big[\frac{2D}{cK}\Big]^{\frac{q}{\nu}}, 2\prt{\frac{\pi^2}{4}}^{p q} \alpha^{\frac{q}{2p\nu +1}}\bra{\frac{2D}{cK}}^{\frac{2p q}{2p\nu +1}}} + y_0 \exp\Big(-\frac{c\alpha^\nu}{2^{\max\{1,p\nu\}}} K \Big),
                \]
                where $D = \max\{1,2p q \pi^2\}$.
\end{enumerate}
           \end{lem}

We continue with a detailed proof of \Cref{lem:exp&cos0}. Upon first reading of the manuscript, the reader may safely skip to \Cref{lem:exp&cos}. \vspace{0.5ex}

           \begin{proof}
               We first prove part (a). Setting $b_k = k$, $s(x) = \frac{1}{c\alpha^\nu} \gamma^{-\nu x}$, $t(x) = \frac{1}{d\alpha^{\nu+q}} \gamma^{-(\nu+q)x}$ (and $\eta(x) = \alpha\gamma^x$), we have 
               \[
               r(x) = \frac{d\alpha^q}{c}\gamma^{qx},\quad u(x) = (\log(r))'(x) s(x) = \frac{q\log(\gamma)}{c\alpha^\nu} \gamma^{-\nu x}.
               \]
               Due to $0 < \gamma < 1$, $r$ is convex and differentiable on $\R$, which meets the requirements on $r$ in \Cref{main-thm}. In addition, $u$ is also a decreasing function on $\R$. Next, choosing $\lambda_k \equiv 2$, condition \eqref{condition2} in \Cref{main-thm}~(b) is satisfied if
               \begin{equation}\label{exp-transition-point}
                   u(x)\geqslant -\frac{1}{2} \quad \iff \quad 
                   \gamma^{-x}\leqslant \alpha\bra{\frac{c}{2q\log(\gamma^{-1})}} = \alpha \bra{\frac{cK}{2pq\log(K/\beta)}}^{\frac{1}{\nu}} \quad \iff \quad  x \leqslant \bar{x},
               \end{equation}
               where $\bar{x}:=\log({\frac{c\alpha^\nu K}{2pq\log(K/\beta)}})\frac{K}{\nu p \log(K/\beta)}$. \\
               \noindent \textbf{Case I}. If $\frac{c\alpha^\nu K}{2pq\log(K/\beta)}<1$, i.e., $\alpha < [{\frac{2pq\log(K/\beta)}{cK}}]^{\frac{1}{\nu}}$, it holds that $\bar{x}<0$. Applying \eqref{recursion-exp&cos0} with $k=0$ and $s(b_0) \geqslant 1$, it follows
                \[
                y_1 \leqslant \frac{1}{t(b_0)} + y_0 \prt{1-\frac{1}{s(b_0)}} \leqslant r(b_0) + y_0 \prt{1-\frac{1}{s(b_0)}}.
                \]
                Using \Cref{lem:extension} and setting $B = r(b_0)$, $C = y_0$, $k_0 = \overline{K} = 0$, we obtain
                \[
                \begin{aligned}
                    y_K&\leqslant r(b_0) + y_0 {\prod}_{k=0}^{K-1}\prt{1-\frac{1}{s(b_k)}} = \frac{d\alpha^q}{c} + y_0 {\prod}_{k=0}^{K-1}({1-c\alpha^\nu \gamma^{\nu k}}) \\
                    &\leqslant \frac{d}{c}\Big[\frac{2pq\log(K/\beta)}{cK}\Big]^{\frac{q}{\nu}} + y_0 {\prod}_{k=0}^{K-1}({1-c\alpha^\nu \gamma^{\nu k}}).
                \end{aligned}
                \]
                \noindent \textbf{Case II}. If $\frac{c\alpha^\nu K}{2pq\log(K/\beta)}\geqslant 1$, we can infer $\bar{x}\geqslant 0$. Setting $\overline{K} := \min\{\lfloor\bar{x}\rfloor,K-1\}$ and according to \Cref{main-thm} and \Cref{lem:extension} (with $\lambda_k \equiv 2$), we have
                \[
                    y_{K}\leqslant 2 r(b_{\overline{K}+1}) + y_0 {\prod}_{k=0}^{K-1}({1-c\alpha^\nu \gamma^{\nu k}}) = 2r(\overline{K}+1)+ y_0 {\prod}_{k=0}^{K-1}({1-c\alpha^\nu \gamma^{\nu k}}).
                \]
                It is evident that $r(\overline{K}+1)\leqslant \max\{r(\lfloor\bar{x}\rfloor+1),r(K)\}$. Furthermore, due to $\lfloor\bar{x}\rfloor+1>\bar{x}$ and \eqref{exp-transition-point}, we have $\gamma^{-(\lfloor\bar{x}\rfloor+1)}>\alpha [\frac{cK}{2pq\log(K/\beta)}]^{{1}/{\nu}}$ which, by the definition of $r$, implies
                \[
                r(\lfloor\bar{x}\rfloor+1)\leqslant \frac{d}{c} \Big[\frac{2pq\log(K/\beta)}{cK}\Big]^{\frac{q}{\nu}}.
                \]
                In addition, it holds that $r(K) = \frac{d\alpha^q}{c}({\frac{\beta}{K}})^{pq}$. Consequently, combining these estimates and the derivations in \textbf{Case I} and \textbf{Case II}, we obtain the following
                \[
                y_K \leqslant \frac{{2}d}{c}\max\cbra{ \bra{\frac{2pq\log(K/\beta)}{cK}}^{\frac{q}{\nu}},{\alpha^q}\prt{\frac{\beta}{K}}^{pq}} + y_0 {\prod}_{k=0}^{K-1}({1-c\alpha^\nu \gamma^{\nu k}}).
                \]
                We continue to estimate the product term. Using $1-x \leqslant \exp(-x)$ and $\gamma = (\frac{\beta}{K})^{\frac{p}{K}}$, we have
                \[
                {\prod}_{k = 0}^{K-1} ({1-c\alpha^\nu \gamma^{\nu k}}) \leqslant \exp\prt{-c\alpha^\nu {\sum}_{k = 0}^{K-1} \gamma^{\nu k}} = \exp\prt{-c\alpha^\nu \frac{1-(\beta/K)^{p\nu}}{1-\gamma^\nu}}. 
                \]
                Due to $\log(\gamma^\nu)\leqslant \gamma^{\nu}-1$, it follows $1-\gamma^\nu \leqslant \log(\gamma^{-\nu}) = \frac{p\nu\log(K/\beta)}{K}$ and we can infer
                \[
                {\prod}_{k = 0}^{K-1}({1-c\alpha^\nu \gamma^{\nu k}}) \leqslant \exp\prt{-\frac{c\alpha^\nu}{p\nu}[1-(\beta/K)^{p\nu}]\frac{K}{\log(K/\beta)}}
                \]
                and
                \[
                y_K \leqslant \frac{{2}d}{c}\max\cbra{ \bra{\frac{2pq\log(K/\beta)}{cK}}^{\frac{q}{\nu}},{\alpha^q}\prt{\frac{\beta}{K}}^{pq}} + y_0 \exp\prt{-\frac{c\alpha^\nu}{p\nu}\frac{[1-(\beta/K)^{p\nu}]K}{\log(K/\beta)}}.
                \]
           Next, we present the proof of part (b). Setting $b_k = [{\frac{1+\cos(\frac{k\pi}{K})}{2}}]^{p q}$, $s(x) = \frac{1}{c\alpha^\nu} x^{-\frac{\nu}{q}}$, $t(x) = \frac{1}{d\alpha^{\nu+q}} x^{-(1+\frac{\nu}{q})}$ (and $\eta(x) = \alpha x^{\frac{1}{q}}$), it holds that
               \[
               r(x) = \frac{d\alpha^q}{c}x, \quad u(x) = (\log(r))'(x) s(x) = \frac{1}{c\alpha^\nu} x^{-(1+\frac{\nu}{q})}
               \]
               and $u(b_k) = \frac{1}{c\alpha^\nu}  [{\frac{1+\cos(\frac{k\pi}{K})}{2}}]^{-p (q+\nu)}$.
               For $K\geqslant 1$, $k = K-1$, using \Cref{lem:Tech-ineqs} (c), we have 
                \[
                (b_{k+1}-b_k) u(b_k) = -b_ku(b_k) = -\frac{1}{c\alpha^\nu} \bra{\frac{1+\cos(\frac{k\pi}{K})}{2}}^{-p \nu} \geqslant -\frac{1}{c\alpha^\nu} \prt{1-\frac{k}{K}}^{-2p \nu}. 
                \]
               For $K\geqslant 2$ and $k\leqslant K-2$, applying \Cref{lem:Tech-ineqs} (b) with $x = \frac{1+\cos(\frac{(k+1)\pi}{K})}{2}$, $y = \frac{1+\cos(\frac{k\pi}{K})}{2}$, and $r = p q$, it follows that
               \[
               b_{k+1}-b_k \geqslant p q \bra{\frac{1+\cos(\frac{k\pi}{K})}{2}}^{p q} \bra{\frac{1+\cos(\frac{(k+1)\pi}{K})}{2}}^{-1} \bra{\frac{\cos(\frac{(k+1)\pi}{K})-\cos(\frac{k\pi}{K})}{2}}
               \]
               and we can infer
               \[
               (b_{k+1}-b_k)u(b_k) \geqslant \frac{p q}{c\alpha^\nu} \bra{\frac{1+\cos(\frac{k\pi}{K})}{2}}^{-p \nu} \bra{\frac{1+\cos(\frac{(k+1)\pi}{K})}{2}}^{-1} \bra{\frac{\cos(\frac{(k+1)\pi}{K})-\cos(\frac{k\pi}{K})}{2}}.
               \]
               Next, using \Cref{lem:Tech-ineqs} (c)--(e), we have for any $0\leq k \leqslant K-2$:
                \[
                (b_{k+1}-b_k)u(b_k) 
                \geqslant \frac{4p q}{c\alpha^\nu} \prt{1-\frac{k}{K}}^{-2p\nu} \prt{1-\frac{k}{K}}^{-2} \bra{-\frac{\pi^2}{2K}\prt{1-\frac{k}{K}}}
                = -\frac{2p q \pi^2}{c\alpha^\nu K} \prt{1-\frac{k}{K}}^{-(2p\nu+1)}.
                \]
                For $k = K-1$, noticing $1-\frac{k}{K} = \frac{1}{K}$, we can combine the previous estimates. In particular, for all $k\in [K-1]$, $K\geqslant 1$, it holds that
                \begin{equation}\label{eq:key-quantity-cos}
                    (b_{k+1}-b_k) u(b_k) \geqslant -\frac{D}{c\alpha^\nu K} \prt{1-\frac{k}{K}}^{-(2p\nu+1)},
                \end{equation}
                where $D = \max\{1,2p q \pi^2\}$. As before and based on \eqref{eq:key-quantity-cos}, we now discuss when condition \eqref{condition2} in \Cref{main-thm}~(b) holds with $\lambda_k \equiv 2$. \\[2mm]
                \noindent \textbf{Case I}. If $K\leqslant \frac{2D}{c\alpha^\nu}$, then we have $\alpha \leqslant (\frac{2D}{cK})^{\frac{1}{\nu}}$ and using \eqref{recursion-exp&cos0} and $s(b_0)\geqslant 1$, it follows that
                \[
                y_1 \leqslant \frac{1}{t(b_0)} + y_0 \Big(1-\frac{1}{s(b_0)}\Big) \leqslant r(b_0) + y_0 \Big(1-\frac{1}{s(b_0)}\Big).
                \]
                Clearly, by the definition of $b_k$ and $r$, it holds that $r(b_k) \leqslant r(b_0)$ for all $k\geqslant 1$. Applying \Cref{lem:extension} with $B = r(b_0)$, $C = y_0$, $k_0 = \overline{K} = 0$, we obtain
                \[
                \begin{aligned}
                    y_K&\leqslant r(b_0) + y_0 {\prod}_{k=0}^{K-1}\Big(1-\frac{1}{s(b_k)}\Big) = \frac{d\alpha^q}{c} + y_0 {\prod}_{k=0}^{K-1}(1-c\alpha^\nu b_k^{\frac{\nu}{q}})\\
                    &\leqslant \frac{d}{c}\left[\frac{2D}{cK}\right]^{\frac{q}{\nu}} + y_0 {\prod}_{k=0}^{K-1}(1-c\alpha^\nu b_k^{\frac{\nu}{q}}).
                \end{aligned}
                \]
                \noindent \textbf{Case II}. In the case $K> \frac{2D}{c\alpha^\nu}$, based on \eqref{eq:key-quantity-cos}, we can infer $(b_{k+1}-b_k) u(b_k) \geqslant -\frac{1}{2}$ if
                \[
                 k\leqslant \bar{x}:= K\prt{1-\bra{\frac{2D}{c\alpha^{\nu}K}}^{\frac{1}{2p\nu +1}}}.
                \]
                Choosing $\overline{K} = \lfloor \bar{x} \rfloor$, we obviously have $0<\bar{x}< K$ and $0\leqslant \overline{K}\leqslant K-1$. According to \Cref{main-thm} and \Cref{lem:extension}, this then yields
                \[
                    y_{K}\leqslant 2 r(b_{\overline{K}+1}) + y_0 \prod_{k=0}^{K-1}({1-c\alpha^\nu b_k^{\frac{\nu}{q}}}) = \frac{2d\alpha^q}{c}\bra{\frac{1+\cos(\frac{(\overline{K}+1)\pi}{K})}{2}}^{p q}+ y_0 \prod_{k=0}^{K-1}(1-c\alpha^\nu b_k^{\frac{\nu}{q}}).
                \]
                Using $\overline{K}+1>\bar{x} = K(1-[\frac{2D}{c\alpha^{\nu}K}]^{\frac{1}{2p\nu +1}})$ and Lemma \ref{lem:Tech-ineqs} (c), we further obtain
				\[	1+\cos\prt{\frac{(\overline{K}+1)\pi}{K}} \leqslant \frac{\pi^2}{2} \prt{1-\frac{\overline{K} + 1}{K}}^2 \leqslant \frac{\pi^2}{2} \bra{\frac{2D}{c\alpha^{\nu}K}}^{\frac{2}{2p\nu +1}}
				\]
                which implies $y_K \leqslant \frac{2d}{c}({\frac{\pi^2}{4}})^{p q} \alpha^{\frac{q}{2p\nu +1}}[{\frac{2D}{cK}}]^{\frac{2p q}{2p\nu +1}} + y_0 \prod_{k=0}^{K-1}({1-c\alpha^\nu b_k^{\frac{\nu}{q}}})$.
                Thus, combining these two cases, it follows that
                \[
                y_{K}\leqslant \frac{d}{c} \max\Big\{2\prt{\frac{\pi^2}{4}}^{p q} \alpha^{\frac{q}{2p\nu +1}}\bra{\frac{2D}{cK}}^{\frac{2p q}{2p\nu +1}},\Big[\frac{2D}{cK}\Big]^{\frac{q}{\nu}}\Big\} + y_0 \prod_{k=0}^{K-1}(1-c\alpha^\nu b_k^{\frac{\nu}{q}}).
                \]
                To handle the product term, we notice that 
				\[
				{\prod}_{k = 0}^{K-1}
				\Big(1-c\alpha^\nu \Big[\frac{1+\cos(\frac{k\pi}{K})}{2}\Big]^{p\nu}\Big) \leqslant \exp\Big(-c\alpha^\nu {\sum}_{k=0}^{K-1}\Big[\frac{1+\cos(\frac{k\pi}{K})}{2}\Big]^{p\nu}\Big).
				\]
	Using \Cref{lem:Tech-ineqs} (f) with $r = p \nu$, we further have $\sum_{k=0}^{K-1}[\frac{1+\cos(\frac{k\pi}{K})}{2}]^{p\nu} \geqslant \frac{K}{2^{\max\{1,p\nu\}}}$. This establishes and verifies the bound in \Cref{lem:exp&cos0} (b) and concludes the proof. 
           \end{proof}

By combining Lemma~\ref{lem:exp&cos0} with the roadmap presented in Section~\ref{sec:recursion&roadmap}, we obtain the following result.
\begin{lem}\label[lemma]{lem:exp&cos}
                Consider a non-negative sequence $\{y_k\}_{0\leqslant k\leqslant K}$ satisfying  
            	\begin{equation} 
            		y_{k+1} \leqslant (1 + \ell_1 \alpha_k^\tau)  y_k  - \ell_2 \alpha_k y_k^{2\theta}+ \ell_3 \alpha_k^\tau
            	\end{equation}
            	for all $k\in [K-1]$, where $\ell_1\geqslant 0$, $\ell_2,\ell_3>0$, $K\in \mathbb{N}$ are given. Let us define 
                \[ \varrho := \frac{2\theta}{(2\theta-1)\tau+1}, \quad \omega:= \frac{\tau-1}{(2\theta-1)\tau+1}, \quad q:= \frac{\omega}{\varrho},\]
                \[
                \zeta := \max\{(2\theta-1)\delta,({\ell_3}/{\ell_2})^{\frac{1}{2\theta}}\}, \quad \xi := \theta \ell_2 \ukc^{2\theta-1} = \max\{\theta\mathfrak{p}(\theta) \ell_2 \delta^{2\theta-1},  \theta \ell_2^{\frac{1}{2\theta}}\ell_3^{\frac{2\theta-1}{2\theta}}\}
                \]
                where $\delta>0$. Then, for any $\alpha$ satisfying $0<\alpha  \leqslant \min\{(\frac{\theta \mathfrak{p}(\theta)\ell_2 \delta^{2\theta-1}}{\ell_1})^{\frac{2\theta}{\tau - 1}},\xi^{-\varrho}\},$
                 the following statements hold:
            	\begin{enumerate}[label=\textup{\textrm{(\alph*)}},leftmargin=8ex,topsep=0pt,itemsep=0ex,partopsep=0ex]
                \vspace{2mm}
                    \item Let $\alpha_k = \alpha\gamma^{k}$ be given with $\gamma = (\frac{\beta}{K})^{\frac{p}{K}}$, $p>0$, and $K > \beta > 0$. Then, it holds that
                    \vspace{-2mm}
                    \[
                    y_K\leqslant 4\max\cbra{\frac{\zeta^\varrho}{\ell_2^\omega}\Big[\frac{2pq\log({K}/{\beta})}{\theta K}\Big]^{\omega},\zeta\alpha^{\frac{\omega}{\varrho}}\prt{\frac{\beta}{K}}^{pq}} + y_0 \exp\Big(- \frac{\varrho \xi\alpha^{\frac{1}{\varrho}}}{p} \frac{[1-(\beta/K)^{p/\varrho}]K}{\log({K}/{\beta})}\Big).
                    \]
                    \item Let $\alpha_k = \alpha[\frac{1+\cos(\frac{k\pi}{K})}{2}]^p$ be given with $p>0$. Then, we have
                    \[
                    \begin{aligned}
                    y_K \leqslant  \max\cbra{\frac{2\zeta^\varrho}{\ell_2^\omega}\bra{\frac{2D}{\theta K}}^{\omega},4\Big(\frac{\pi^2}{4}\Big)^{pq}\bra{\frac{\zeta^{2p\varrho+\varrho}\alpha^\omega}{\ell_2^{2p\omega}}}^{\frac{1}{2p+\varrho}}\bra{\frac{2D}{\theta K}}^{\frac{2p\omega}{2p+\varrho}}} + y_0 \exp\prt{-\frac{\xi\alpha^{\frac{1}{\varrho}}K}{2^{\max\{1,p/\varrho\}}}},
                    \end{aligned}
                    \]
                    where $D = \max\{1,2pq\pi^2\}$.
            	\end{enumerate}
            \end{lem}
            \begin{proof}
                The results follow directly from the roadmap provided in Section~\ref{sec:recursion&roadmap}. More specifically, based on \eqref{eq:relaxed-PL-final} and \eqref{eq:relaxed-PL-condition}, we can set $c = \xi$, $d = 2 \zeta \xi$, $\nu = \frac{1}{\varrho}$, and $q = \tau - \frac{1}{\varrho} = \frac{\tau - 1}{2\theta} = \frac{\omega}{\varrho}$ in \Cref{lem:exp&cos0}. Here, we also use the identity $\zeta/\xi^\omega = \zeta^{\varrho}/(\theta\ell_2)^\omega$. 
            \end{proof}

The following key results on the convergence behavior of $\SGD$ and $\RR$ with exponential and cosine step sizes can now be obtained directly from Lemma~\ref{lem:exp&cos}.

\begin{thm}[Rate of $\SGD$ with Exponential Step Sizes]\label{thm:sgd_exp}
	Assume that \ref{A1}--\ref{A2} and \ref{P} hold. Let $\{\vx^k\}_{0 \leq k \leqslant K}$, with $K \geqslant 1$ and $K>\beta>0$, be generated by $\SGD$ using the step size $\alpha_k = \alpha \gamma^{k}$, where $\gamma = (\frac{\beta}{K})^{\frac{p}{K}}$. Define $\varrho_s$, $\omega_s$, $\zeta_s$, $\xi_s$ and $\bar{\alpha}_s$ as in equations \eqref{eq:notation-sgd} and \eqref{eq:stepsize-ub}. If $0<\alpha  \leqslant \bar{\alpha}_s$, then the following result holds:
    \vspace{-1mm}
	\begin{equation*}\label{eq:sgd_exp1}
		y_K\leqslant 4\max\cbra{\frac{\zeta_s^{\varrho_s}}{\mu^{\omega_s}}\Big[\frac{p\log({K}/{\beta})}{\theta^2 K}\Big]^{\omega_s},\zeta_s \alpha^{\frac{1}{2\theta}}\prt{\frac{\beta}{K}}^{\frac{p}{2\theta}}} + y_0 \exp\Big(- \frac{\varrho_s \xi_s\alpha^{\frac{1}{\varrho_s}}}{p}\cdot \frac{[1-(\beta/K)^{p/\varrho_s}]K}{\log({K}/{\beta})}\Big).
	\end{equation*}
    Moreover, if $p\geqslant \varrho_s$, then
    \begin{equation}\label{eq:sgd_exp2}
    y_K = \cO\prt{\Big[\frac{\log({K}/{\beta})}{K}\Big]^{\omega_s}}.
    \end{equation}
\end{thm}

The full proof is detailed in \Cref{proof:exp}. For the specific case $\theta = \frac{1}{2}$, we obtain a rate of $\cO({\log(K)}/{K})$. This is sharper than the rate $\cO({\log^2(K)}/{K})$ established in \cite[Theorem 1]{li2021second} and \citep{vaswani2022towards}, with the latter one assuming strong convexity. In \cite[Theorem 4.1]{wang2021convergence}, a matching bound is derived under strong convexity and bounded variance assumptions. Furthermore, our result covers the broader range $\theta \in (\frac12, 1]$, which is not addressed in  previous works. 

\begin{thm}[Rate of $\RR$ with Exponential Step Sizes]\label{thm:rr_exp}
	Assume that \ref{B1}--\ref{B3} and \ref{P} hold. Let $\{\vx^k\}_{0 \leq k \leqslant K}$, with $K \geqslant 1$ and $K>\beta>0$, be generated by $\RR$ using the step size $\alpha_k = \alpha \gamma^{k}$, where $\gamma = (\frac{\beta}{K})^{\frac{p}{K}}$. Define $\varrho_r$, $\omega_r$, $\zetarc$, $\xirc$, and $\bar\alpha_r$ as in equations \eqref{eq:notation-rr} and \eqref{eq:stepsize-ub}. If $0<\alpha  \leqslant \bar{\alpha}_r$, then it holds that:
	\begin{equation*}\label{eq:rr_exp1}
		y_K\leqslant 4\max\cbra{\frac{\zetarc^{\varrho_r}}{({\mu}/{2})^{\omega_r}}\Big[\frac{2p\log({K}/{\beta})}{\theta^2 \sqrt{N}K}\Big]^{\omega_r},\frac{\zetarc\alpha^{\frac{1}{\theta}}}{N^{\frac{1}{2\theta}}}\prt{\frac{\beta}{K}}^{\frac{p}{\theta}}} + y_0 \exp\Big(- \frac{\varrho_r \xirc\alpha^{\frac{1}{\varrho_r}}}{p N^{1-\frac{1}{2\theta}}} \frac{[1-(\beta/K)^{p/\varrho_r}]K}{\log({K}/{\beta})}\Big),
	\end{equation*}
    Furthermore, if $p\geqslant \varrho_r$, $\frac{K}{\log(K/\beta)} \geqslant \frac{2p\log(\sqrt{N}K)}{\theta\xirc \alpha^{{1}/{\varrho_r}}} N^{1-\frac{1}{2\theta}}$, and $K \geqslant 2\beta > 0$, then we have
        \begin{equation*}\label{eq:rr_exp3}
		y_K = \cO\prt{\bra{\frac{\log({K}/{\beta})}{\sqrt{N}K}}^{\omega_r}}.
	\end{equation*}
\end{thm}  

The proof of Theorem~\ref{thm:rr_exp} is provided in \Cref{proof:exp}. To the best of our knowledge, no existing literature discusses $\RR$ with exponential step sizes, not even for the standard case $\theta = \frac12$. \Cref{thm:rr_exp} fills this gap by providing a complete convergence characterization for all $\theta \in [\frac{1}{2}, 1]$ under the $(\theta,\mu)$-PL condition.

\begin{thm}[Rate of $\SGD$ with Cosine Step Sizes]\label{thm:sgd_cos}
	Let the conditions \ref{A1}--\ref{A2} and \ref{P} hold and let $\{\vx^k\}_{0 \leq k \leqslant K}$, with $K \geqslant 2$, be generated by $\SGD$ using cosine step sizes $\alpha_k = \alpha[\frac{1+\cos(\frac{k\pi}{K})}{2}]^p$ with $p>0$. Define $\varrho_s$, $\omega_s$, $\zeta_s$, $\xi_s$ and $\bar{\alpha}_s$ as in equations \eqref{eq:notation-sgd} and \eqref{eq:stepsize-ub}. If $0<\alpha  \leqslant \bar{\alpha}_s$, then it holds that:
    \vspace{-1mm}
	\begin{equation*}\label{eq:sgd_cos1}
		y_K \leqslant  \max\cbra{\frac{2\zeta_s^{\varrho_s}}{\mu^{\omega_s}}\bra{\frac{2D}{\theta K}}^{\omega_s},4\Big(\frac{\pi}{2}\Big)^{\frac{p}{\theta}}\bra{\frac{\zeta_s^{2p\varrho_s+\varrho_s}\alpha^{\omega_s}}{\mu^{2p\omega_s}}}^{\frac{1}{2p+\varrho_s}}\bra{\frac{2D}{\theta K}}^{\frac{2p\omega_s}{2p+\varrho_s}}} + y_0 \exp\prt{-\frac{\xi_s\alpha^{\frac{1}{\varrho_s}}}{2^{\max\{1,p/\varrho_s\}}}K},
	\end{equation*}
    where $D = \max\{1,\frac{p\pi^2}{\theta}\}$ and we have
    \begin{equation}\label{eq:sgd_cos2}
    y_K = \cO({K^{-\frac{2p\omega_s}{2p+\varrho_s}}}).
    \end{equation}
    In the setting where $\theta$ is known, choosing $\alpha= [\frac{\beta\log(K)}{K}]^{\varrho_s}$ and assuming $\beta\geqslant {2^{\max\{1,p/\varrho_s\}}\omega_s}/{\xi_s}$ and $K\geqslant {\beta\log(K)}/{\bar{\alpha}_s^{{1}/{\varrho_s}}}$, it follows that
    \begin{equation}\label{eq:sgd_cos3}
    y_K = \cO\prt{\bra{\frac{(\log(K))^{\frac{\varrho_s}{2p+\varrho_s}}}{K}}^{\omega_s}}.
    \end{equation}
\end{thm}

The proof is deferred to \Cref{proof:cos}. To our knowledge, the only existing convergence analysis of $\SGD$ with cosine step sizes under a PL-type condition appears in \citep{li2021second}. For the specific case $p = 1$ and $\theta = \frac12$, \cite[Theorem 2]{li2021second} establishes a complexity bound of $\cO(K^{-\frac{2}{3}})$, which is recovered by our result in \eqref{eq:sgd_cos2}. However, fixing $\theta = \frac{1}{2}$ implies that the parameter is known. In this scenario, \Cref{thm:sgd_cos} demonstrates that a tailored selection of the step size parameter $\alpha$ yields an improved rate of $\cO(\log^{1/3}(K)/K)$; cf.\ \eqref{eq:sgd_cos3} with $\varrho_s = \omega_s = 1$. Furthermore, the results for the range $\theta\in (\frac{1}{2},1]$ appear to be novel.

\begin{thm}[Rate of $\RR$ with Cosine Step Sizes]\label{thm:rr_cos}
Assume that \ref{B1}--\ref{B3} and \ref{P} hold. Let $\{\vx^k\}_{0 \leq k \leqslant K}$, with $K \geqslant 2$, be generated by $\RR$ using the step size $\alpha_k = \alpha[\frac{1+\cos(\frac{k\pi}{K})}{2}]^p$ with $p>0$. Let $\varrho_r$, $\omega_r$, $\zetarc$, $\xirc$ and $\bar{\alpha}_r$ be given as in \eqref{eq:notation-rr}--\eqref{eq:stepsize-ub}. If $0<\alpha  \leqslant \bar{\alpha}_r$, then we have:
\begin{equation*}
    y_K \leqslant  \max\cbra{\frac{2\zetarc^{\varrho_r}}{(\mu/2)^{\omega_r}}\bra{\frac{2D_r^{(1)}}{\theta \sqrt{N} K}}^{\omega_r},D_r^{(2)}\bra{\frac{\alpha}{\sqrt{N}}}^{\frac{\omega_r}{2p+\varrho_r}}\bra{\frac{2D_r^{(1)}}{\theta\sqrt{N}K}}^{\frac{2p\omega_r}{2p+\varrho_r}}} + y_0 \exp\prt{-\frac{\xirc\alpha^{\frac{1}{\varrho_r}}N^{\frac{1}{2\theta}-1}K}{2^{\max\{1,p/\varrho_r\}}}},
\end{equation*}
where $D_r^{(1)} = \max\{1,\frac{2p\pi^2}{\theta}\}$ and $D_r^{(2)}=4(\frac{\pi^2}{4})^{\frac{1}{\theta}}[{\frac{\zetarc^{2p\varrho+\varrho}}{(\mu/2)^{2p\omega}}}]^{\frac{1}{2p+\varrho}}$. Furthermore, if 
$K\geqslant N^{1-\frac{1}{2\theta}}\cdot \frac{2^{\max\{1,p/\varrho_r\}}\omega_r \log(\sqrt{N}K)}{\xirc \alpha^{1/\varrho_r}} $, it holds that
\begin{equation}\label{eq:rr_cos1}
    y_K = \cO({\sqrt{N}^{-\frac{\omega_r}{2p+\varrho_r}}(\sqrt{N}K)^{-\frac{2p\omega_r}{2p+\varrho_r}}}).
\end{equation}
In the setting where the {\L}ojasiewicz exponent $\theta$ is known, selecting $\alpha= [\frac{\beta\log(\sqrt{N}K)}{K} N^{1-\frac{1}{2\theta}}]^{\varrho_r}$ and assuming $\beta\geqslant { 2^{\max\{1,p/\varrho\}}\omega_r}/{\xirc}$ and $K\geqslant {\beta \log(\sqrt{N}K) N^{1-\frac{1}{2\theta}}}/{\bar{\alpha}_r^{{1}/{\varrho_r}}}$, we obtain
	\begin{equation}\label{eq:rr_cos2}
    y_K = \cO\Big( \Big[\frac{(\log(NK))^\frac{\varrho_r}{2p+\varrho_r}}{\sqrt{N}K}\Big]^{\omega_r}\Big).
	\end{equation}	
\end{thm}

The proof can be found in \Cref{proof:cos}. The results in \Cref{thm:rr_cos} seem to be new. 

\subsection{Analyzing Constant and Polynomial Step Sizes}\label{sec}

For completeness, we conduct the aforementioned analysis for constant and polynomial step size rules. The findings are presented in Theorems~\ref{thm:sgd_const}, \ref{thm:rr_const}, \ref{thm:sgd_poly}, and \ref{thm:rr_poly}, accompanied by an in-depth comparison with previous research.

\subsubsection{Constant Step Size}  
We now characterize the convergence of $\SGD$ and $\RR$ in the constant step size regime.

        \begin{thm}[Rate for $\SGD$ with Constant Step Size]\label{thm:sgd_const}
	Let conditions \ref{A1}--\ref{A2} and \ref{P} hold and let the sequence $\{\vx^k\}_{0\leqslant k\leqslant K}$, $K\geqslant 2$, be generated by $\SGD$ with $\alpha_k  = \alpha$ where $\alpha>0$. Define $\varrho_s$, $\omega_s$, $\zeta_s$, $\xi_s$ and $\bar{\alpha}_s$ as in \eqref{eq:notation-sgd} and \eqref{eq:stepsize-ub}. If $\alpha \leqslant
            		\bar{\alpha}_s$, then we have: 
	\begin{equation}\label{eq:sgd_const1}
		y_{K} \leqslant 2\zeta_s \alpha^{\frac{1}{2\theta}} + y_0 \exp(-\xi_s \alpha^{1/\varrho_s}K). 
	\end{equation}
	If $\theta$ is known, selecting $\alpha = [\frac{\beta\log(K)}{K}]^{\varrho_s}$ and assuming $\beta \geqslant \frac{\omega_s}{\xi_s}$ and $K \geqslant \frac{\beta\log(K)}{\bar{\alpha}_s^{1/\varrho_s}}$, it holds that
	\begin{equation}\label{eq:sgd_const2}
		y_{K} \leqslant 2\xi_s \beta^{\omega_s} \Big[\frac{\log(K)}{K}\Big]^{\omega_s} + y_0 K^{-\beta \xi_s} = \cO\prt{\Big[\frac{\log(K)}{K}\Big]^{\omega_s}}.
	\end{equation}
\end{thm}

The proof is provided in \Cref{proof:constant}. $\SGD$ under constant step sizes has been studied extensively, but mostly limited to the case where $\theta = \frac{1}{2}$ (or under strong convexity); see, e.g.,  \citep{gower2019sgd,bottou2018optimization,needell2014stochastic,ShaZha2013}. We can recover the existing results and obtain new complexity bounds for general exponents $\theta\in (\frac{1}{2},1]$.

\begin{thm}[Rate for $\RR$ with Constant Step Size]\label{thm:rr_const}
    Assume that \ref{B1}--\ref{B3} and \ref{P} hold. Let $\{\vx^k\}_{0 \leqslant k \leqslant K}$, with $K \geqslant 2$, be generated by $\RR$ using the step size $\alpha_k  = \alpha$. Define $\varrho_r$, $\omega_r$, $\zetarc$, $\xirc$ and $\bar{\alpha}_r$ as in \eqref{eq:notation-rr} and \eqref{eq:stepsize-ub}. If $0<\alpha  \leqslant \bar{\alpha}_r$, then the following result holds:
	\begin{equation}\label{eq:rr_const1}
		y_{K} \leqslant \frac{2\zetarc \alpha^{\frac{1}{\theta}}}{N^{\frac{1}{2\theta}}}  + y_0 \exp\Big(-\frac{\xirc \alpha^{1/\varrho_r}K}{N^{1-\frac{1}{2\theta}}}\Big). 
	\end{equation}
	If $\theta$ is known, setting $\alpha= [\frac{\beta\log(\sqrt{N}K)}{K} N^{1-\frac{1}{2\theta}}]^{\varrho_r}$ and assuming $K\geqslant \frac{\beta \log(\sqrt{N}K) N^{1-\frac{1}{2\theta}}}{\bar{\alpha}_r^{{1}/{\varrho_r}}}$ and $\beta\geqslant {\omega_r}/{\xirc}$, we have
	\begin{equation}\label{eq:rr_const2}
		y_{K} \leqslant 2 \zetarc \bra{\frac{\beta \log(\sqrt{N}K)}{\sqrt{N}K}}^{\omega_r} + y_0 (\sqrt{N}K)^{-\beta\xirc} = \cO\prt{\bra{\frac{\log(\sqrt{N}K)}{\sqrt{N}K}}^{\omega_r}}.
	\end{equation}
\end{thm}

The proof is deferred to \Cref{proof:constant}. Previous work on $\RR$ with constant step sizes has mainly focused on the specific setting $\theta = \frac{1}{2}$ or on strong convexity, see \citep{haochen2019,ahn2020sgd,Mis2020,nguyen2020unified}. Notably, \cite[Theorem 1]{nguyen2020unified} provides the complexity bound $\cO(\frac{\sL^2\sigma^2}{\mu^3} \frac{\log^2(NK)}{NK^2} + y_0\frac{1}{NK^2})$ for strongly convex functions, which aligns with our result when $\theta = \frac{1}{2}$. We extend this analysis to the range $\theta\in (\frac{1}{2},1]$, a regime that currently appears to be unexplored.

\subsubsection{Results on Polynomial Step Sizes.}  \label{sec:poly-lem}
We now study the convergence of $\SGD$ and $\RR$ under the general $(\theta,\mu)$-PL condition when using polynomial step sizes. We first proceed with a technical lemma to provide the convergence analysis for sequences $\{y_k\}_k$ that follow the key algorithmic recursion \eqref{PL-recursion}.

   \begin{lem}\label{lem:poly}
	Let $\{y_k\}_{k}$ be a given sequence satisfying the following recursion: 
	\begin{equation} 
		y_{k+1} \leqslant (1 + \ell_1 \alpha_k^\tau)  y_k  - \ell_2 \alpha_k y_k^{2\theta}+ \ell_3 \alpha_k^\tau, \label{eq:am-nice-02}
	\end{equation}
	where $\ell_1,\ell_2,\ell_3,p>0$, $\tau>1$, $\theta\in [\frac{1}{2},1]$ and $\alpha_k = \frac{\alpha}{(k+\gamma)^p}$, $\alpha,\gamma> 0$. Let us set 
    \[ \varrho := \frac{2\theta}{(2\theta-1)\tau+1}, \;\; \mathcal{u}_1 = \frac{p(\tau-1)}{2\theta}, \;\; \mathcal{u}_2 = \frac{\tau-1}{(2\theta-1)\tau+1}, \;\; \mathcal{u}_3 = \frac{1-p}{2\theta-1}, \]
    \[\ukc = \max\{(2\theta - 1)\delta, ({\ell_3}/{\ell_2})^{\frac{1}{2\theta}}\}, \quad \xi := \theta \ell_2 \ukc^{2\theta-1} = \max\{\theta\mathfrak{p}(\theta) \ell_2 \delta^{2\theta-1},  \theta \ell_2^{\frac{1}{2\theta}}\ell_3^{\frac{2\theta-1}{2\theta}}\}. \] 
    where $\delta>0$. The following statements hold:
	\begin{enumerate}[label=\textup{\textrm{(\alph*)}},leftmargin=8ex,topsep=2pt,itemsep=0ex,partopsep=0ex]
			\item If $p\in (0,\varrho)$, $\frac{\alpha}{\gamma^p} \leqslant \min\{(\frac{\theta \mathfrak{p}(\theta)\ell_2\delta^{2\theta-1}}{\ell_1})^{\frac{2\theta}{\tau - 1}},\xi^{-\varrho}\}$, and $\gamma \geqslant (\frac{2\mathcal{u}_1}{\xi\alpha^{1/\varrho}})^{\frac{1}{1-p/\varrho}}$, then we have
            \[
                y_{k+1} \leqslant 4\zeta \alpha^{\frac{\tau-1}{2\theta}} (k+1+\gamma)^{-\mathcal{u}_1} + y_0 \exp\prt{-\frac{\xi\alpha^{1/\varrho}(k+1)}{(k+1+\gamma)^{p/\varrho}}}.
            \] 
			\item If $p = \varrho$, $\frac{\alpha}{\gamma^p} \leqslant \min\{(\frac{\theta \mathfrak{p}(\theta)\ell_2\delta^{2\theta-1}}{\ell_1})^{\frac{2\theta}{\tau - 1}},\xi^{-\varrho}\}$, and $\alpha \geqslant (\frac{2\mathcal{u}_2}{\xi})^\varrho$, then it holds that
            \[
            y_{k+1} \leqslant 4\zeta \alpha^{\frac{\tau-1}{2\theta}} (k+1+\gamma)^{-\mathcal{u}_2} + y_0 (\gamma^{-1}(k+1+\gamma))^{-\xi\alpha^{1/\varrho}}.
            \]
			\item If $\theta \in (\frac{1}{2},1]$, $p\in (\varrho,1)$, $\alpha \geqslant \frac{2\mathcal{u}_3}{\theta \ell_2}$, and $\gamma \geqslant \max\{(\frac{\alpha^{\tau-1} \ell_1}{\theta \ell_2})^{\frac{1}{\tau p-1}}, (\frac{\alpha^{\tau-1} \ell_3}{\ell_2})^{\frac{1}{\tau p-\mathcal{u}_3-1}},\alpha\theta \ell_2\}$, then: 
            \begin{equation*} y_{k+1} \leqslant 4 (k+1+\gamma)^{-\mathcal{u}_3} + y_0 (\gamma^{-1}(k+1+\gamma))^{-\alpha\theta \ell_2}. \end{equation*}
			\item Suppose $\theta \in (\frac{1}{2},1]$, $p=1$. Then, for every $\alpha\geqslant \frac{2}{\theta(2\theta-1)\ell_2}$, there is 
   $\gamma_0>0$, which only depends on $\alpha,\theta,\ell_2$, such that if $\gamma\geqslant \gamma_0$, we have
            \begin{equation*} y_{k+1} \leqslant 4\log(k+1+\gamma)^{-\frac{1}{2\theta-1}} + y_0 \Big(\frac{\log(k+1+\gamma)}{\log(\gamma)}\Big)^{-\alpha\theta \ell_2}. \end{equation*}
	\end{enumerate}
\end{lem}

The proof mainly follows the roadmap in \Cref{sec:recursion&roadmap}, but requires certain refinements. We provide a detailed discussion in \Cref{proof:lem:poly}.

\begin{thm}[Rate for $\SGD$ with Polynomial Step Sizes]\label{thm:sgd_poly}
	Let \ref{A1}--\ref{A2} and \ref{P} hold and let the sequence $\{\vx^k\}_{0\leqslant k\leqslant K}$, with $K\geqslant 1$, be generated by $\SGD$ with $\alpha_k=\frac{\alpha}{(k+\gamma)^{p}}$. Under certain conditions on the parameters $\alpha$ and $\gamma$, the following result holds:
	\begin{equation}\label{eq:sgd_poly1}
		y_{K} = \begin{cases}
            \cO(K^{-\min\{\frac{p}{2\theta},\frac{1-p}{2\theta-1}\}}),\quad & \text{if $\theta \in (\frac{1}{2},1]$, $p\in (0,1)$ or $\theta = \frac{1}{2}$, $p\in (0,1]$,}\\
		    \cO(\log(K)^{-\frac{1}{2\theta-1}}), & \text{if $\theta\in (\frac{1}{2},1]$, $p=1$.}
		\end{cases}
	\end{equation}
	In the setting where $\theta$ is known, choosing $p = \varrho_s$, $\alpha \geqslant (\frac{2\omega_s}{\xi_s})^{\varrho_s}$, $\frac{\alpha}{\gamma^\varrho}\leqslant \bar{\alpha}_s$, it holds that
	\begin{equation}\label{eq:sgd_poly2}
		y_{K} \leqslant 4\ukc_s \alpha^{\frac{1}{2\theta}} (K+\gamma)^{-\omega_s} + y_0 (\gamma^{-1}(K+\gamma))^{-2\omega_s} = \cO(K^{-\omega_s}).
	\end{equation}
\end{thm}

The proof is detailed in \Cref{proof:sgd_rr_poly}. Previously, \cite[Theorem 12]{fontaine2021convergence} analyzed convergence for polynomial step sizes $\frac{\alpha}{(k+\gamma)^p}$ with $p\in (0,1)$ under a generalized weakly quasi-convex condition (see \cite[Assumption F3]{fontaine2021convergence}) and bounded variance. While their Assumption F3 is more general than our $(\theta,\mu)$-PL condition, they require additional assumptions on the iterate upper bounds to derive rates outside the $(\theta,\mu)$-PL setting. Our \Cref{thm:sgd_poly} recovers their findings for $(\theta,\mu)$-PL functions. Moreover, \Cref{thm:sgd_poly} extends to the case where $\theta \neq \frac{1}{2}$ and $p = 1$. This scenario, which is not addressed in \citep{fontaine2021convergence}, yields the rate $\cO((\log(K))^{-\frac{1}{2\theta-1}})$ for $\SGD$. Additionally, \cite[Corollary 1]{Fatkhullin2022} studies $\SGD$ under the same $(\theta,\mu)$-PL condition with a more general variance assumption, but limits the analysis to optimal polynomial step sizes. \Cref{thm:sgd_poly} recovers the canonical cases in that work, specifically where $h(t) = t$ and $b_k = 1$ in \cite[Assumption 4]{Fatkhullin2022} (or equivalently, setting $\beta = 1$ and $\tau=0$ in \cite[Corollary 1]{Fatkhullin2022}).

We can also extract some other interesting non-asymptotic results for specific $\theta$ from Lemma \ref{lem:poly} and \Cref{thm:sgd_poly}. If $\theta = \frac{1}{2}$, we can obtain the rate $\cO(\frac{\sL\sigma^2}{\mu^2 K} + y_0 [\frac{\sA\sL}{\mu^2 K}]^2)$. This result matches exactly the rate presented in \cite[Theorem 4.8]{gower2021sgd}. If $\theta = 1$ and $p \in [\frac{1}{2},1]$, the convergence rates provided in Lemma \ref{lem:poly} match the results for convex functions derived in \cite[Theorem 4]{Moulines2011}. 

\begin{thm}[Rate for $\RR$ with Polynomial Step Sizes]\label{thm:rr_poly}
	Suppose \ref{B1}--\ref{B3} and \ref{P} hold. Let the sequence $\{\vx^k\}_{0\leqslant k\leqslant K}$, $K\geqslant 3$, be generated by $\RR$ with $\alpha_k=\frac{\alpha}{(k+\gamma)^{p}}$. For fixed $N$, under certain conditions on the parameters $\alpha$ and $\gamma$, the following result holds:
	\begin{equation}\label{eq:rr_poly1}
		y_{K} = \begin{cases}
            \cO(K^{-\min\{\frac{p}{\theta},\frac{1-p}{2\theta-1}\}}),\quad & \text{if $\theta \in (\frac{1}{2},1]$, $p\in (0,1)$ or $\theta = \frac{1}{2}$, $p\in (0,1]$,}\\
		    \cO(\log(K)^{-\frac{1}{2\theta-1}}), & \text{if $\theta\in (\frac{1}{2},1]$, $p=1$.}
		\end{cases}
	\end{equation}
	In the setting where $\theta$ is known, choosing $p = \varrho_r$, $\alpha =  [\beta \log(\sqrt{N}K) N^{1-\frac{1}{2\theta}}]^{\varrho_r}$, $\beta \geqslant {2\omega_r}/{\xirc}$, $\gamma\geqslant \frac{\beta \log(\sqrt{N}K) N^{1-\frac{1}{2\theta}}}{\bar{\alpha}_r^{{1}/{\varrho_r}}}$, and $K\geqslant 2\gamma$, we have
	\begin{equation}\label{eq:rr_poly2}
	    y_{K} \leqslant 4\zetarc \beta^{\omega_r} \bra{\frac{\log(\sqrt{N}K)}{\sqrt{N}(K+\gamma)}}^{\omega_r} + y_0 (\sqrt{N}K)^{-2\omega_r} = \cO\prt{\bra{\frac{\log(\sqrt{N}K)}{\sqrt{N}K}}^{\omega_r}}.
	\end{equation}
\end{thm}

The proof is deferred to \Cref{proof:sgd_rr_poly}. To our knowledge, \Cref{thm:rr_poly} provides the first non-asymptotic convergence rates for $\RR$ under the general $(\theta,\mu)$-PL condition with polynomial step sizes. \vspace{1mm}

\begin{figure}[t]
	\centering
	\begin{tikzpicture}
	\node[right] at (0,0) {\includegraphics[width=7cm]{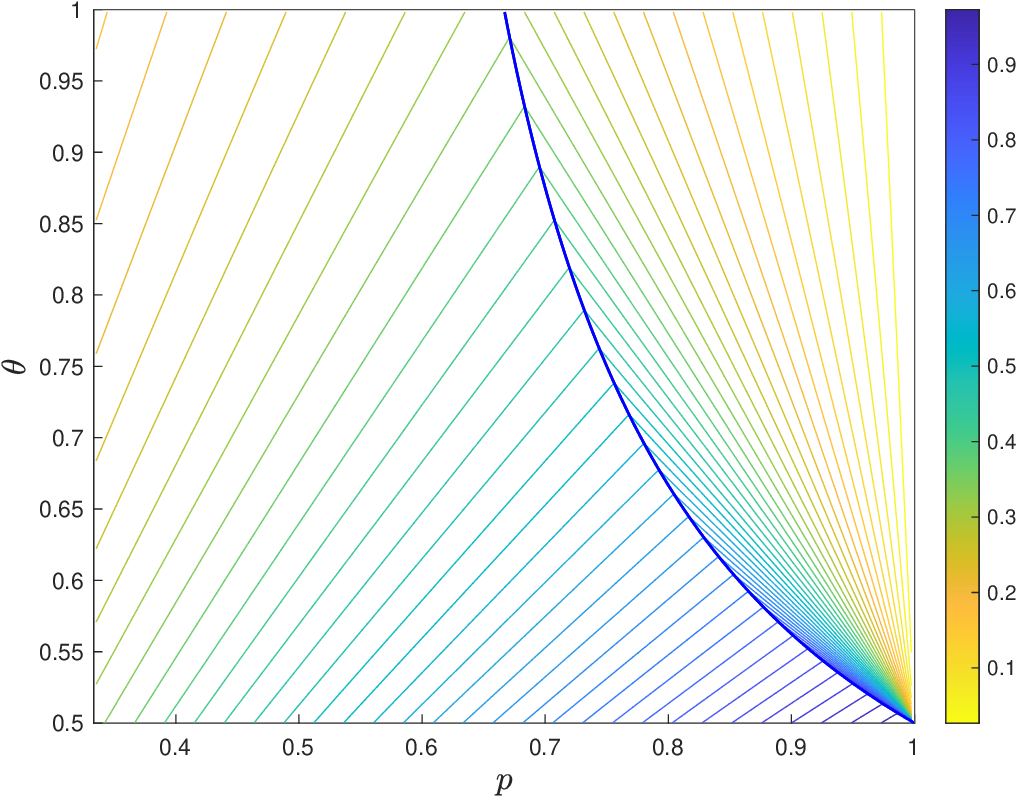}};
    \node at (3.5,3) {$\mathcal{w}_s$ ($\SGD$)};
	\node[right] at (7.5,0) {\includegraphics[width=7cm]{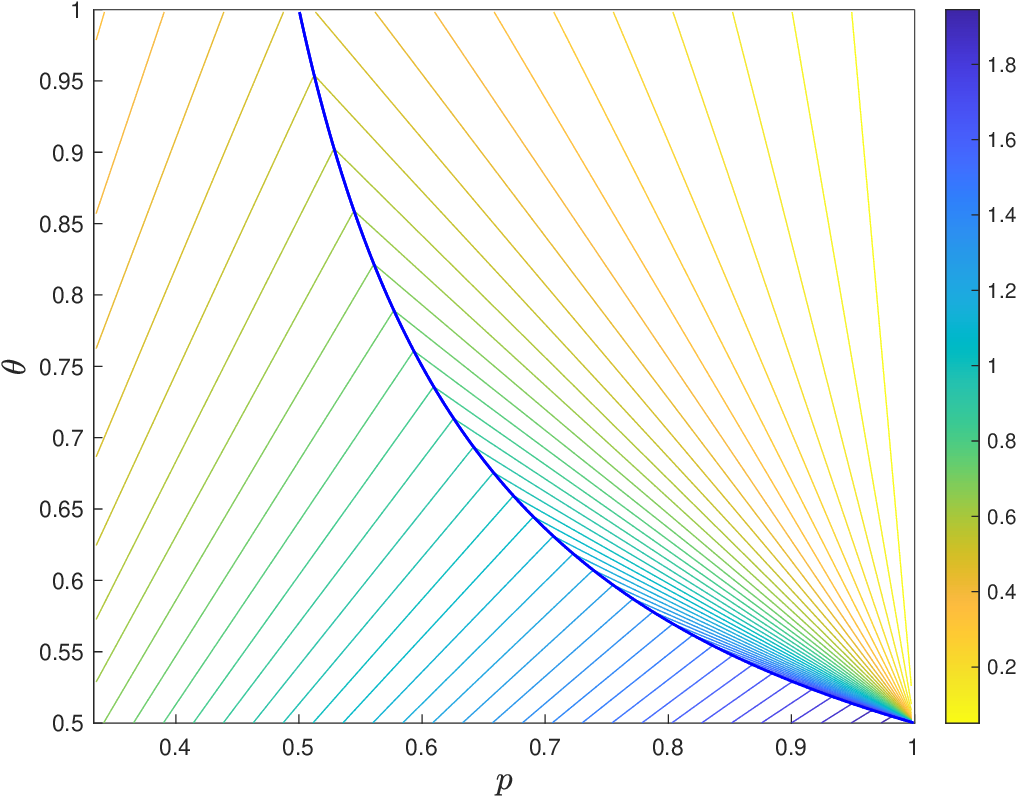}}; 
    \node at (11,3) {$\mathcal{w}_r$ ($\RR$)};
	\end{tikzpicture}
	\caption{Heatmap of the rates $\mathcal{w}_s$, $\mathcal{w}_r$ w.r.t. $p,\theta$. The blue lines depict the optimal choice of $p$ w.r.t. $\theta$; left: $p = \varrho_s = \frac{2\theta}{4\theta-1}$ ($\SGD$), right: $p = \varrho_r = \frac{\theta}{3\theta-1}$ ($\RR$). }
	\label{fig:rate-01}
\end{figure}

\noindent\textbf{Dependence of the Rate on $p$ and $\theta$.}  It is evident from Theorems~\ref{thm:sgd_poly} and \ref{thm:rr_poly} that the convergence rates depend significantly on the parameters $p$ and $\theta$. Focusing on the asymptotic dependence on $K$, in the regimes where $\theta \in (\frac{1}{2},1]$ with $p\in (0,1)$, or where $\theta = \frac{1}{2}$ with $p\in (0,1]$, $\SGD$ and $\RR$ converge with upper bounds of $\cO(K^{-\mathcal{w}_s(p,\theta)})$ and $\cO(K^{-\mathcal{w}_r(p,\theta)})$, respectively. The exponents $\mathcal{w}_s(p,\theta)$ and $\mathcal{w}_r(p,\theta)$ are defined as follows:
\[
\mathcal{w}_s(p,\theta) = \min\cbra{\frac{p}{2\theta},\frac{1-p}{2\theta-1}},
\quad
\mathcal{w}_r(p,\theta) = \min\cbra{\frac{p}{\theta},\frac{1-p}{2\theta-1}}.
\]
A visualization of the rate mappings $\mathcal{w}_s$ and $\mathcal{w}_r$ is provided in \Cref{fig:rate-01}.

\subsection{Landscape Adaptivity and Noise Adaptivity of Exponential Step Sizes}
In this section, we discuss several notable phenomena emerging from our analysis. Specifically, our results indicate that exponential step sizes demonstrate superior adaptivity to both the objective landscape and the gradient noise compared to the other three step size schedules. The following discussion focuses primarily on the $\SGD$ algorithm. 

\subsubsection{Landscape Adaptivity}\label{sec:landscape-adaptivity}
The convergence rates of $\SGD$ employing constant, polynomial, and cosine step size rules rely significantly on the choice of parameters. Achieving the optimal convergence rate necessitates a careful selection of these parameters based on $\theta$. Recalling $\varrho_s = \frac{2\theta}{4\theta-1}$, the following table summarizes how the optimal parameters for these step sizes depend on $\theta$:

\begin{table}[H] \label[table]{table:op-parameter}
\centering
\setlength{\tabcolsep}{7pt}
\NiceMatrixOptions{cell-space-limits=2pt}
\begin{NiceTabular}{|c|c|c|c|c|}%
 [ 
   code-before = 
    \rectanglecolor{lavender!40}{1-2}{2-2}
    \rectanglecolor{lavender!40}{1-4}{2-4}
 ]
\toprule
\Block{1-1}{\textbf{Step Sizes}} & \Block{1-1}{Exponential} & \Block{1-1}{Cosine} & \Block{1-1}{Constant} & \Block{1-1}{Polynomial}  \\ \Hline
\Block{1-1}{Optimal Parameters} & \Block{1-1}{--} & \Block{1-1}{$\alpha\sim [\frac{\log(K)}{K}]^{\varrho_s}$} & \Block{1-1}{$\alpha\sim [\frac{\log(K)}{K}]^{\varrho_s}$} & \Block{1-1}{$p = \varrho_s$} \\
\bottomrule
\end{NiceTabular}
\end{table}

\begin{figure}[t]
	\centering
	\begin{tikzpicture}
	\node[right] at (0,0) {\includegraphics[width=5.5cm]{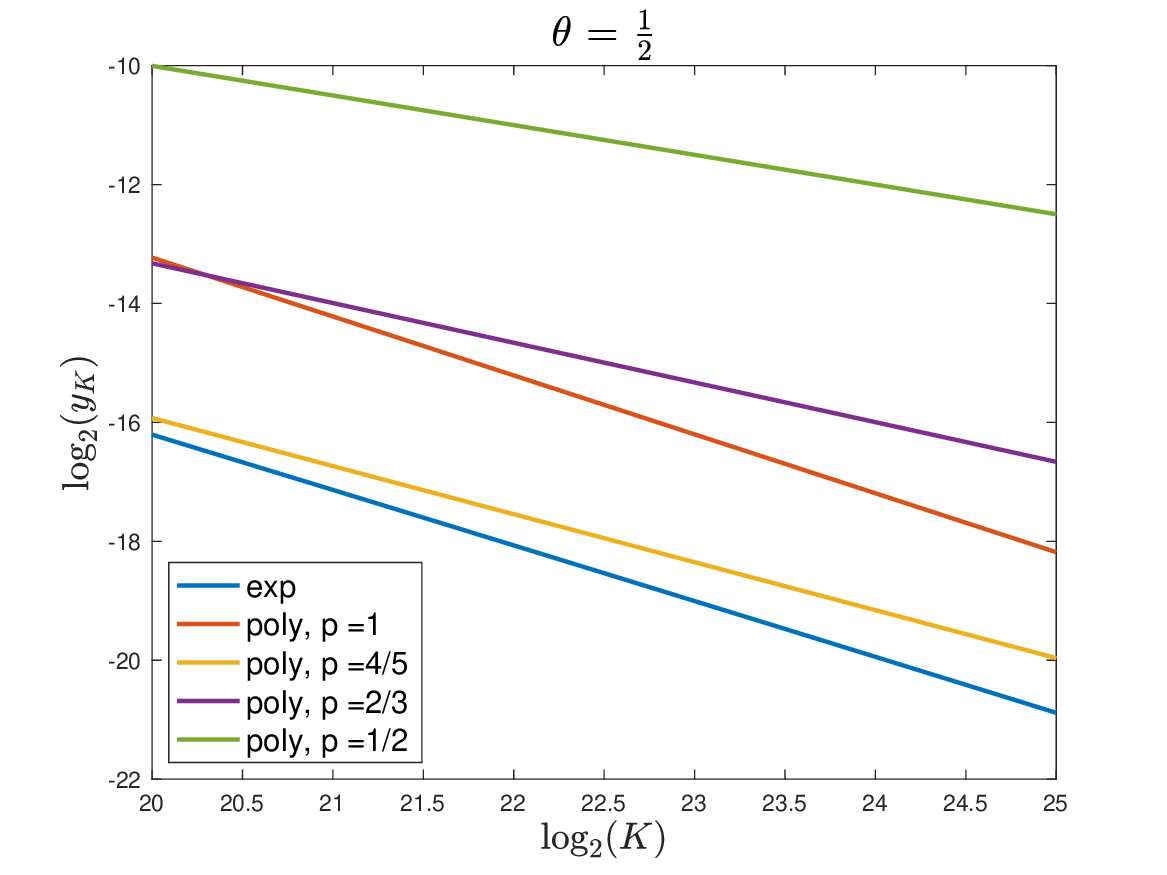}};
	\node[right] at (5,0) {\includegraphics[width=5.5cm]{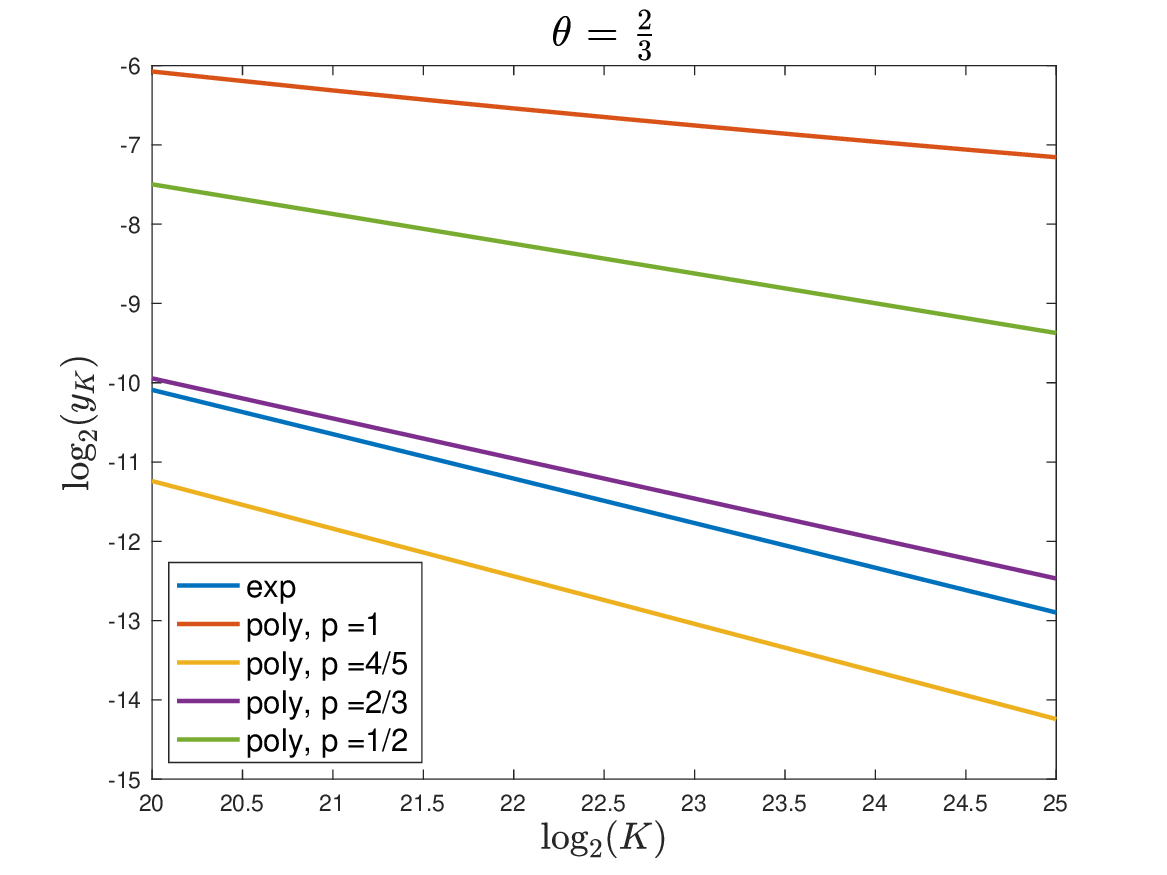}};
        \node[right] at (10,0){\includegraphics[width=5.5cm]{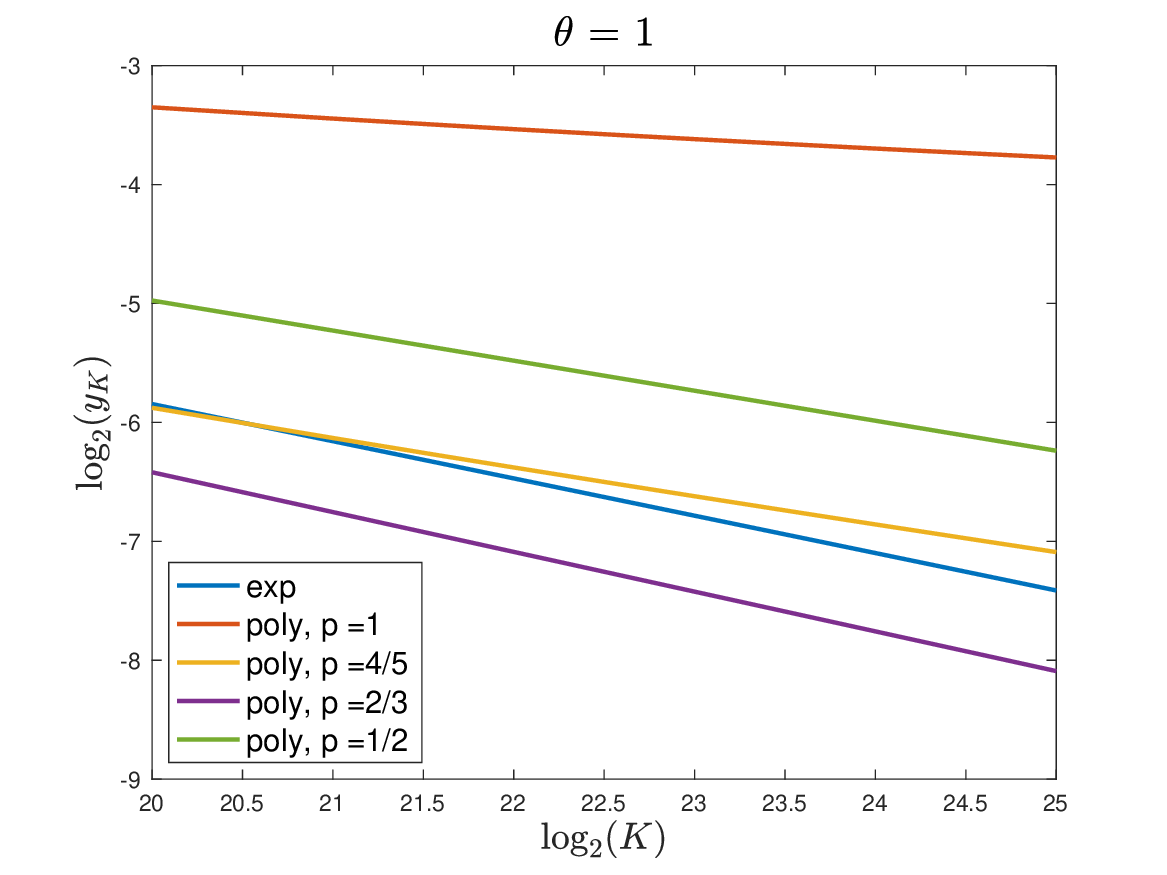}};
	\end{tikzpicture}
        \vspace{-5mm}
	\caption{Relation between $\log_2(y_K)$ and $\log_2(K)$ for different $\theta$ with polynomial and exponential step sizes. A steeper slope corresponds to a faster sublinear rate.}
	\label{fig:rate-02}
\end{figure}

\noindent(see Theorems~\ref{thm:sgd_exp}, \ref{thm:sgd_cos}, \ref{thm:sgd_const}, \ref{thm:sgd_poly}). The exponential step size rule, $\alpha_k = \alpha (\frac{\beta}{K})^{\frac{pk}{K}}$, exhibits insensitivity w.r.t. $\alpha$, $\beta$, and $p$, consistently matching the optimal convergence rate of other step size rules for any $\theta\in [\frac{1}{2},1]$ (up to a logarithmic term). We perform a simple numerical simulation to visually demonstrate this observation. The simulation involves the following recursion, representing the worst-case scenario of \eqref{PL-recursion} with $\ell_1 = \ell_2 = \ell_3 = 1$, $\tau = 2$:
\begin{equation*}
	y_{k+1} = (1+ \alpha_k^2) y_k -  \alpha_k y_k^{2\theta} + \alpha_k^2,\quad k\in [K].
\end{equation*}
 We calculate $y_K$ using polynomial step sizes $\alpha_k = \frac{\alpha}{(k+\gamma)^p}$ with $p \in \{1,\frac{4}{5},\frac{2}{3},\frac{1}{2}\}$, and an exponential step size $\alpha_k = \alpha (\frac{1}{K})^{{k}/{K}}$, for $\theta \in \{\frac{1}{2}, \frac{2}{3},1\}$, considering various total numbers of iterations $K$. In \Cref{fig:rate-02}, the relationship between $\log_2(y_K)$ and $\log_2(K)$ is depicted for large $K$. The slope of each line represents the sublinear convergence rate corresponding to each step size rule, where a steeper slope indicates a faster sublinear convergence rate. Two important facts are illustrated in this figure:
 \begin{itemize}
     \vspace{-1mm}
     \item Polynomial step sizes achieve the optimal rate with different values of $p$ for different $\theta$. The optimal choice in our simulation is $p = 1,\frac{4}{5},\frac{2}{3}$ for $\theta = \frac{1}{2},\frac{2}{3},1$ respectively. This matches the optimal choice $p = \varrho_s = \frac{2\theta}{4\theta - 1}$ in \Cref{thm:sgd_poly}.
     \vspace{-1mm}
     \item For different values of $\theta$,  exponential step sizes automatically attain the optimal rate that can be achieved by the polynomial step size with optimal $p$.
 \end{itemize}

\subsubsection{Noise Adaptivity} \label{sec:noise-adaptivity}
As discussed in Section~\ref{sec:related-work}, the noise adaptivity of exponential step sizes in $\SGD$ has been previously established for the special case $\theta = \frac{1}{2}$, \citep{li2021second,vaswani2022towards}. Our analysis demonstrates that this observation extends to the regime where $\theta \in (\frac{1}{2}, 1]$. To isolate the convergence behavior, we first present a lemma regarding the noise-free case:

\begin{lem}\label[lemma]{lem:noise-free}
    Let \ref{A1} and \ref{P} hold and let $\{\vx^k\}_{0 \leqslant k \leqslant K}$, $K\geqslant 1$, be generated by the gradient descent method with step sizes $\{\alpha_k\}_{0 \leqslant k \leqslant K-1} \subset \R_{++}$, where $\alpha_k\leqslant \frac{1}{\sL}$. Then we have
    \begin{equation}
        f(\vx^k)-f^* \leqslant \begin{cases}
        (f(\vx^0)-f^*) \exp(-\mu \sum_{k=0}^{K-1} \alpha_k),\quad &\text{if $\theta = \frac{1}{2}$,}\\
        [(f(\vx^0)-f^*)^{1-2\theta}+(2\theta-1)\mu\sum_{k=0}^{K-1}\alpha_k]^{-\frac{1}{2\theta-1}},\quad &\text{if $\theta\in (\frac{1}{2},1]$.}
    \end{cases}
    \end{equation}
\end{lem}

The proof can be found in \Cref{proof:noise-free}. \Cref{lem:noise-free} implies that, under the condition $\alpha_k \leqslant \frac{1}{\sL}$, the convergence rate is governed principally by the cumulative step sizes $\sum_{k=0}^{K-1} \alpha_k$. Based on \Cref{lem:noise-free}, the optimal achievable rates are $\cO(\exp(-\Theta(K)))$ for $\theta = \frac{1}{2}$, and $\cO(K^{-{1}/{(2\theta - 1)}})$ for $\theta \in (\frac{1}{2}, 1]$. In particular, a step size schedule achieves these optimal rates if and only if $\sum_{k=0}^{K-1} \alpha_k = \Theta(K)$. The following table summarizes estimates of $\sum_{k=0}^{K-1} \alpha_k$ for the four step size variants with parameters optimally tuned for the noisy setting:

\begin{table}[H] \label[table]{table:sum-op-step}
\centering
\setlength{\tabcolsep}{6pt}
\NiceMatrixOptions{cell-space-limits=2pt}
\begin{NiceTabular}{|c|c|c|c|c|}%
 [ 
   code-before = 
    \rectanglecolor{lavender!40}{1-2}{2-2}
    \rectanglecolor{lavender!40}{1-4}{2-4}
 ]
\toprule
\Block{1-1}{\textbf{Step Sizes}} & \Block{1-1}{Exponential} & \Block{1-1}{Cosine} & \Block{1-1}{Constant} & \Block{1-1}{Polynomial}  \\ \Hline
\Block{1-1}{$\sum_{k=0}^{K-1} \alpha_k$} & \Block{1-1}{$\geqslant \Theta(\frac{K}{\log(K/\beta)})$} & \Block{1-1}{$\Theta(K^{1-\varrho_s} [\log(K)]^{\varrho_s})$} & \Block{1-1}{$\Theta(K^{1-\varrho_s} [\log(K)]^{\varrho_s})$} & \Block{1-1}{$\Theta(K^{1-\varrho_s})$} \\
\bottomrule
\end{NiceTabular}
\end{table}

\noindent see \Cref{proof:noise-free} for additional details. We observe that $\varrho_s$ takes values in the interval $[\frac{2}{3}, 1]$ as $\theta$ varies from $\frac{1}{2}$ to $1$. According to Lemma~\ref{lem:noise-free}, parameters optimized for the noisy setting render cosine, constant, and polynomial step sizes suboptimal in the noise-free case. By contrast, exponential step sizes maintain a (near) optimal convergence rate under the same conditions. Consequently, we conclude that the exponential schedule demonstrates significantly superior noise adaptivity compared to the other three approaches.

\acks{Xiao Li was partly supported 
by the Shenzhen Science and Technology Program under grant No$.$ RCYX20221008093033010. Andre Milzarek was partly supported by the National Natural Science Foundation of China under grant No$.$ 12150410304 and W2532005, by the Shenzhen Science and Technology Program under grant No$.$ RCYX20221008093033010, 
and by the Guangdong Provincial Key Laboratory of Mathematical Foundations for Artificial Intelligence (2023B1212010001).}





\appendix

\section*{}

\section{Preparatory Tools} 

In this first section, we list well-known mathematical facts and tools that will be used repeatedly in the later analysis.

\begin{prop}[Integral Test] \label[proposition]{prop:int} Let $I = [\iota,\infty)$, $\iota \in \R$, be a given interval and let $f:I \to \R$ be continuous and non-increasing. Then, for all $a, b \in \N \cap I$ with $a \leq b$, we have
\[ \int_a^{b+1}f(x)\,\mathrm{d}x \leq \sum_{k = a}^b f(k) \leq f(a) + \int_a^bf(x)\,\mathrm{d}x. \]
\end{prop}

\begin{lem}[Technical Inequalities]\label[lemma]{lem:Tech-ineqs}
    The following inequalities hold true:
    \begin{enumerate}[label=\textup{\textrm{(\alph*)}},leftmargin=8ex,topsep=0pt,itemsep=0ex,partopsep=0ex]
        \item $\log(1+x)\leqslant x$, $\prod_{i=1}^k (1+x_i) \leqslant \exp({\sum_{i=1}^k x_i})$, for all $x,x_i\geqslant -1$.
        \item $x^r-y^r \geqslant r y^r x^{-1}(x-y)$ for any $x,y,r>0$.
        \item $2({1-\frac{k}{K}})^2 \leqslant 1+\cos(\frac{k\pi}{K}) \leqslant \frac{\pi^2}{2} ({1-\frac{k}{K}})^2$, for all $k\in [K]$, $K\geqslant 1$. 
        \item $1+\cos(\frac{(k+1)\pi}{K})\geqslant \frac{1}{2}({1-\frac{k}{K}})^2$, for all $k\in [K-2]$, $K\geqslant 2$.
        \item $\cos(\frac{(k+1)\pi}{K})-\cos(\frac{k\pi}{K}) \geqslant -\frac{\pi^2}{K}({1-\frac{k}{K}})$, for all $k\in [K-1]$, $K\geqslant 1$.
        \item $\sum_{k=0}^{K-1} [{\frac{1+\cos(\frac{k\pi}{K})}{2}}]^r \geqslant \frac{K}{2^{\max\{1,r\}}}$, for all $K \geqslant 1$, $r>0$.
    \end{enumerate}
\end{lem}
\begin{proof}
    The first inequality in (a) is well-known. The second inequality in (a) can be derived by observing that $\prod_{i=1}^k (1+x_i) = \exp({\sum_{i=1}^k \log(1+x_i)})$ and subsequently applying the first inequality.
    
Part (b) follows from the weighted AM-GM inequality, $\frac1w(w_1x_1+w_2x_2) \geq (x_1^{w_1}x_2^{w_2})^{1/w}$, for all $x_1, x_2, w_1, w_2 \geqslant 0$ and $w = w_1+w_2 > 0$. In particular, setting $x_1 = (\frac{x}{y})^r$, $x_2 = \frac{y}{x}$, $w_1 = 1$, $w_2 = r$, we have $(\frac{x}{y})^r + \frac{ry}{x} \geqslant (1 + r) ((\frac{x}{y})^r(\frac{y}{x})^r)^{1/(1+r)} = 1+r$ for any $x, y, r > 0$.
    
To prove (c), we use the double-angle formula $1-\cos(2a) = 2\sin^2(a)$. This implies $1+\cos(\frac{k\pi}{K}) = 1-\cos({\frac{(K-k)\pi}{K}}) = 2\sin^2({\frac{\pi}{2}({1-\frac{k}{K}})})$ for $k \in [K]$. Given that $x \mapsto \frac{\sin(x)}{x}$ is monotonically decreasing on $[0, \frac{\pi}{2}]$, we have $\frac{2}{\pi}x \leqslant \sin(x) \leqslant x$ for all $x \in [0, \frac{\pi}{2}]$. Substituting $x = \frac{\pi}{2}({1-\frac{k}{K}})$, we complete the proof of (c). Combining (c) with the inequality $(1-\frac{k+1}{K}) \geqslant \frac{1}{2}({1-\frac{k}{K}})$ for all $k \in [K-2]$ and $K \geqslant 2$, we establish (d). 

For part (e), we employ the sum-to-product identity $\cos(a) - \cos(b) = -2\sin({\frac{a+b}{2}})\sin({\frac{a-b}{2}})$ alongside the properties $\sin(x) = \sin(\pi-x)$ and $\sin(x) \leq x$. This yields 
\[ \cos\prt{\frac{(k+1)\pi}{K}} - \cos\prt{\frac{k\pi}{K}} = -2\sin\prt{\frac{(2k+1)\pi}{2K}}\sin\prt{\frac{\pi}{2K}} \geqslant -\frac{\pi^2}{K}\prt{1-\frac{k}{K}}. \]

To establish (f), we first prove the inequality 
\begin{equation}\label{sum-pow-ineq}
    {\sum}_{k=0}^{K-1} x_k^r \geqslant \min\cbra{K^{1-r}\prt{{\sum}_{k=0}^{K-1} x_k}^r, {\sum}_{k=0}^{K-1} x_k}
\end{equation}
    for all $0 \leqslant x_k \leqslant 1$, $K \geqslant 1$. When $r \geqslant 1$, Jensen's inequality yields 
\[ {\sum}_{k=0}^{K-1} x_k^r \geqslant K\prt{K^{-1} {\sum}_{k=0}^{K-1} x_k}^r = K^{1-r} \prt{{\sum}_{k=0}^{K-1} x_k}^r. \]
For $0 < r < 1$, the non-increasing nature of the function $r\mapsto x^r$ on $\R_{++}$, given $x\in [0,1]$, implies that $x_k^r \geqslant x_k$, leading to $\sum_{k=0}^{K-1} x_k^r \geqslant \sum_{k=0}^{K-1} x_k$. Combining the two cases, this verifies the inequality \eqref{sum-pow-ineq}. Next, setting $x_k = \frac{1+\cos(\frac{k\pi}{K})}{2}$, we apply Lagrange's trigonometric identity $\sum_{k = 0}^n \cos(ka) = \frac{1}{2} + \frac{\sin(\frac{(2n+1)a}{2})}{2\sin(\frac{a}{2})}$ to compute 
\[ {\sum}_{k=0}^{K-1} x_k = \frac{K}{2} + \frac{1}{2}{\sum}_{k=0}^{K-1} \cos\Big(\frac{k\pi}{K}\Big) = \frac{K}{2} + \frac{1}{2} \bra{\frac{1}{2} + \frac{\sin((K-\frac12)\frac{\pi}{K})}{2\sin(\frac{\pi}{2K})}} = \frac{K + 1}{2}. \]
Since $\frac{K+1}{2} \geqslant \frac{K}{2}$, substituting this into \eqref{sum-pow-ineq} concludes the proof of (f).
\end{proof}

        \section{Further Derivations for Exponential and Cosine Step Sizes} \label{proof:exp&cos-supp}

      We now study the convergence of $\SGD$ and $\RR$ under the general $(\theta,\mu)$-PL condition when using exponential and cosine step size schemes. Here, we provide detailed verifications of Theorems \ref{thm:sgd_exp} to \ref{thm:rr_cos} using Lemma \ref{lem:exp&cos}.

\subsection{\texorpdfstring{Exponential Step Sizes: Verifying Theorems~\ref{thm:sgd_exp} and \ref{thm:rr_exp}}{Exponential Step Sizes: Verifying Theorems 3.7 and 3.8}} \label{proof:exp} 

\noindent \textbf{Proof of \Cref{thm:sgd_exp}.} By assumption, we have $\alpha \leqslant \bar{\alpha}_s\leqslant \frac{1}{\sL}$, $K>\beta>0$, $\gamma<1$, and $\alpha_k = \alpha \gamma^k \leqslant \alpha \leqslant \frac{1}{\sL}$ for all $k \in [K]$. Hence, the descent property \eqref{eq:sgd-recursion} is applicable for all $k \in [K]$. The remaining conditions on $\alpha$, $\beta$, and $K$ ensure that \Cref{lem:exp&cos}~(a) is applicable. The first inequality follows directly by setting $\varrho = \varrho_s$, $\omega = \omega_s$, $\zeta = \zeta_s$, $\xi = \xi_s$, $\delta = 1$, $\ell_2 = \mu$, and $q = \frac{1}{2\theta}$ in \Cref{lem:exp&cos}~(a). If $p\geqslant \varrho_s$, \eqref{eq:sgd_exp2} can be obtained by noting that $\frac{p}{2\theta} \geqslant \frac{\varrho_s}{2\theta} = \omega_s$.
\hfill \BlackBox \\[2mm]
             
\noindent \textbf{Proof of \Cref{thm:rr_exp}.} As before, the inequality $\alpha \leqslant \bar{\alpha}_r \leqslant  \frac{1}{2\sL}$ ensures $\alpha_k \leqslant \frac{1}{2\sL}$ and thus, the descent condition \eqref{eq:rr-recursion} is applicable. The first bound in \Cref{thm:rr_exp} can be obtained by using Lemma~\ref{lem:exp&cos}~(a) with $\varrho = \varrho_r$, $\omega = \omega_r$, $\zeta = \zeta_r$, $\xi = \xi_r$, $\delta = N^{-\frac{1}{2\theta}}$, $\ell_2 = \frac{\mu}{2}$, and $q = \frac{\omega_r}{\varrho_r} = \frac{1}{\theta}$. 
If $p\geqslant \varrho_r$, then $\frac{p}{\theta}\geqslant \frac{\varrho_r}{\theta} = \omega_r$. Moreover, it holds that $\frac{1}{2\theta}\geqslant \frac{1}{2(3\theta-1)}= \frac{\omega_r}{2}$ and $N^{-\frac{1}{2\theta}}K^{-\frac{p}{\theta}} = \cO([\frac{\log({K}/{\beta})}{\sqrt{N}K}]^{\omega_r})$, which implies that the max-term is of order $\cO([\frac{\log({K}/{\beta})}{\sqrt{N}K}]^{\omega_r})$. For the second term, by $p\geqslant \varrho_r$ and $K\geqslant 2\beta$, we have $1-(\beta/K)^{p/\varrho_r} \leqslant \frac{1}{2}$ and
    \[
    \exp\Big(- \frac{\varrho_r \xirc\alpha^{\frac{1}{\varrho_r}}}{p N^{1-\frac{1}{2\theta}}}\cdot \frac{[1-(\beta/K)^{p/\varrho_r}]K}{\log({K}/{\beta})}\Big) \leqslant \exp\Big(- \frac{\varrho_r \xirc\alpha^{\frac{1}{\varrho_r}}}{p N^{1-\frac{1}{2\theta}}}\cdot \frac{K}{2\log({K}/{\beta})}\Big) \leqslant (\sqrt{N}K)^{\omega_r},
    \]
    where the last step follows from $\frac{K}{\log(K/\beta)} \geqslant \frac{2p\log(\sqrt{N}K)}{\theta\xirc \alpha^{{1}/{\varrho_r}}} N^{1-\frac{1}{2\theta}} = \frac{2p\omega_r\log(\sqrt{N}K)}{\varrho_r\xirc \alpha^{{1}/{\varrho_r}}} N^{1-\frac{1}{2\theta}}$.
    \hfill \BlackBox

\subsection{\texorpdfstring{Cosine Step Sizes: Verifying Theorems~\ref{thm:sgd_cos} and \ref{thm:rr_cos}}{Cosine Step Sizes: Verifying Theorems 3.9 and 3.10}} \label{proof:cos}
        \noindent \textbf{Proof of \Cref{thm:sgd_cos}.} By the assumption $\alpha \leqslant \bar{\alpha}_s$ and the definition of $\bar{\alpha}_s$, the step sizes satisfy $\alpha_k \leqslant \alpha \leqslant \frac{1}{\sL}$ for all $k \in [K]$. Consequently, the descent property \eqref{eq:sgd-recursion} holds for all $k \in [K]$. Applying Lemma~\ref{lem:exp&cos}~(b) to \eqref{eq:sgd-recursion} with parameters $\varrho = \varrho_s$, $\omega = \omega_s$, $\zeta = \zeta_s$, $\xi = \xi_s$, $\delta = 1$, $\ell_2 = \mu$, and $q = \frac{\omega_s}{\varrho_s} = \frac{1}{2\theta}$ yields the first bound in \Cref{thm:sgd_cos} and \eqref{eq:sgd_cos2}. 

    
    Next, consider the specific choice $\alpha= [\frac{\beta\log(K)}{K}]^{\varrho_s}$. The condition $K\geqslant \frac{\beta\log(K)}{\bar{\alpha}_s^{{1}/{\varrho_s}}}$ guarantees that $\alpha \leqslant \bar{\alpha}_s$, and thus the previously derived bound remains valid. Observe that
    \[
    \alpha^{\frac{\omega_s}{2p+\varrho_s}} K^{-{\frac{2p\omega_s}{2p+\varrho_s}}} = (\beta\log(K))^{\frac{\varrho_s\omega_s}{2p+\varrho_s}} K^{-{\frac{\varrho_s\omega_s}{2p+\varrho_s}}-{\frac{2p\omega_s}{2p+\varrho_s}}} = \bra{\frac{(\beta\log(K))^{\frac{\varrho_s}{2p+\varrho_s}}}{K}}^{\omega_s}.
    \]
    This implies that the first max-term is of order $\cO([\frac{(\log(K))^{\frac{\varrho_s}{2p+\varrho_s}}}{K}]^{\omega_s})$. Furthermore, using the definition of $\alpha$ and $\beta\geqslant {2^{\max\{1,p/\varrho_s\}}\omega_s}/{\xi_s}$, it follows
    \[
    \exp\prt{-\frac{\xi_s\alpha^{\frac{1}{\varrho_s}}}{2^{\max\{1,p/\varrho_s\}}}K} = \exp\prt{-\frac{\beta\xi_s}{2^{\max\{1,p/\varrho_s\}}}\log(K)} \leqslant K^{-\omega_s}.
    \]
    This verifies \eqref{eq:sgd_cos3}.
 \hfill \BlackBox \vspace{2mm}

\noindent \textbf{Proof of \Cref{thm:rr_cos}.} Inequality $\alpha \leqslant \bar{\alpha}_r \leqslant \frac{1}{2\sL}$ guarantees $\alpha_k \leqslant \frac{1}{2\sL}$, ensuring that the descent condition \eqref{eq:rr-recursion} holds. The first bound follows from Lemma~\ref{lem:exp&cos}~(b) with parameters $\varrho = \varrho_r$, $\omega = \omega_r$, $\zeta = \zeta_r$, $\xi = \xi_r$, $\delta = N^{-\frac{1}{2\theta}}$, $\ell_2 = \frac{\mu}{2}$, and $q = \frac{\omega_r}{\varrho_r} = \frac{1}{\theta}$. To establish \eqref{eq:rr_cos1} under 
\begin{equation}\label{eq:rr_cos_cond1}
    K\geqslant \frac{2^{\max\{1,p/\varrho_r\}}\omega_r \log(\sqrt{N}K)}{\xirc \alpha^{1/\varrho_r}} N^{1-\frac{1}{2\theta}},
\end{equation}
it suffices to estimate the order of the two terms on the right-hand side of the first bound. The first term is of order $\cO(\max\{(\sqrt{N} K)^{-\omega_r},\sqrt{N}^{-\frac{\omega_r}{2p+\varrho_r}}(\sqrt{N}K)^{-\frac{2p\omega_r}{2p+\varrho_r}}\})$. Note that \eqref{eq:rr_cos_cond1} implies that $K^{-1} \leqslant \Theta([\log(\sqrt{N}K)N^{1-\frac{1}{2\theta}}]^{-1}) = \cO(N^{-(1-\frac{1}{2\theta})})$. Then we have
    \[
    \begin{aligned}
    (\sqrt{N}K)^{-\frac{\varrho_r\omega_r}{2p+\varrho_r}} &= \sqrt{N}^{-\frac{\varrho_r\omega_r}{2p+\varrho_r}} (K^{-1})^{\frac{\varrho_r\omega_r}{2p+\varrho_r}}\\
    &= \sqrt{N}^{-\frac{\varrho_r\omega_r}{2p+\varrho_r}} \cdot \cO(N^{-(1-\frac{1}{2\theta})\frac{\varrho_r\omega_r}{2p+\varrho_r}}) = \cO(\sqrt{N}^{-(3-\frac{1}{\theta})\frac{\varrho_r\omega_r}{2p+\varrho_r}}).
    \end{aligned}
    \]
Recalling $\varrho_r = \frac{\theta}{3\theta-1}$, it follows that $(3-\frac{1}{\theta})\varrho_r = 1$ and $(\sqrt{N}K)^{-\frac{\varrho_r\omega_r}{2p+\varrho_r}} = \cO(\sqrt{N}^{-\frac{\omega_r}{2p+\varrho_r}})$. This yields
\[
(\sqrt{N} K)^{-\omega_r} = (\sqrt{N} K)^{-\frac{\varrho_r\omega_r}{2p+\varrho_r}} (\sqrt{N} K)^{-\frac{2p\omega_r}{2p+\varrho_r}} = \cO(\sqrt{N}^{-\frac{\omega_r}{2p+\varrho_r}} (\sqrt{N} K)^{-\frac{2p\omega_r}{2p+\varrho_r}}).
\]
Consequently, the first term simplifies to $\cO(\sqrt{N}^{-\frac{\omega_r}{2p+\varrho_r}}(\sqrt{N}K)^{-\frac{2p\omega_r}{2p+\varrho_r}})$. Moreover, invoking \eqref{eq:rr_cos_cond1} again, we have
    \begin{equation}\label{eq:rr-cos-pf1}
    \exp\prt{-\frac{\xirc\alpha^{\frac{1}{\varrho_r}}N^{\frac{1}{2\theta}-1}}{2^{\max\{1,p/\varrho_r\}}}K} \leqslant (\sqrt{N} K)^{-\omega_r},
    \end{equation}
    which indicates that the second term is of order $\cO((\sqrt{N} K)^{-\omega_r})$. Combining these estimates establishes \eqref{eq:rr_cos1}. 
    
    Finally, consider $\alpha= [\frac{\beta\log(\sqrt{N}K)}{K} N^{1-\frac{1}{2\theta}}]^{\varrho_r}$. If $\beta\geqslant \frac{ 2^{\max\{1,p/\varrho\}}\omega_r}{\xirc}$ and $K\geqslant \frac{\beta \log(\sqrt{N}K) N^{1-\frac{1}{2\theta}}}{\bar{\alpha}_r^{{1}/{\varrho_r}}}$, it holds that $\alpha\leqslant \bar{\alpha}_r$ and \eqref{eq:rr_cos_cond1} is satisfied; thus, the bound \eqref{eq:rr_cos1} remains valid. The final result \eqref{eq:rr_cos2} follows by substituting $\alpha$ into the term $({\alpha}/{\sqrt{N}})^{\frac{\omega_r}{2p+\varrho_r}}(\sqrt{N}K)^{-\frac{2p\omega_r}{2p+\varrho_r}}$ and using the identity $2(1-\frac{1}{2\theta})\varrho_r-1 = -\varrho_r$.  \hfill \BlackBox

    \section{Derivations for Constant Step Sizes}\label{app:constant}
    We now study convergence of $\SGD$ and $\RR$ under the general $(\theta,\mu)$-PL condition when utilizing constant step size schemes. We first proceed with a technical lemma and apply Chung's Lemma to sequences $\{y_k\}_k$ that follow the key algorithmic recursion 
    %
    %
    \eqref{PL-recursion}. In subsection~\ref{proof:constant}, we use the obtained result to provide proofs for Theorems~\ref{thm:sgd_const} and \ref{thm:rr_const}.

    \subsection{A Technical Lemma for Constant Step Size}

    \begin{lem} \label{lem:constant}
        Consider a non-negative sequence $\{y_k\}_{0\leqslant k\leqslant K}$ satisfying the recursion 
            	\begin{equation*} 
            		y_{k+1} \leqslant (1 + \ell_1 \alpha_k^\tau)  y_k  - \ell_2 \alpha_k y_k^{2\theta}+ \ell_3 \alpha_k^\tau
            	\end{equation*}
            	for all $k\in [K-1]$, where $\ell_1,\ell_2,\ell_3>0$, $K\in \mathbb{N}$ and $K\geqslant 1$. Let us define 
                \[ \varrho := \frac{2\theta}{(2\theta-1)\tau+1}, \quad \omega := \frac{\tau-1}{(2\theta-1)\tau+1}, \]
                $\ukc := \max\{(2\theta - 1)\delta, ({\ell_3}/{\ell_2})^{\frac{1}{2\theta}}\}$, and $\xi := \theta \ell_2 \ukc^{2\theta-1} = \max\{\theta\mathfrak{p}(\theta) \ell_2 \delta^{2\theta-1},  \theta \ell_2^{\frac{1}{2\theta}}\ell_3^{\frac{2\theta-1}{2\theta}}\}$, where $\delta>0$. Let us set $\alpha_k  \equiv \alpha > 0$, and let $\alpha$ satisfy $\alpha \leqslant
            		\min\{(\frac{\theta \mathfrak{p}(\theta)\ell_2\delta^{2\theta-1}}{\ell_1})^{\frac{2\theta}{\tau - 1}},\xi^{-\varrho}\}$; then: 
                    %
                    \begin{equation*} y_{K} \leqslant 2\zeta \alpha^{\frac{\omega}{\varrho}} + y_0 \exp({-\xi \alpha^{1/\varrho}K}). \end{equation*}
        \end{lem}
        \begin{proof}
            This follows from the roadmap presented in Section~\ref{sec:recursion&roadmap} (cf.\ \eqref{eq:relaxed-PL-final}), Example~\ref{ex:chung1}, $\tau-\frac{1}{\varrho} = \frac{\omega}{\varrho}$, and $(1-\xi\alpha^{1/\varrho})^K \leq \exp(-\xi\alpha^{1/\varrho}K)$.
        \end{proof}

\subsection{\texorpdfstring{Verifying Theorems~\ref{thm:sgd_const} and \ref{thm:rr_const}}{Verifying Theorems 3.11 and 3.12}} \label{proof:constant}

\noindent \textbf{Proof of \Cref{thm:sgd_const}\ }
Since $\alpha \leqslant \bar{\alpha}_s \leqslant\frac{1}{\sL}$, we have $\alpha_k \leqslant \frac{1}{\sL}$, and thus the recursion \eqref{eq:sgd-recursion} holds. The first inequality, \eqref{eq:sgd_const1}, follows directly by applying Lemma~\ref{lem:constant} with parameters $\varrho = \varrho_s$, $\omega = \omega_s$, $\zeta = \zeta_s$, $\xi = \xi_s$, $\delta = 1$, and $\frac{\omega_s}{\varrho_s} = \frac{1}{2\theta}$. Now consider the choice $\alpha = [\frac{\beta\log(K)}{K}]^{\varrho_s}$. The condition $K \geqslant {\beta\log(K)}/{\bar{\alpha}_s^{1/\varrho_s}}$ yields $\alpha\leqslant \bar{\alpha}_s$, ensuring that the bound \eqref{eq:sgd_const1} remains valid. Furthermore, we have $\exp(-\xi_s \alpha^{1/\varrho_s}K) = \exp(-\beta\xi_s \log(K))\leqslant K^{-\omega_s}$, where the last step follows from the assumption $\beta \geqslant \frac{\omega_s}{\xi_s}$. \hfill \BlackBox\\[2mm]

\noindent\textbf{Proof of \Cref{thm:rr_const}\ }
Since the condition $\alpha_k = \alpha \leqslant \frac{1}{2\sL}$ is satisfied, the recursion \eqref{eq:rr-recursion} is applicable. Applying Lemma~\ref{lem:constant} with the parameters $\varrho = \varrho_r$, $\omega = \omega_r$, $\zeta = \zeta_r$, $\xi = \xi_r$, $\delta = N^{-\frac{1}{2\theta}}$, and $\frac{\omega_r}{\varrho_r} = \frac{1}{\theta}$ yields the first bound \eqref{eq:rr_const1}. Next, we consider the special choice $\alpha= [\frac{\beta\log(\sqrt{N}K)}{K} N^{1-\frac{1}{2\theta}}]^{\varrho_r}$. Due to $K\geqslant {\beta \log(\sqrt{N}K) N^{1-\frac{1}{2\theta}}}/{\bar{\alpha}_r^{{1}/{\varrho_r}}}$, we have $\alpha\leqslant \bar{\alpha}_r$, and thus \eqref{eq:rr_const1} remains valid. Substituting this choice of $\alpha$ into \eqref{eq:rr_const1} and utilizing the relations $\frac{\varrho_r}{\theta} = \omega_r$, $(1-\frac{1}{2\theta})\omega_r-\frac{1}{2\theta} = -\frac{\omega_r}{2}$, and $\beta\xirc\geqslant \omega_r$, we obtain the final bound \eqref{eq:rr_const2}.
                \hfill \BlackBox
            
	\section{Derivations for Polynomial Step Sizes} \label{proof:poly}

    Finally, we study the convergence of $\SGD$ and $\RR$ when using polynomial step size schemes. We first provide a proof for the key technical result, Lemma \ref{lem:poly}. In \Cref{proof:sgd_rr_poly}, we then give proofs for Theorems~\ref{thm:sgd_poly} and \ref{thm:rr_poly}. Before proceeding, we first provide a useful variant of Lemma \ref{lem:non-asymptotic-Chung-main} (a).

    \begin{lem}\label{lem:chung-complement}
       Let $\{a_k\}_k \subseteq \R_{++}$ be a given sequence satisfying the following recursion:
        \begin{equation*}
        	a_{k+1} \leqslant \Big(1-\frac{c}{(k+\gamma)^\nu}\Big) a_k + \frac{d}{(k+\gamma)^{\nu+q}},\quad \nu\in(0,1],~ q>0, ~c,d > 0,~\gamma \geqslant c^{\frac{1}{\nu}}.
        \end{equation*}
        If $\nu\in (0,1)$, $\gamma \geqslant c^{\frac{1}{\nu}}$ and $\gamma \geqslant (\frac{(1+\varsigma)q}{c})^{\frac{1}{1-\nu}}$ for some $\varsigma>0$, then for all $k \in [K]$, we have
		\begin{equation} \label{eq:class-recursion} a_{k+1} \leqslant \frac{(1+\varsigma) d}{\varsigma c} (k+1+\gamma)^{-q} +  \Big(a_0 - \frac{\lambda d}{c\gamma^{q}} \Big)^+ \exp\Big(-\frac{c(k+1)}{(k+1+\gamma)^{\nu}}\Big).
		\end{equation}
    \end{lem}
    
    \begin{proof}
        Based on the result in Lemma~\ref{lem:non-asymptotic-Chung-main}~(a), if we additionally assume that $\gamma \geqslant (\frac{(1+\varsigma)q}{c})^{\frac{1}{1-\nu}}$ for some $\varsigma>0$, then it follows $\gamma^{1-\nu} \geqslant \frac{(1+\varsigma)q}{c}$ and $\lambda = \frac{c\gamma^{1-\nu}}{c\gamma^{1-\nu}-q}\leqslant \frac{1+\varsigma}{\varsigma}$. Using the concavity of the mapping $x\mapsto x^{1-\nu}$, it holds that
		\[
		(k+1+\gamma)^{1-\nu} - \gamma^{1-\nu} \geqslant \frac{1-\nu}{(k+1+\gamma)^\nu} [(k+1+\gamma)-\gamma] = (1-\nu) \frac{k+1}{(k+1+\gamma)^\nu}.
		\]
		Based on this estimate and following the derivations in \textbf{Case I} of the proof of \ref{lem:non-asymptotic-Chung-main}~(a), we readily obtain \eqref{eq:class-recursion}.
    \end{proof}
    \vspace{-8mm}

    \subsection{\texorpdfstring{Proof of Lemma \ref{lem:poly}}{Proof of Lemma 3.13}}\label{proof:lem:poly}
        \noindent For polynomial step sizes, the roadmap shown in \Cref{sec:recursion&roadmap} remains valid when $ p \leqslant \varrho $, but may no longer apply when $ p > \varrho $. This is because, when $ p > \varrho $, we have $ \alpha_k u_k^{2\theta-1} \sim k^{-p/\varrho} $ with $ p/\varrho > 1 $. In this scenario, Lemma~\ref{lem:non-asymptotic-Chung-main} is no longer applicable, and a new auxiliary sequence $ \{u_k\}_k $ must be constructed. As a result, we divide the proof into two main cases: $ p \leq \varrho $ and $ p > \varrho $. We note that $ p > \varrho $ can occur only when $ \theta \in \left(\frac{1}{2}, 1\right] $. \\[1mm]
    \noindent\textbf{Case I:} $p\leqslant \varrho$. Following the roadmap in Section~\ref{sec:recursion&roadmap}, we only need to analyze
    \[
    y_{k+1} \leqslant \Big(1-\frac{\xi\alpha^{1/\varrho}}{(k+\gamma)^{p/\varrho}}\Big) y_k + \frac{2\zeta\xi\alpha^{\tau}}{(k+\gamma)^{p\tau}}
    \]
    under the condition
    \begin{equation}\label{eq:poly-para-cond1}
    {\alpha}/{\gamma^p} \leqslant \min\{({\theta \mathfrak{p}(\theta)\ell_2\delta^{2\theta-1}}/{\ell_1})^{\frac{2\theta}{\tau - 1}},\xi^{-\varrho}\}.
    \end{equation}
    \noindent\textbf{Sub-Case I-1:} $p\in (0,\varrho)$. We simply substitute the following into Lemma~\ref{lem:chung-complement}: 
    $$
    c =\xi\alpha^{1/\varrho},\quad d = 2\zeta\xi\alpha^{\tau}, \quad \nu = p/\varrho, \quad q = p \cdot ({\tau-1})/{2\theta} = \mathcal{u}_1, \quad \varsigma = 1.
    $$
    Note that the condition $ \gamma \geqslant c^{1/\nu} $ is already implied by \eqref{eq:poly-para-cond1}. Thus, we only need to assume $\gamma \geqslant ({2\mathcal{u}_1}/({\xi\alpha^{1/\varrho}}))^{\frac{1}{1-p/\varrho}}$.  By Lemma~\ref{lem:chung-complement} and using $\tau-\frac{1}{\varrho} = \frac{\tau-1}{2\theta}$, we obtain
    \[
        y_{k+1} \leqslant 4\zeta \alpha^{\frac{\tau-1}{2\theta}} (k+1+\gamma)^{-\mathcal{u}_1} + y_0 \exp\prt{-\frac{\xi\alpha^{1/\varrho}(k+1)}{(k+1+\gamma)^{p/\varrho}}}.
    \]
    \noindent\textbf{Sub-Case I-2:} $p=\varrho$. We simply substitute the following terms into Lemma~\ref{lem:non-asymptotic-Chung-main} (b): 
    \[ c =\xi\alpha^{1/\varrho},\quad d = 2\zeta\xi\alpha^{\tau}, \quad \nu = 1, \quad q = p\tau-1 = \varrho\tau-1 = \mathcal{u}_2. \]
    Note that $ \gamma \geqslant c $ is ensured by \eqref{eq:poly-para-cond1}; we only need to additionally satisfy $ c > q $. For simplicity, we assume $ c \geqslant 2q $, which is equivalent to $\alpha \geqslant ({2\mathcal{u}_2}/{\xi})^\varrho$. Then, by Lemma~\ref{lem:non-asymptotic-Chung-main} (b):
    \[
        y_{k+1} \leqslant 4\zeta \alpha^{\frac{\tau-1}{2\theta}} (k+1+\gamma)^{-\mathcal{u}_2} + y_0 (\gamma^{-1}(k+1+\gamma))^{-\xi\alpha^{1/\varrho}}.
    \]
    \noindent\textbf{Case II:} $p\in (\varrho,1]$. Note that this can occur only when $\theta\in (\frac{1}{2},1]$. We again consider the sub-cases $p\in (\varrho,1)$ and $p = 1$ separately.\\[1mm]
    \noindent \textbf{Sub-Case II-1:} $\theta\in (\frac{1}{2},1]$ and $p\in (\varrho,1)$. To ensure applicability of Lemma~\ref{lem:non-asymptotic-Chung-main}, we select the auxiliary sequence $\{u_k\}_k$ such that $ \alpha_k u_k^{2\theta - 1} \sim \frac{1}{k} $. This requirement leads to the choice $ u_k = (k + \gamma)^{-\mathcal{u}_3} $, where $\mathcal{u}_3 = (1-p)/(2\theta-1)$. Then, \eqref{eq:roadmap-step2-1} can be written as
				\[
				y_{k+1} \leqslant \Big[1 + \frac{\ell_1 \alpha^\tau}{(k+\gamma)^{\tau p}} - \frac{2\theta \ell_2 \alpha}{k+\gamma}\Big] y_k + \Big[\frac{(2\theta-1)\ell_2\alpha}{(k+\gamma)^{1+\mathcal{u}_3}} +\frac{\ell_3 \alpha^\tau}{(k+\gamma)^{\tau p}}\Big].
				\]
				If $\alpha \geqslant \frac{2\mathcal{u}_3}{\theta \ell_2}$, $\gamma \geqslant \max\{(\frac{\alpha^{\tau-1} \ell_1}{\theta \ell_2})^{\frac{1}{\tau p-1}}, (\frac{\alpha^{\tau-1} \ell_3}{\ell_2})^{\frac{1}{\tau p-\mathcal{u}_3-1}},\alpha\theta \ell_2\}$ and noting $\tau p > 1+\mathcal{u}_3$, it follows
				\begin{equation*}
						\frac{\ell_1 \alpha^\tau}{(k+\gamma)^{\tau p}} - \frac{2\theta \ell_2 \alpha}{k+\gamma} \leqslant -\frac{\theta \ell_2 \alpha}{k+\gamma}, \quad \frac{\ell_3 \alpha^\tau}{(k+\gamma)^{\tau p}} \leqslant \frac{\ell_2\alpha}{(k+\gamma)^{1+\mathcal{u}_3}},\quad \frac{\theta \ell_2 \alpha}{k+\gamma} \leqslant 1, 
				\end{equation*}
				and $\alpha \theta \ell_2 \geqslant 2\mathcal{u}_3$ for all $k\geqslant 0$. Consequently, we can infer
				\[
				y_{k+1} \leqslant \Big[1- \frac{\theta \ell_2\alpha}{k+\gamma}\Big] y_k + \frac{2\theta \ell_2 \alpha}{(k+\gamma)^{1+\mathcal{u}_3}}.
				\]
				Thus, applying Lemma~\ref{lem:non-asymptotic-Chung-main} (b) with $c = \theta \ell_2\alpha$, $d = 2\theta \ell_2\alpha$, $q = \mathcal{u}_3$ and using $c-q\geqslant \frac{c}{2}$, we obtain $y_{k+1} \leqslant 4 (k+1+\gamma)^{-\mathcal{u}_3} + y_0 [\gamma^{-1}(k+1+\gamma)]^{-\theta \ell_2\alpha}$. \\[1mm]
    \noindent \textbf{Sub-Case II-2:} $\theta\in (\frac{1}{2},1]$ and $p=1$. Here, we introduce the auxiliary sequence $u_k:=\log(k+\gamma)^{-\frac{1}{2\theta-1}}$. Mimicking the earlier cases, 
there exists $\gamma_0 > 0$, which solely depends on $\alpha,\theta,\ell_2$, such that for any $\gamma\geqslant \gamma_0$, we have $\gamma \log(\gamma) \geqslant \alpha\theta \ell_2$ and \eqref{eq:roadmap-step2-1} can be simplified to 
			\[
			y_{k+1} \leqslant \Big[1 - \frac{\alpha\theta \ell_2}{(k+\gamma) \log(k+\gamma)}\Big] y_k + \frac{2\alpha\theta \ell_2 }{(k+\gamma)\log(k+\gamma)^{{2\theta}/({2\theta-1})}}.
			\]
            (We will omit an explicit construction of such $\gamma_0$ here). Let us now define $s(x) := \frac{1}{\alpha\theta \ell_2} x \log(x)$, $t(x) := \frac{1}{2\alpha\theta \ell_2}x \log(x)^{\frac{2\theta}{2\theta-1}}$, $r(x) := 2 \log(x)^{-\frac{1}{2\theta-1}}$, and $b_k = k+\gamma$. Then, it holds that
			\begin{align*}
			u(x) = (\log(r))'(x) s(x) = 
   -({\alpha\theta(2\theta-1)\ell_2})^{-1}.
			\end{align*}
			Thus, if $\alpha\geqslant \frac{2}{\theta(2\theta-1)\ell_2}$, it follows $(b_{k+1}-b_k) u(b_k) \geqslant -\frac{1}{2}$. Since the rate function $r$ is clearly convex on $I = [\gamma,\infty)$, all assumptions in \Cref{main-thm} (b) are satisfied with $\lambda_k\equiv 2$. Hence, applying \Cref{main-thm} (b), we can infer
			\[
			\begin{aligned}
				y_{k+1} \leqslant ~&4\log(k+1+\gamma)^{-\frac{1}{2\theta-1}} + y_0 \exp\Big(-{\sum}_{i=0}^{k} \frac{\alpha\theta \ell_2}{(i+\gamma) \log(i+\gamma)}\Big)\\
				\leqslant ~&4\log(k+1+\gamma)^{-\frac{1}{2\theta-1}} + y_0 ({\log(k+1+\gamma)}/{\log(\gamma)})^{-\alpha\theta \ell_2},
			\end{aligned}
			\]
            where we utilize the integral test $\sum_{i=0}^k \frac{1}{(k+\gamma)\log(k+\gamma)} \geqslant \int_0^{k+1} \frac{1}{(x+\gamma)\log(x+\gamma)}\,\mathrm{d}x = \log(\log(k+1+\gamma))-\log(\log(\gamma))$ in the last line. 
            \hfill\BlackBox

        \subsection{Results for Polynomial Step Sizes}\label{proof:sgd_rr_poly}

\noindent \textbf{Proof of \Cref{thm:sgd_poly}\ }
                The result \eqref{eq:sgd_poly1} is obtained directly from Lemma~\ref{lem:poly} by setting $\tau = 2$, $\ell_1 = \frac{\sA \sL}{2}$, $\ell_2 = \mu$, and $\ell_3 = \frac{\sL \sigma^2}{2}$. In this case, the parameters in Lemma~\ref{lem:poly} correspond to $\varrho = \varrho_s$, $\omega = \omega_s$, $\zeta = \zeta_s$, $\xi = \xi_s$, and $\delta = 1$. Next, by specifying $p = \varrho_s$, the condition $\frac{\alpha}{\gamma^\varrho}\leqslant \bar{\alpha}_s$ ensures that $\alpha_k\leqslant \frac{1}{\sL}$, thereby implicitly validating the recursion estimate \eqref{eq:sgd-recursion}. Substituting $p = \varrho_s$, ${\mathcal u}_2 = \omega_s$, and $\alpha \geqslant (\frac{2\omega_s}{\xi_s})^{\varrho_s}$ into Lemma~\ref{lem:poly}~(b) yields \eqref{eq:sgd_poly2}. 
                \hfill \BlackBox \\[2mm]
               
\noindent \textbf{Proof of \Cref{thm:rr_poly}\ }
    The bound \eqref{eq:rr_poly1} follows from Lemma~\ref{lem:poly} by setting $\tau = 3$, $\ell_1 = \frac{\sA\sL^2}{2N}$, $\ell_2 = \frac{\mu}{2}$, $\ell_3 = \frac{\sL^2 \sigma^2}{2N}$, and $\delta = N^{-\frac{1}{2\theta}}$. The parameters in Lemma~\ref{lem:poly} then reduce to $\varrho = \varrho_r$, $\omega = \omega_r$, $\zeta = \zeta_r$, $\xi = \xi_r$. With the specific choices $p = \varrho_r$ and $\alpha = [\beta \log(\sqrt{N}K) N^{1-\frac{1}{2\theta}}]^{\varrho_r}$, the condition $\gamma\geqslant {\beta \log(\sqrt{N}K) N^{1-\frac{1}{2\theta}}}/{\bar{\alpha}_r^{{1}/{\varrho_r}}}$ guarantees $\frac{\alpha}{\gamma^p}\leqslant \bar{\alpha}_r$. Consequently, we have $\alpha_k\leqslant \frac{\alpha}{\gamma^p}\leqslant\bar{\alpha}_r\leqslant \frac{1}{2\sL}$, which validates the recursion \eqref{eq:rr-recursion}. Moreover, the bound $\frac{\alpha}{\gamma^p}\leqslant \bar{\alpha}_r$ directly implies the first requirement of Lemma~\ref{lem:poly}~(b), namely $\frac{\alpha}{\gamma^p} \leqslant \min\{(\frac{\theta \mathfrak{p}(\theta)\ell_2\delta^{2\theta-1}}{\ell_1})^{\frac{2\theta}{\tau - 1}},\xi^{-\varrho}\}$. Furthermore, provided that $\beta \geqslant {2\omega_r}/{\xirc}$ and $K\geqslant 3$, the second requirement, $\alpha \geqslant (\frac{2\mathcal{u}_2}{\xi})^\varrho$, is also satisfied (thanks to $\mathcal{u}_2 = \omega_r$). Substituting these parameters into Lemma~\ref{lem:poly}~(b) and using $\frac{\varrho_r}{\theta} = \omega_r$, we have:
    \[
    y_{K} \leqslant 4\zetarc \beta^{\omega_r} \bra{\frac{\log(\sqrt{N}K)}{\sqrt{N}(K+\gamma)}}^{\omega_r} + y_0 (\gamma^{-1}(K+\gamma))^{-\xirc\beta \log(\sqrt{N}K)}.
\]
Finally, we observe that $(\gamma^{-1}(K+\gamma))^{-\xirc\beta \log(\sqrt{N}K)} = (\sqrt{N}K)^{-\frac{\log(\gamma^{-1}(K+\gamma))}{\log(\sqrt{N}K)}\cdot \xirc\beta \log(\sqrt{N}K)} \leqslant (\sqrt{N}K)^{-2\omega_r}$,
where the last inequality holds because $K+\gamma\geqslant 3\gamma$ implies $\log(\gamma^{-1}(K+\gamma))\geqslant \log(3)>1$ and it holds that $\beta \geqslant {2\omega_r}/{\xirc}$.
\hfill\BlackBox

\section{Derivations for the Noise-free Case}\label[appendix]{proof:noise-free}
We first present a self-contained proof of \Cref{lem:noise-free}, which establishes the convergence rate for the deterministic gradient descent method under the $(\theta,\mu)$-PL condition. Here, let us also refer to Theorem 4 in \cite{xiao2020policy} for the special case $\alpha_k \equiv \alpha$.

\noindent \textbf{Proof of \Cref{lem:noise-free}\ }
    By the $\sL$-smoothness of $f$ and the descent lemma, we obtain:
    \[
    f(\vx^{k+1})-f^* \leqslant f(\vx^k)-f^* - \frac{\alpha_k}{2} \|\nabla f(\vx^k)\|^2.
    \]
    Setting $a_k = f(\vx^k) - f^*$ and applying assumption~\ref{P} to the previous inequality, we have
    \[
    a_{k+1} \leqslant a_k - {\alpha_k\mu} a_k^{2\theta}.
    \]
     If $\theta = \frac{1}{2}$, then for any $k$ where $a_k > 0$, this inequality reduces to $a_{k+1}\leqslant (1-\alpha_k\mu) a_k$. Since $a_k = 0$ implies $a_i = 0$ for all $i \geqslant k$, we conclude that
    \[
    a_k\leqslant {\prod}_{i=0}^{k-1} (1-\alpha_i\mu) a_0 \leqslant \exp(-\mu{\sum}_{i=0}^{k-1} \alpha_i).
    \]
    If $\theta \in (\frac{1}{2},1]$, then for any $k$ with $a_k>0$, using the convexity of the function $t\mapsto t^{1-2\theta}$ on $(0, \infty)$, it follows that
    $a_{k+1}^{1-2\theta}-a_k^{1-2\theta} \geqslant (1-2\theta) a_k^{-2\theta} (a_{k+1}-a_k) \geqslant (2\theta-1)\alpha_k\mu$, where the second inequality is due to $a_{k+1} - a_k \leqslant -\alpha_k\mu a_k^{2\theta}$ and $1-2\theta < 0$. Summing this relation over $i=0, \dots, k-1$ yields $a_k \leqslant [a_0^{1-2\theta}+(2\theta-1)\mu{\sum}_{i=0}^{k-1} \alpha_i]^{-\frac{1}{2\theta-1}}$. \hfill\BlackBox \vspace{0.5ex}

Finally, we provide details about the step size estimates shown in \Cref{sec:noise-adaptivity}. According to Theorems~\ref{thm:sgd_const} and \ref{thm:sgd_poly}, the optimally tuned choices for the constant and polynomial step sizes are given by $\alpha_k = \alpha = \Theta([\frac{\log(K)}{K}]^{\varrho_s})$ and $\alpha_k = \Theta(k^{-\varrho_s})$. Applying \Cref{prop:int}, this yields ${\sum}_{k=0}^{K-1} \alpha_k = \Theta([\log(K)]^{\varrho_s}K^{1-\varrho_s})$ and ${\sum}_{k=0}^{K-1} \alpha_k = \Theta(K^{1-\varrho_s})$, respectively. Similarly, based on \Cref{thm:sgd_cos}, the optimal cosine step sizes satisfy $\alpha_k = \alpha [({1+\cos(\tfrac{k\pi}{K})})/{2}]^p$, $\alpha = [\tfrac{\beta\log(K)}{K}]^{\varrho_s}$, and using \Cref{lem:Tech-ineqs} (f), we have ${\sum}_{k=0}^{K-1} \alpha_k = \alpha {\sum}_{k=0}^{K-1} [({1+\cos(\tfrac{k\pi}{K})})/{2}]^p = \Theta([\log(K)]^{\varrho_s}K^{1-\varrho_s})$. For exponential step sizes, we can choose $\alpha_k = \alpha \gamma^k$, $\gamma = (\tfrac{\beta}{K})^{\frac{p}{K}}$ and following the proof of \Cref{lem:exp&cos0}, we obtain ${\sum}_{k=0}^{K-1} \alpha_k = \alpha \tfrac{1-\gamma^K}{1-\gamma} \geq \tfrac{\alpha [1-(\beta/K)^p]}{p}\tfrac{K}{\log(K/\beta)}$.

\vskip 0.2in
\bibliography{reference}

\begin{thebibliography}{50}
\providecommand{\natexlab}[1]{#1}
\providecommand{\url}[1]{\texttt{#1}}
\expandafter\ifx\csname urlstyle\endcsname\relax
  \providecommand{\doi}[1]{doi: #1}\else
  \providecommand{\doi}{doi: \begingroup \urlstyle{rm}\Url}\fi

\bibitem[Abadi et~al.(2015)Abadi, Agarwal, Barham, Brevdo, Chen, Citro, Corrado, Davis, et~al.]{tensorflow2015}
M.~Abadi, A.~Agarwal, P.~Barham, E.~Brevdo, Z.~Chen, C.~Citro, G.~S. Corrado, A.~Davis, et~al.
\newblock {TensorFlow}: Large-scale machine learning on heterogeneous systems, 2015.

\bibitem[Agarwal et~al.(2021)Agarwal, Kakade, Lee, and Mahajan]{agarwal2021theory}
A.~Agarwal, S.~M. Kakade, J.~D. Lee, and G.~Mahajan.
\newblock On the theory of policy gradient methods: Optimality, approximation, and distribution shift.
\newblock \emph{J. Mach. Learn. Res.}, 22\penalty0 (1):\penalty0 4431--4506, 2021.

\bibitem[Ahn et~al.(2020)Ahn, Yun, and Sra]{ahn2020sgd}
K.~Ahn, C.~Yun, and S.~Sra.
\newblock {SGD} with shuffling: {O}ptimal rates without component convexity and large epoch requirements.
\newblock In \emph{Advances in Neural Information Processing Systems}, volume~33, pages 17526--17535, 2020.

\bibitem[Attia and Koren(2024)]{attia2024free}
A.~Attia and T.~Koren.
\newblock How free is parameter-free stochastic optimization?
\newblock \emph{arXiv preprint arXiv:2402.03126}, 2024.

\bibitem[Bolte et~al.(2014)Bolte, Sabach, and Teboulle]{BolSabTeb14}
J.~Bolte, S.~Sabach, and M.~Teboulle.
\newblock Proximal alternating linearized minimization for nonconvex and nonsmooth problems.
\newblock \emph{Math. Program.}, 146\penalty0 (1-2, Ser. A):\penalty0 459--494, 2014.

\bibitem[Bottou et~al.(2018)Bottou, Curtis, and Nocedal]{bottou2018optimization}
L.~Bottou, F.~E. Curtis, and J.~Nocedal.
\newblock Optimization methods for large-scale machine learning.
\newblock \emph{SIAM Rev.}, 60\penalty0 (2):\penalty0 223--311, 2018.

\bibitem[Chung(1954)]{Chu54}
K.~L. Chung.
\newblock On a stochastic approximation method.
\newblock \emph{Ann. Math. Statist.}, pages 463--483, 1954.

\bibitem[Derman(1956)]{derman1956application}
C.~Derman.
\newblock An application of {C}hung's lemma to the {K}iefer-{W}olfowitz stochastic approximation procedure.
\newblock \emph{Ann. Math. Statist.}, pages 532--536, 1956.

\bibitem[Duchi et~al.(2011)Duchi, Hazan, and Singer]{duchi2011adaptive}
J.~Duchi, E.~Hazan, and Y.~Singer.
\newblock Adaptive subgradient methods for online learning and stochastic optimization.
\newblock \emph{J. Mach. Learn. Res.}, 12\penalty0 (7):\penalty0 2121--2159, 2011.

\bibitem[Fabian(1968)]{fabian1968asymptotic}
V.~Fabian.
\newblock On asymptotic normality in stochastic approximation.
\newblock \emph{Ann. Math. Statist.}, pages 1327--1332, 1968.

\bibitem[Fatkhullin et~al.(2022)Fatkhullin, Etesami, He, and Kiyavash]{Fatkhullin2022}
I.~Fatkhullin, J.~Etesami, N.~He, and N.~Kiyavash.
\newblock Sharp analysis of stochastic optimization under global {K}urdyka-{{\L}}ojasiewicz inequality.
\newblock In \emph{Advances in Neural Information Processing Systems}, volume~35, pages 15836--15848, 2022.

\bibitem[Fontaine et~al.(2021)Fontaine, De~Bortoli, and Durmus]{fontaine2021convergence}
X.~Fontaine, V.~De~Bortoli, and A.~Durmus.
\newblock Convergence rates and approximation results for {SGD} and its continuous-time counterpart.
\newblock In \emph{Conference on Learning Theory}, pages 1965--2058. PMLR, 2021.

\bibitem[Ge et~al.(2019)Ge, Kakade, Kidambi, and Netrapalli]{ge2019step}
R.~Ge, S.~M. Kakade, R.~Kidambi, and P.~Netrapalli.
\newblock The step decay schedule: A near optimal, geometrically decaying learning rate procedure for least squares.
\newblock In \emph{Advances in Neural Information Processing Systems}, volume~32, pages 14977--14988, 2019.

\bibitem[Goffin(1977)]{goffin1977convergence}
J.-L. Goffin.
\newblock On convergence rates of subgradient optimization methods.
\newblock \emph{Math. Program.}, 13:\penalty0 329--347, 1977.

\bibitem[Gower et~al.(2021)Gower, Sebbouh, and Loizou]{gower2021sgd}
R.~Gower, O.~Sebbouh, and N.~Loizou.
\newblock {SGD} for structured nonconvex functions: Learning rates, minibatching and interpolation.
\newblock In \emph{International Conference on Artificial Intelligence and Statistics}, pages 1315--1323. PMLR, 2021.

\bibitem[Gower et~al.(2019)Gower, Loizou, Qian, Sailanbayev, Shulgin, and Richt{\'a}rik]{gower2019sgd}
R.~M. Gower, N.~Loizou, X.~Qian, A.~Sailanbayev, E.~Shulgin, and P.~Richt{\'a}rik.
\newblock {SGD}: General analysis and improved rates.
\newblock In \emph{International Conference on Machine Learning}, pages 5200--5209. PMLR, 2019.

\bibitem[G{\"u}rb{\"u}zbalaban et~al.(2019)G{\"u}rb{\"u}zbalaban, Ozdaglar, and Parrilo]{gurbu2015}
M.~G{\"u}rb{\"u}zbalaban, A.~Ozdaglar, and P.~A. Parrilo.
\newblock Convergence rate of incremental gradient and incremental {N}ewton methods.
\newblock \emph{SIAM J. Optim.}, 29\penalty0 (4):\penalty0 2542--2565, 2019.

\bibitem[G{\"u}rb{\"u}zbalaban et~al.(2021)G{\"u}rb{\"u}zbalaban, Ozdaglar, and Parrilo]{gurbu2019}
M.~G{\"u}rb{\"u}zbalaban, A.~Ozdaglar, and P.~Parrilo.
\newblock Why random reshuffling beats stochastic gradient descent.
\newblock \emph{Math. Program.}, 186\penalty0 (1-2):\penalty0 49--84, 2021.

\bibitem[{HaoChen} and Sra(2019)]{haochen2019}
J.~Z. {HaoChen} and S.~Sra.
\newblock Random shuffling beats {SGD} after finite epochs.
\newblock In \emph{International Conference on Machine Learning}, pages 2624--2633. PMLR, 2019.

\bibitem[Huang et~al.(2021)Huang, Yuan, Mao, and Yin]{huang2021improved}
X.~Huang, K.~Yuan, X.~Mao, and W.~Yin.
\newblock An improved analysis and rates for variance reduction under without-replacement sampling orders.
\newblock In \emph{Advances in Neural Information Processing Systems}, volume~34, pages 3232--3243, 2021.

\bibitem[Karimi et~al.(2016)Karimi, Nutini, and Schmidt]{KarNutSch16}
H.~Karimi, J.~Nutini, and M.~Schmidt.
\newblock Linear convergence of gradient and proximal-gradient methods under the {P}olyak-{{\L}}ojasiewicz condition.
\newblock In \emph{Machine Learning and Knowledge Discovery in Databases: European Conference}, pages 795--811. Springer, 2016.

\bibitem[Khaled and Richt{\'a}rik(2020)]{khaled2020better}
A.~Khaled and P.~Richt{\'a}rik.
\newblock Better theory for {SGD} in the nonconvex world.
\newblock \emph{arXiv preprint arXiv:2002.03329}, 2020.

\bibitem[Kingma(2014)]{kingma2014adam}
D.~P. Kingma.
\newblock Adam: A method for stochastic optimization.
\newblock \emph{arXiv preprint arXiv:1412.6980}, 2014.

\bibitem[Li et~al.(2021)Li, Zhuang, and Orabona]{li2021second}
X.~Li, Z.~Zhuang, and F.~Orabona.
\newblock A second look at exponential and cosine step sizes: Simplicity, adaptivity, and performance.
\newblock In \emph{International Conference on Machine Learning}, pages 6553--6564. PMLR, 2021.

\bibitem[Li et~al.(2023)Li, Milzarek, and Qiu]{li2021convergence}
X.~Li, A.~Milzarek, and J.~Qiu.
\newblock Convergence of random reshuffling under the {K}urdyka--{{\L}}ojasiewicz inequality.
\newblock \emph{SIAM J. Optim.}, 33\penalty0 (2):\penalty0 1092--1120, 2023.

\bibitem[Loizou et~al.(2021)Loizou, Vaswani, Laradji, and Lacoste-Julien]{loizou2021stochastic}
N.~Loizou, S.~Vaswani, I.~H. Laradji, and S.~Lacoste-Julien.
\newblock Stochastic {P}olyak step-size for {SGD}: An adaptive learning rate for fast convergence.
\newblock In \emph{International Conference on Artificial Intelligence and Statistics}, pages 1306--1314. PMLR, 2021.

\bibitem[{\L}ojasiewicz(1959)]{lojasiewicz1959}
S.~{\L}ojasiewicz.
\newblock Sur le probl{\`e}me de la division.
\newblock \emph{Studia Mathematica}, 18:\penalty0 87--136, 1959.

\bibitem[{\L}ojasiewicz(1963)]{lojasiewicz1963}
S.~{\L}ojasiewicz.
\newblock Une propri{\'e}t{\'e} topologique des sous-ensembles analytiques r{\'e}els.
\newblock \emph{Les {\'e}quations aux d{\'e}riv{\'e}es partielles}, 117:\penalty0 87--89, 1963.

\bibitem[Loshchilov and Hutter(2017)]{loshchilov2016sgdr}
I.~Loshchilov and F.~Hutter.
\newblock {SGDR}: Stochastic gradient descent with warm restarts.
\newblock In \emph{International Conference on Learning Representations}, 2017.

\bibitem[Mei et~al.(2020)Mei, Xiao, Szepesvari, and Schuurmans]{mei2020global}
J.~Mei, C.~Xiao, C.~Szepesvari, and D.~Schuurmans.
\newblock On the global convergence rates of softmax policy gradient methods.
\newblock In \emph{International Conference on Machine Learning}, pages 6820--6829. PMLR, 2020.

\bibitem[Mishchenko et~al.(2020)Mishchenko, Khaled, and Richtarik]{Mis2020}
K.~Mishchenko, A.~Khaled, and P.~Richtarik.
\newblock Random reshuffling: Simple analysis with vast improvements.
\newblock In \emph{Advances in Neural Information Processing Systems}, volume~33, pages 17309--17320, 2020.

\bibitem[Moulines and Bach(2011)]{Moulines2011}
E.~Moulines and F.~Bach.
\newblock Non-asymptotic analysis of stochastic approximation algorithms for machine learning.
\newblock In \emph{Advances in Neural Information Processing Systems}, volume~24, pages 451--459, 2011.

\bibitem[Needell et~al.(2014)Needell, Ward, and Srebro]{needell2014stochastic}
D.~Needell, R.~Ward, and N.~Srebro.
\newblock Stochastic gradient descent, weighted sampling, and the randomized {K}aczmarz algorithm.
\newblock In \emph{Advances in Neural Information Processing Systems}, volume~27, pages 1017--1025, 2014.

\bibitem[Nemirovski et~al.(2009)Nemirovski, Juditsky, Lan, and Shapiro]{NemJudLanSha09}
A.~Nemirovski, A.~Juditsky, G.~Lan, and A.~Shapiro.
\newblock Robust stochastic approximation approach to stochastic programming.
\newblock \emph{SIAM J. Optim.}, 19\penalty0 (4):\penalty0 1574--1609, 2009.

\bibitem[Nguyen et~al.(2021)Nguyen, Tran-Dinh, Phan, Nguyen, and van Dijk]{nguyen2020unified}
L.~M. Nguyen, Q.~Tran-Dinh, D.~T. Phan, P.~H. Nguyen, and M.~van Dijk.
\newblock A unified convergence analysis for shuffling-type gradient methods.
\newblock \emph{J. Mach. Learn. Res.}, 22:\penalty0 1--44, 2021.

\bibitem[Orabona and P{\'a}l(2021)]{orabona2021parameter}
F.~Orabona and D.~P{\'a}l.
\newblock Parameter-free stochastic optimization of variationally coherent functions.
\newblock \emph{arXiv preprint arXiv:2102.00236}, 2021.

\bibitem[Paszke et~al.(2019)Paszke, Gross, Massa, Lerer, Bradbury, Chanan, Killeen, Lin, Gimelshein, Antiga, et~al.]{paszke2019pytorch}
A.~Paszke, S.~Gross, F.~Massa, A.~Lerer, J.~Bradbury, G.~Chanan, T.~Killeen, Z.~Lin, N.~Gimelshein, L.~Antiga, et~al.
\newblock Pytorch: An imperative style, high-performance deep learning library.
\newblock In \emph{Advances in Neural Information Processing Systems}, volume~32, pages 8026--8037, 2019.

\bibitem[Polyak(1963)]{polyak1963}
B.~T. Polyak.
\newblock Gradient methods for the minimisation of functionals.
\newblock \emph{USSR Computational Mathematics and Mathematical Physics}, 3\penalty0 (4):\penalty0 864--878, 1963.

\bibitem[Rakhlin et~al.(2012)Rakhlin, Shamir, and Sridharan]{rakhlin2012making}
A.~Rakhlin, O.~Shamir, and K.~Sridharan.
\newblock Making gradient descent optimal for strongly convex stochastic optimization.
\newblock In \emph{International Conference on Machine Learning}, pages 449--456. PMLR, 2012.

\bibitem[Robbins and Monro(1951)]{RobMon51}
H.~Robbins and S.~Monro.
\newblock A stochastic approximation method.
\newblock \emph{Ann. Math. Statist.}, pages 400--407, 1951.

\bibitem[Sacks(1958)]{Sac58}
J.~Sacks.
\newblock Asymptotic distribution of stochastic approximation procedures.
\newblock \emph{Ann. Math. Statist.}, 29\penalty0 (2):\penalty0 373--405, 1958.

\bibitem[Shamir(2016)]{shamir2016}
O.~Shamir.
\newblock Without-replacement sampling for stochastic gradient methods.
\newblock In \emph{Advances in Neural Information Processing Systems}, volume~29, pages 46--54, 2016.

\bibitem[Shamir and Zhang(2013)]{ShaZha2013}
O.~Shamir and T.~Zhang.
\newblock Stochastic gradient descent for non-smooth optimization: Convergence results and optimal averaging schemes.
\newblock In \emph{International Conference on Machine Learning}, pages 71--79. PMLR, 2013.

\bibitem[Sohl-Dickstein et~al.(2014)Sohl-Dickstein, Poole, and Ganguli]{sohl2014fast}
J.~Sohl-Dickstein, B.~Poole, and S.~Ganguli.
\newblock Fast large-scale optimization by unifying stochastic gradient and quasi-{N}ewton methods.
\newblock In \emph{International Conference on Machine Learning}, pages 604--612. PMLR, 2014.

\bibitem[Touvron et~al.(2023)Touvron, Martin, Stone, Albert, Almahairi, Babaei, Bashlykov, Batra, Bhargava, Bhosale, et~al.]{touvron2023llama}
H.~Touvron, L.~Martin, K.~Stone, P.~Albert, A.~Almahairi, Y.~Babaei, N.~Bashlykov, S.~Batra, P.~Bhargava, S.~Bhosale, et~al.
\newblock Llama 2: Open foundation and fine-tuned chat models.
\newblock \emph{arXiv preprint arXiv:2307.09288}, 2023.

\bibitem[Tran et~al.(2021)Tran, Nguyen, and Tran-Dinh]{tran2021smg}
T.~H. Tran, L.~M. Nguyen, and Q.~Tran-Dinh.
\newblock {SMG}: A shuffling gradient-based method with momentum.
\newblock In \emph{International Conference on Machine Learning}, pages 10379--10389. PMLR, 2021.

\bibitem[Vaswani et~al.(2022)Vaswani, Dubois-Taine, and Babanezhad]{vaswani2022towards}
S.~Vaswani, B.~Dubois-Taine, and R.~Babanezhad.
\newblock Towards noise-adaptive, problem-adaptive (accelerated) stochastic gradient descent.
\newblock In \emph{International Conference on Machine Learning}, pages 22015--22059. PMLR, 2022.

\bibitem[Wang et~al.(2021)Wang, Magn{\'u}sson, and Johansson]{wang2021convergence}
X.~Wang, S.~Magn{\'u}sson, and M.~Johansson.
\newblock On the convergence of step decay step-size for stochastic optimization.
\newblock In \emph{Advances in Neural Information Processing Systems}, volume~34, pages 14226--14238, 2021.

\bibitem[Wang and Zou(2022)]{wang2022policy}
Y.~Wang and S.~Zou.
\newblock Policy gradient method for robust reinforcement learning.
\newblock In \emph{International Conference on Machine Learning}, pages 23484--23526. PMLR, 2022.

\bibitem[Xiao(2022)]{xiao2020policy}
L.~Xiao.
\newblock On the convergence rates of policy gradient methods.
\newblock \emph{J. Mach. Learn. Res.}, 23:\penalty0 1--36, 2022.

\end{thebibliography}

\end{document}